\documentclass[12pt]{article}
\PassOptionsToPackage{hyphens}{url}
\usepackage[pagebackref=true,colorlinks]{hyperref}
\usepackage{amsmath,amssymb,amsthm,mathtools,amsrefs}
\usepackage{csquotes}
\usepackage{etex} 
\usepackage{tikz}
\usetikzlibrary{calc}
\usepackage{pgfplots}
\pgfplotsset{compat=newest}
\usepackage{adjustbox}
\usepackage{stmaryrd}
\usepackage{fancyhdr}
\usepackage{soul}
\usepackage{dsfont}
\usepackage[normalem]{ulem}
\usepackage{longtable}
\usepackage[bottom]{footmisc}
\usepackage[margin=1cm]{caption}
\usepackage{subfig} 
\usepackage{pictexwd, dcpic}
\usepackage{float}
\usepackage{enumerate}
\usepackage[all]{xy} 
\hypersetup{
	colorlinks=true,
    linktoc=all,
    linkcolor=blue,  
    citecolor=blue,
    }

\makeatletter
\renewcommand*\env@matrix[1][\arraystretch]{%
  \edef\arraystretch{#1}%
  \hskip -\arraycolsep
  \let\@ifnextchar\new@ifnextchar
  \array{*\c@MaxMatrixCols c}}
\makeatother

\makeatletter
\pgfmathdeclarefunction{erf}{1}{%
  \begingroup
    \pgfmathparse{#1 > 0 ? 1 : -1}%
    \edef\sign{\pgfmathresult}%
    \pgfmathparse{abs(#1)}%
    \edef\x{\pgfmathresult}%
    \pgfmathparse{1/(1+0.3275911*\x)}%
    \edef\t{\pgfmathresult}%
    \pgfmathparse{%
      1 - (((((1.061405429*\t -1.453152027)*\t) + 1.421413741)*\t 
      -0.284496736)*\t + 0.254829592)*\t*exp(-(\x*\x))}%
    \edef\y{\pgfmathresult}%
    \pgfmathparse{(\sign)*\y}%
    \pgfmath@smuggleone\pgfmathresult%
  \endgroup
}
\makeatother

\theoremstyle{theorem}
\newtheorem{theorem}[equation]{Theorem}
\newtheorem{lemma}[equation]{Lemma}
\newtheorem{proposition}[equation]{Proposition}
\newtheorem{corollary}[equation]{Corollary}
\theoremstyle{definition}
\newtheorem{definition}[equation]{Definition}
\newtheorem{construction}[equation]{Construction}
\newtheorem{question}[equation]{Question}
\newtheorem{problem}[equation]{Problem}
\newtheorem{example}[equation]{Example}
\newtheorem{exercise}[equation]{Exercise}
\newtheorem*{answer}{Answer}
\newtheorem*{solution}{Solution}
\newtheorem{remark}[equation]{Remark}

\newtheorem{notation}[equation]{Notation}

\numberwithin{equation}{section}

\let\a=\alpha \let\b=\beta  \let\de=\delta \let\e=\epsilon
\let\z=\zeta  \let\q=\theta  \let\k=\kappa
\let\l=\lambda \let\r=\rho
\let\s=\sigma \let\t=\tau    
\let\w=\omega      \let\G=\Gamma  \let\Q=\Theta 
     
\let\C=\Chi \let\W=\Omega
\def\vf{\varphi}

\newcommand{\be}{\begin{equation}}
\newcommand{\ee}{\end{equation}}
\def\ba{\begin{align}} 
\def\ea{\end{align}}
\newcommand{\bea}{\begin{eqnarray}}
\newcommand{\eea}{\end{eqnarray}}
\newcommand{\bx}{\begin{example}}
\newcommand{\ex}{\end{example}}
\newcommand{\bex}{\begin{exercise}}
\newcommand{\eex}{\end{exercise}}
\newcommand{\ban}{\begin{answer}}
\newcommand{\ean}{\end{answer}}
\newcommand{\bt}{\begin{theorem}}
\newcommand{\et}{\end{theorem}}
\newcommand{\bc}{\begin{corollary}}
\newcommand{\ec}{\end{corollary}}
\newcommand{\blem}{\begin{lemma}}
\newcommand{\elem}{\end{lemma}}
\newcommand{\bp}{\begin{problem}}
\newcommand{\ep}{\end{problem}}
\newcommand{\bn}{\begin{proposition}}
\newcommand{\en}{\end{proposition}}
\newcommand{\bd}{\begin{definition}}
\newcommand{\ed}{\end{definition}}
\newcommand{\bcon}{\begin{construction}}
\newcommand{\econ}{\end{construction}}
\newcommand{\bq}{\begin{question}}
\newcommand{\eq}{\end{question}}
\newcommand{\bprf}{\begin{proof}}
\newcommand{\eprf}{\end{proof}}
\newcommand{\br}{\begin{remark}}
\newcommand{\er}{\end{remark}}
\newcommand{\bs}{\begin{solution}}
\newcommand{\es}{\end{solution}}
\newcommand{\beqs}{\begin{eqnarray}}
\newcommand{\eeqs}{\end{eqnarray}}

 \let\ov=\overline

\newcommand{\<}{\langle}
\renewcommand{\>}{\rangle}

\newcommand{\id}{\mathrm{id}}

\newcommand{\mC}{\mathcal{C}}
\newcommand{\mD}{\mathcal{D}}

\newcommand{\up}{\uparrow}
\newcommand{\dn}{\downarrow}

\newcommand{\aand}{\qquad \& \qquad}
\newcommand{\tr}{{\rm tr} }

\def\ev{{{\mathrm{ev}}}}

\def\R{{{\mathbb R}}}
\def\C{{{\mathbb C}}}

\def\K{{{\mathbb K}}}

\def\N{{{\mathbb N}}}
\def\Z{{{\mathbb Z}}}
\def\Q{{{\mathbb Q}}}

\def\Hi{{{\mathcal{H}}}}
\newcommand{\CAlg}{\mathbf{C\text{*-}Alg}}
\newcommand{\CPCAlg}{\mathbf{C\text{*-}AlgCP}}
\newcommand{\CPfdCAlg}{\mathbf{fdC\text{*-}AlgCP}}
\newcommand{\cCAlg}{\mathbf{cC\text{*-}Alg}}
\newcommand{\fdcCAlg}{\mathbf{fdcC\text{*-}Alg}}
\newcommand{\fdCAlg}{\mathbf{fdC\text{*-}Alg}}
\newcommand{\cCAlgPos}{\mathbf{cC\text{*-}AlgPos}}
\newcommand{\fdcCAlgPos}{\mathbf{fdcC\text{*-}AlgPos}}

\def\mA{{{\mathcal{A}}}}
\def\mS{{{\mathcal{S}}}}
\def\mB{{{\mathcal{B}}}}

\newcommand{\FinSetFun}{\mathbf{FinSetFun}}
\newcommand{\FinSetStoch}{\mathbf{FinSetStoch}}

\newcommand{\op}{\mathrm{op}}

\newcommand{\FinProb}{\mathbf{FinProb}}

\newcommand{\cH}{\mathbf{cHaus}}
\newcommand{\cHStoch}{\mathbf{cHausStoch}}

\newcommand{\stoch}{\;\xy0;/r.25pc/:\ar@{~>}(3,0);(-3,0);\endxy\;}

\newcommand{\ds}{\displaystyle}
\newcommand{\ben}{\renewcommand{\theenumi}{\alph{enumi}} 
\renewcommand{\labelenumi}{(\theenumi)}\begin{enumerate}}
\newcommand{\een}{\end{enumerate}}

\newcommand{\dieone}{
\begin{tikzpicture}[scale=0.4,baseline=-0.6ex]
\draw[rounded corners=2] (-0.5,-0.5) rectangle (0.5,0.5);
\node[scale=0.5] at (0,0) {$\bullet$};
\end{tikzpicture}
}
\newcommand{\diethree}{
\begin{tikzpicture}[scale=0.4,baseline=-0.6ex]
\draw[rounded corners=2] (-0.5,-0.5) rectangle (0.5,0.5);
\node[scale=0.5] at (-0.25,-0.25) {$\bullet$};
\node[scale=0.5] at (0,0) {$\bullet$};
\node[scale=0.5] at (0.25,0.25) {$\bullet$};
\end{tikzpicture}
}
\newcommand{\diefive}{
\begin{tikzpicture}[scale=0.4,baseline=-0.6ex]
\draw[rounded corners=2] (-0.5,-0.5) rectangle (0.5,0.5);
\node[scale=0.5] at (-0.25,-0.25) {$\bullet$};
\node[scale=0.5] at (-0.25,0.25) {$\bullet$};
\node[scale=0.5] at (0,0) {$\bullet$};
\node[scale=0.5] at (0.25,-0.25) {$\bullet$};
\node[scale=0.5] at (0.25,0.25) {$\bullet$};
\end{tikzpicture}
}
\newcommand{\dietwo}{
\begin{tikzpicture}[scale=0.4,baseline=-0.6ex]
\draw[rounded corners=2] (-0.5,-0.5) rectangle (0.5,0.5);
\node[scale=0.5] at (-0.25,0.25) {$\bullet$};
\node[scale=0.5] at (0.25,-0.25) {$\bullet$};
\end{tikzpicture}
}
\newcommand{\diefour}{
\begin{tikzpicture}[scale=0.4,baseline=-0.6ex]
\draw[rounded corners=2] (-0.5,-0.5) rectangle (0.5,0.5);
\node[scale=0.5] at (-0.25,-0.25) {$\bullet$};
\node[scale=0.5] at (-0.25,0.25) {$\bullet$};
\node[scale=0.5] at (0.25,-0.25) {$\bullet$};
\node[scale=0.5] at (0.25,0.25) {$\bullet$};
\end{tikzpicture}
}
\newcommand{\diesix}{
\begin{tikzpicture}[scale=0.4,baseline=-0.6ex]
\draw[rounded corners=2] (-0.5,-0.5) rectangle (0.5,0.5);
\node[scale=0.5] at (-0.25,-0.25) {$\bullet$};
\node[scale=0.5] at (-0.25,0.25) {$\bullet$};
\node[scale=0.5] at (-0.25,0) {$\bullet$};
\node[scale=0.5] at (0.25,0) {$\bullet$};
\node[scale=0.5] at (0.25,-0.25) {$\bullet$};
\node[scale=0.5] at (0.25,0.25) {$\bullet$};
\end{tikzpicture}
}

\newcommand\blfootnote[1]{%
  \begingroup
  \renewcommand\thefootnote{}\footnote{#1}%
  \addtocounter{footnote}{-1}%
  \endgroup
}

\title{Discrete probabilistic and algebraic dynamics:\\
a stochastic commutative Gelfand-Naimark Theorem}
\author{Arthur J. Parzygnat}
\date{}

\begin{document}
\maketitle 

\begin{abstract}
We introduce a category of stochastic maps 
(certain Markov kernels) on 
compact Hausdorff spaces, 
construct a stochastic analogue of the
Gelfand spectrum functor, 
and prove a stochastic version of the commutative Gelfand-Naimark Theorem.
This relates concepts from algebra and operator theory 
to concepts from topology and probability theory. 
For completeness, we review stochastic matrices, 
their relationship to positive maps on commutative $C^*$-algebras, 
and the Gelfand-Naimark Theorem. 
No knowledge of probability theory 
nor $C^*$-algebras is assumed 
and several examples are drawn from physics.
\end{abstract}

\tableofcontents

\section{An algebraic perspective on probability theory}
\label{sec:fancywords}

\begin{displayquote}
\emph{That's another thing about categories, not that they give you the proofs, but they tell you what to prove!}---Mikhail Gromov%
\footnote{This quote is taken from Mikhail Gromov's 
Lecture 1 on ``Probability, symmetry, linearity''
given at the Institut des Hautes \'Etudes Scientifiques on October 10, 2014.}
\end{displayquote}

\subsection{Brief background and motivation}

The fact that every compact Hausdorff space can be recovered from its
$C^*$-algebra of complex-valued continuous functions was discovered
by Gelfand and Naimark in 1943 (see Lemma 1 in \cite{GN43}). 
This involved constructing a space out of a commutative $C^*$-algebra and a
continuous map out of a $*$-homomorphism of algebras.
Extending earlier work of Riesz and Markov (and many others), 
Kakutani proved that finite regular measures on such spaces
correspond to positive linear functionals on the algebras
(see Theorem 9 in \cite{Ka41}). 
One can combine these ideas in the following way. 
Continuous maps between spaces can be viewed as deterministic
assignments in the sense that their values on inputs are specified precisely.
Instead, we can imagine a non-deterministic
analogue as a ``smeared out'' function, or more precisely, 
an assignment sending inputs to probability measures. 
Such ``stochastic maps'' correspond to (completely) positive linear
maps between the algebras. 
The latter are well-defined for non-commutative $C^*$-algebras and 
describe non-deterministic dynamics in quantum systems \cite{Li76}. 
Although this correspondence is 
well-understood to a large extent \cite{Ci06}, \cite{EFHN15} 
some important features are missed
without considering the \emph{category} of such spaces and stochastic maps. 
For example, how are the compositions of stochastic maps defined, 
are compositions respected under the correspondence, 
what is the explicit procedure to go back from a positive map to a stochastic 
map, and is there a way to formulate all of these dualities in the same context? 

Our aim is to introduce these ideas in a suitable language
and answer these questions. In the spirit
of the Rosetta Stone \cite{BaSt11}, we will use category theory to do this, 
though we will only use the bare minimum necessary. 
In particular, we only assume the reader knows what 
categories, functors, and natural transformations are and how to 
compose them though the reader may also understand the results
without even knowing these definitions verbatim from memory. 
Using category theory has the
immediate benefit of formulating our questions precisely thereby 
making the goals clear: we must identify the appropriate 
topological, algebraic, deterministic, and non-deterministic categories, construct
mappings between them, and prove the above mentioned equivalences.

This work fits into the general context of categorical probability theory, which 
was pioneered by Lawvere who first defined a category of stochastic maps in 
1962 \cite{La62}. 
Further work in the area of category-theoretic aspects of measure 
theory was developed by Linton in 1963 \cite{Li63}
and the Polish school in the early 1970's, 
notably Semadeni \cite{Se73} and \'Swirszcz \cite{Sw74}.
In the early 1980's, Giry further developed these ideas by describing the monadic
structures available in probability theory on a larger class of spaces \cite{Gi82}. 
More recently, this area was revisited in the area of computer science 
where functional programming makes use of category-theoretic ideas,
particularly monads. The relationship to the category of ordinary relations,
used in the theory of computing, was generalized to the stochastic
setting by Panangaden in 1998 \cite{Pa99} (an earlier version of this
article was titled ``Probabilistic Relations'') based on earlier work of Kozen
in 1981 \cite{Ko81}.
A monadic and algebraic viewpoint on the questions raised above were worked
out by Furber and Jacobs recently \cite{FuJa13}. Our goals are similar, 
but we approach the problem from an analytic perspective. 

\subsection{Overview of results}

We prove a generalization of the commutative Gelfand-Naimark Theorem
(henceforth referred to as just the Gelfand-Naimark Theorem)
that is valid for continuous stochastic maps on compact Hausdorff spaces. 
This is done at three different levels of increasing complication, 
the first two of which are well known.  
In Section \ref{sec:stochasticmatrices}, 
we review the equivalence between stochastic matrices
on finite sets and positive maps on finite-dimensional commutative (unital)
$C^*$-algebras, which are all basically of the form $\C^{n}$ for some $n\in\N$
with pointwise algebraic structure. Specializing to the case of functions instead of 
stochastic matrices corresponds to $*$-homomorphisms on the
$C^*$-algebra side. In Section \ref{sec:GelfandNaimark}, 
we review the Gelfand-Naimark
Theorem phrased in an appropriate categorical context as an equivalence
between the category of compact Hausdorff spaces with continuous maps
and the category of commutative unital $C^*$-algebras and unit-preserving
$*$-homomorphisms. We do not prove the fundamental theorems 
here since they are well documented in standard references, 
but we do introduce the necessary concepts needed to understand 
the statements in what follows. 

Finally, in Section \ref{sec:stochasticGelfandNaimark}, we extend the 
Gelfand-Naimark Theorem to an equivalence between 
compact Hausdorff spaces with continuous stochastic maps 
and commutative $C^*$-algebras with positive unital maps. 
To even state this theorem, an appropriate category of such
spaces and stochastic maps must be defined. Although versions 
of this category have been defined in the literature
(see \cite{Gi82} and \cite{Pa99} for example),
the spaces are typically either 
measurable spaces
or
Polish spaces.
There is an important distinction between these two cases. 
Briefly, a \uline{\emph{Markov kernel}} from a measurable space $(X,\mathcal{M}_{X})$
to another one $(Y,\mathcal{M}_{Y})$ is a function 
$X\times\mathcal{M}_{Y}\to[0,1]$ such that fixing the left variable 
gives a probability measure and fixing the right variable gives a measurable 
function.
A Markov kernel may equivalently be defined as a stochastic map, i.e. 
a measurable function $X\to\mathrm{ProbMeas}(Y),$ where 
$\mathrm{ProbMeas}(Y)$ is the set of all probability measures on $Y$
equipped with the smallest $\s$-algebra for which the evaluation function
$\mathrm{ev}_{E}:\mathrm{ProbMeas}(Y)\to[0,1],$
defined by 
$\mathrm{ProbMeas}(Y)\ni\nu\mapsto\mathrm{ev}_{E}(\nu):=\nu(E),$ 
is measurable for all $E\in\mathcal{M}_{Y}.$ 
In other words, in the case of measurable spaces, 
stochastic maps are in one-to-one correspondence with Markov kernels.
However, in the case of Polish spaces, many applications demand a 
more restrictive class of Markov kernels that have additional continuity 
properties, but demanding that the Markov kernel is continuous
is too restrictive in several contexts 
(see Remarks \ref{rmk:LSCversusC} and \ref{rmk:lackofcontinuity}).
Fortunately, 
one can use the second perspective and define stochastic maps
as continuous functions $X\to\mathrm{ProbMeas}(Y),$ 
where $\mathrm{ProbMeas}(Y)$ is now equipped with the vague topology.
As a result, it is a-priori unclear whether the associated kernel function
is Borel measurable after fixing the right variable. 
Nevertheless, it follows from the fact that the spaces are Polish \cite{Gi82}. 

If Polish spaces are replaced by compact Hausdorff spaces, 
one can still define stochastic maps, but the above proof of measurability 
of the associated kernels fails. To circumvent this difficulty, 
we restrict our probability measures to be regular (this additional
assumption is automatically satisfied for Polish spaces). 
We first prove that the stochastic maps 
on such spaces have measurable Markov kernels in Lemmas
\ref{lem:evaluationofmeasurecontinuous} and \ref{lem:evaluationmeasurable}. 
Lemma \ref{lem:evaluationofmeasurecontinuous} extends 
results due to Billingsley (see proof of Theorem 2.1 in \cite{Bi56}). 
This leads to an explicit definition of composition of such stochastic maps 
in Proposition \ref{prop:compositionstochastic} without referring
to the Riesz-Markov-Kakutani Representation theorem 
as is done in \cite{FuJa13} and \cite{EFHN15}. 
We then prove that stochastic maps 
on compact Hausdorff spaces
form a category, denoted by $\cHStoch$, in Theorem \ref{thm:cHStoch}. 

In addition, we provide
an explicit and geometric construction of the Gelfand spectrum functor from 
commutative $C^*$-algebras and positive maps to the category of 
compact Hausdorff spaces and continuous stochastic maps
without using the using algebras of functions and without 
appealing to the Riesz-Markov-Kakutani 
Representation Theorem.
Instead, this ``stochastic'' Gelfand spectrum functor
is constructed using techniques from convex analysis and Choquet theory
in Theorems \ref{thm:regularprobmeasonspectrumfromstate} and 
\ref{thm:stochasticspectrum}. 
Finally, we prove
that the stochastic Gelfand spectrum functor exhibits a categorical inverse 
to the continuous functions functor in Theorem 
\ref{thm:stochasticGN}. 

Our results can be succinctly summarized in that we complete the diagrammatic cube 
(the ${}^{\op}$ superscript refers to flipping the directionality of all morphisms)
\be
\label{eq:thecube}
\xy0;/r.35pc/:
(-15,20)*+{\cH^{\op}}="cHaus";
(-15,-5)*+{\color{blue}\cHStoch^{\op}}="cHausStoch";
(15,20)*+{\cCAlg}="cCAlg";
(15,-5)*+{\cCAlgPos}="cCAlgPos";
(-30,5)*+{\FinSetFun^{\op}}="FinSetFun";
(-30,-20)*+{\FinSetStoch^{\op}}="FinSetStoch";
(3,5)*+{\fdcCAlg}="fdcCAlg";
(3,-20)*+{\fdcCAlgPos}="fdcCAlgPos";
{\ar@<1.0ex>"cHaus";"cCAlg"};
{\ar@<1.0ex>"cCAlg";"cHaus"};
{\ar@<0.95ex>|(0.6){\hole}@[blue]"cHausStoch";"cCAlgPos"};
{\ar@<1.0ex>|(0.6){\hole}@[blue]"cHausStoch";"cCAlgPos"};
{\ar@<1.05ex>|(0.6){\hole}@[blue]"cHausStoch";"cCAlgPos"};
{\ar@<0.95ex>|(0.4){\hole}@[blue]"cCAlgPos";"cHausStoch"};
{\ar@<1.0ex>|(0.4){\hole}@[blue]"cCAlgPos";"cHausStoch"};
{\ar@<1.05ex>|(0.4){\hole}@[blue]"cCAlgPos";"cHausStoch"};
{\ar@<1.0ex>"FinSetFun";"fdcCAlg"};
{\ar@<1.0ex>"fdcCAlg";"FinSetFun"};
{\ar@<1.0ex>"FinSetStoch";"fdcCAlgPos"};
{\ar@<1.0ex>"fdcCAlgPos";"FinSetStoch"};
{\ar@{^{(}->}"FinSetFun";"cHaus"};
{\ar@<-0.05ex>@{^{(}->}@[blue]"FinSetStoch";"cHausStoch"};
{\ar@{^{(}->}@[blue]"FinSetStoch";"cHausStoch"};
{\ar@<0.05ex>@{^{(}->}@[blue]"FinSetStoch";"cHausStoch"};
{\ar@{^{(}->}"fdcCAlg";"cCAlg"};
{\ar@{^{(}->}"fdcCAlgPos";"cCAlgPos"};
{\ar"FinSetFun";"FinSetStoch"};
{\ar@<-0.05ex>|(0.55){\hole}|(0.65){\hole}@[blue]"cHaus";"cHausStoch"};
{\ar|(0.55){\hole}|(0.65){\hole}@[blue]"cHaus";"cHausStoch"};
{\ar@<0.05ex>|(0.55){\hole}|(0.65){\hole}@[blue]"cHaus";"cHausStoch"};
{\ar"fdcCAlg";"fdcCAlgPos"};
{\ar"cCAlg";"cCAlgPos"};
\endxy
\ee
by adding $\cHStoch$ and the functors into and from it. 
The arrows between the left and right faces describe the equivalences
between topology and analysis on the left and algebra on the right. 
The arrows from top to bottom describe the generalization of going from
deterministic dynamics to stochastic dynamics. 
The arrows from front to back describe the inclusion of finite systems
to possibly infinite systems. 
The arrows from the left to the right are particularly important in 
physics for the following reason. 
Much of classical and quantum dynamics, whether
deterministic or non-deterministic, can be formulated in a single
category: $C^*$-algebras and completely positive maps.
Although, classical evolution is typically thought of as being quite distinct from quantum
evolution, they both
take place in the \emph{same} category allowing a precise comparison
between dynamics. In particular, it allows mixing classical and quantum systems.
For example, measurement of quantum systems by 
macroscopic (classical) beings can be thought of as a map of $C^*$-algebras
and pulling back states from a non-commutative $C^*$-algebra to a 
commutative one (see Example \ref{ex:measurement}). 
The importance of commutative subalgebras
in the theory of measurement and their relationship to quantum 
mechanics has been emphasized
and summarized nicely in Heunen's work \cite{He17}. 

In summary, we first analyze the front and right face
in Section \ref{sec:stochasticmatrices}, 
then the top face in Section \ref{sec:GelfandNaimark}, 
and finally the 
category
$\cHStoch$ and how it fits into this cube via a stochastic 
Gelfand-Naimark Theorem
in Section \ref{sec:stochasticGelfandNaimark}. 
Concluding remarks, including a more detailed discussion
of the relationship between our work and others, are in Section \ref{sec:closing}. 
For the reader who has all the necessary background, our main results
are contained in Sections \ref{sec:cHausStoch} through 
\ref{sec:theclassicalcube}. The reader may find the 
\hyperref[indexofnotation]{Index of notation} at the end helpful. 

Probability theory evolved before the formulation of measure theory, 
and, in particular, before the birth of category theory. 
As such, it has developed its own language and culture, which
the author is largely unfamiliar with resulting in a 
potentially contentious presentation. Caveat lector. 

\vspace{3mm}
\noindent
\textbf{Acknowledgements.}

This work is based on a series of lectures given
at the Analysis Learning Seminar at the University of
Connecticut, Storrs and
the Mathematical Physics, Fourier Analysis and Applications Seminar
at the CUNY Graduate Center in New York
during the Spring 2017 semester. 
The author is grateful for the invitations
by Matthew Badger and Vyron Vellis to the former and 
by Azita Mayeli and Max Yarmolinsky to the latter. 
The author has greatly benefited from discussions with 
Stefan Andronache,
Iddo Ben-Ari,
Behrang Forghani, 
Marcelo Nomura,
Benjamin Russo,
Ambar Sengupta, 
Scott O. Wilson, 
and Yun Yang.
The author also thanks Markus Haase for several insightful comments
and suggestions on an earlier version of this work.
Parts of Sections \ref{sec:stochasticmatrices} and \ref{sec:GelfandNaimark}
were worked out when the author was partially supported by 
the CUNY Graduate Center Capelloni
Dissertation Fellowship and NSF grant PHY-1213380.

\section{Positive maps and stochastic matrices}
\label{sec:stochasticmatrices}

\subsection{Introduction}

We start off with basic probability theory on finite sets
and linear algebra. We first learned about this perspective
from the work of Baez and Fritz on entropy \cite{BaFr14}.
We find their notation, particularly distinguishing between
deterministic and non-deterministic processes, incredibly lucid. 
The former are depicted with straight arrows $\longrightarrow$ 
while the latter with curvy arrows $\stoch.$ 
We introduce the concepts of $C^*$-algebras
and states on $C^*$-algebras. There are two main categories described
here both of whose objects consist of $C^*$-algebras.
The difference occurs at the level of morphisms. 
In one category, the morphisms consist of algebra $*$-homomorphisms
and in the other case they consist of completely positive maps. 
We will show that there is a correspondence 
between these categories when restricted to finite-dimensional 
commutative $C^*$-algebras and ordinary concepts in probability theory, 
namely stochastic matrices. 
Thus, the category of operations in quantum mechanics can be viewed
as a non-commutative extension of the category of stochastic
(non-deterministic) processes in classical mechanics. 

In this section, we introduce
\begin{enumerate}[i.]
\item
$\FinSetFun$, the category of finite sets and functions (deterministic processes),
\item
$\FinSetStoch$, the category of finite sets and stochastic matrices 
(non-deterministic processes),
\item
$\fdcCAlg,$ the category of finite-dimensional commutative $\C^*$-algebras 
and $*$-homomorphisms 
(the non-commutative analogue of deterministic processes), and 
\item
$\fdcCAlgPos$, the category of finite-dimensional commutative 
$\C^*$-algebras and completely positive maps 
(the non-commutative analogue of stochastic maps)
\end{enumerate}
and show that there are functors (the 
$\op$ superscript will be explained in 
Theorem \ref{thm:FinProbtoCP})
\be
\xy0;/r.25pc/:
(-17.5,7.5)*+{\FinSetFun^{\op}}="FinProb";
(-17.5,-7.5)*+{\FinSetStoch^{\op}}="FinStoch";
(17.5,7.5)*+{\fdcCAlg}="states";
(17.5,-7.5)*+{\fdcCAlgPos}="CPstates";
{\ar@{^{(}->}@<-0.5ex>"FinProb";"states"};
{\ar@{^{(}->}@<-0.3ex>"FinStoch";"CPstates"};
{\ar"FinProb";"FinStoch"};
{\ar"states";"CPstates"};
\endxy
\ee
so that the horizontal functors are fully faithful and so that the diagram 
commutes. 

These concepts will be rigorously defined and sufficient intuition
will be provided to illustrate that they are nothing more than 
simple ideas in probability theory. 
However, we do assume the reader is (at least vaguely) familiar with the terms
category, functor, and natural transformation.

\subsection{Some categories for finite probability theory}

We begin with probability theory as it might be more intuitive than algebra.

\bd
A \emph{\uline{finite probability space}} is a pair $(X,p)$ consisting of 
a finite set $X$ and a function $p:X\to\R$ satisfying
\be
p(x)\ge0\qquad\forall\;x\in X
\ee
and
\be
\sum_{x\in X}p(x)=1.
\ee
The elements of $X$ are called \emph{\uline{events}} 
and $p$ is called a \emph{\uline{probability measure}} on $X.$
The notation $p_{x}:=p(x)$ will often be used. 
\ed

Note that $p$ associates a number to any subset $E\subseteq X$ given by 
\be
p(E):=\sum_{x\in E}p_{x}
\ee
and is the reason we refer to $p$ as a \emph{measure}
(more on this will be discussed in Section \ref{sec:stochasticGelfandNaimark}).
Technically, these \emph{subsets} are called events.

\bd
\label{defn:measpresfnct}
Let $(X,p)$ and $(Y,q)$ be two finite probability spaces. A 
\emph{\uline{probability-}}\-\emph{\uline{preserving}}
function from $(X,p)$ to $(Y,q)$
is a function $f : X \to Y$
satisfying
\be
\label{eq:measurepreserving}
q_{y} = \sum_{x \in f^{-1}(y)} p_{x}
\ee
for all $y \in Y.$
\ed

Composition is given by the usual composition of functions
and the composition of two probability-preserving functions is
immediately seen to be probability-preserving. 
The meaning of a probability-preserving function between
probability spaces can be seen nicely in the following example.

\bx
\label{ex:die} 
Consider the set 
\be
X:=\left\{
\begin{tikzpicture}[baseline=-0.5ex]
\draw[rounded corners=2] (-0.5,-0.5) rectangle (0.5,0.5);
\node[scale=1.25] at (0,0) {$\bullet$};
\end{tikzpicture}
\;,\;
\begin{tikzpicture}[baseline=-0.5ex]
\draw[rounded corners=2] (-0.5,-0.5) rectangle (0.5,0.5);
\node[scale=1.25] at (-0.25,0.25) {$\bullet$};
\node[scale=1.25] at (0.25,-0.25) {$\bullet$};
\end{tikzpicture}
\;,\;
\begin{tikzpicture}[baseline=-0.5ex]
\draw[rounded corners=2] (-0.5,-0.5) rectangle (0.5,0.5);
\node[scale=1.25] at (-0.25,-0.25) {$\bullet$};
\node[scale=1.25] at (0,0) {$\bullet$};
\node[scale=1.25] at (0.25,0.25) {$\bullet$};
\end{tikzpicture}
\;,\;
\begin{tikzpicture}[baseline=-0.5ex]
\draw[rounded corners=2] (-0.5,-0.5) rectangle (0.5,0.5);
\node[scale=1.25] at (-0.25,-0.25) {$\bullet$};
\node[scale=1.25] at (-0.25,0.25) {$\bullet$};
\node[scale=1.25] at (0.25,-0.25) {$\bullet$};
\node[scale=1.25] at (0.25,0.25) {$\bullet$};
\end{tikzpicture}
\;,\;
\begin{tikzpicture}[baseline=-0.5ex]
\draw[rounded corners=2] (-0.5,-0.5) rectangle (0.5,0.5);
\node[scale=1.25] at (-0.25,-0.25) {$\bullet$};
\node[scale=1.25] at (-0.25,0.25) {$\bullet$};
\node[scale=1.25] at (0,0) {$\bullet$};
\node[scale=1.25] at (0.25,-0.25) {$\bullet$};
\node[scale=1.25] at (0.25,0.25) {$\bullet$};
\end{tikzpicture}
\;,\;
\begin{tikzpicture}[baseline=-0.5ex]
\draw[rounded corners=2] (-0.5,-0.5) rectangle (0.5,0.5);
\node[scale=1.25] at (-0.25,-0.25) {$\bullet$};
\node[scale=1.25] at (-0.25,0.25) {$\bullet$};
\node[scale=1.25] at (-0.25,0) {$\bullet$};
\node[scale=1.25] at (0.25,0) {$\bullet$};
\node[scale=1.25] at (0.25,-0.25) {$\bullet$};
\node[scale=1.25] at (0.25,0.25) {$\bullet$};
\end{tikzpicture}
\right\} 
\ee
with the probability measure given by $\frac{1}{6}$ for each element. 
Consider the set 
\be
Y:=\{O,E\}
\ee
consisting of just two elements ($O$ stands for ``odd'' and $E$ stands for ``even'')
with probability measure given by $\frac{1}{2}$ for each element. 
Also consider the function 
$f:X\to Y$ defined by 
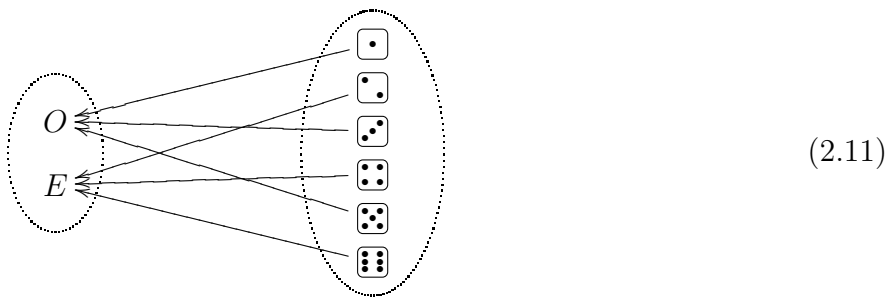
\be
\xy 0;/r.25pc/:
(20,0)*\xycircle(9.0,18.0){.};
(20,13.75)*+{\dieone}="1";
(20,8.25)*+{\dietwo}="2";
(20,2.75)*+{\diethree}="3";
(20,-2.75)*+{\diefour}="4";
(20,-8.25)*+{\diefive}="5";
(20,-13.75)*+{\diesix}="6";
(-20,0)*\xycircle(6.0,10.0){.};
(-20,4)*+{O}="a";
(-20,-4)*+{E}="b";
{\ar"1";"a"};
{\ar"3";"a"};
{\ar"5";"a"};
{\ar"2";"b"};
{\ar"4";"b"};
{\ar"6";"b"};
\endxy
\ee
sending 
\begin{tikzpicture}[scale=0.4,baseline=-0.6ex]
\draw[rounded corners=2] (-0.5,-0.5) rectangle (0.5,0.5);
\node[scale=0.5] at (0,0) {$\bullet$};
\end{tikzpicture}
,
\begin{tikzpicture}[scale=0.4,baseline=-0.6ex]
\draw[rounded corners=2] (-0.5,-0.5) rectangle (0.5,0.5);
\node[scale=0.5] at (-0.25,-0.25) {$\bullet$};
\node[scale=0.5] at (0,0) {$\bullet$};
\node[scale=0.5] at (0.25,0.25) {$\bullet$};
\end{tikzpicture}
, and 
\begin{tikzpicture}[scale=0.4,baseline=-0.6ex]
\draw[rounded corners=2] (-0.5,-0.5) rectangle (0.5,0.5);
\node[scale=0.5] at (-0.25,-0.25) {$\bullet$};
\node[scale=0.5] at (-0.25,0.25) {$\bullet$};
\node[scale=0.5] at (0,0) {$\bullet$};
\node[scale=0.5] at (0.25,-0.25) {$\bullet$};
\node[scale=0.5] at (0.25,0.25) {$\bullet$};
\end{tikzpicture}
to $O$ and sending 
\begin{tikzpicture}[scale=0.4,baseline=-0.6ex]
\draw[rounded corners=2] (-0.5,-0.5) rectangle (0.5,0.5);
\node[scale=0.5] at (-0.25,0.25) {$\bullet$};
\node[scale=0.5] at (0.25,-0.25) {$\bullet$};
\end{tikzpicture}
, 
\begin{tikzpicture}[scale=0.4,baseline=-0.6ex]
\draw[rounded corners=2] (-0.5,-0.5) rectangle (0.5,0.5);
\node[scale=0.5] at (-0.25,-0.25) {$\bullet$};
\node[scale=0.5] at (-0.25,0.25) {$\bullet$};
\node[scale=0.5] at (0.25,-0.25) {$\bullet$};
\node[scale=0.5] at (0.25,0.25) {$\bullet$};
\end{tikzpicture}
, and
\begin{tikzpicture}[scale=0.4,baseline=-0.6ex]
\draw[rounded corners=2] (-0.5,-0.5) rectangle (0.5,0.5);
\node[scale=0.5] at (-0.25,-0.25) {$\bullet$};
\node[scale=0.5] at (-0.25,0.25) {$\bullet$};
\node[scale=0.5] at (-0.25,0) {$\bullet$};
\node[scale=0.5] at (0.25,0) {$\bullet$};
\node[scale=0.5] at (0.25,-0.25) {$\bullet$};
\node[scale=0.5] at (0.25,0.25) {$\bullet$};
\end{tikzpicture}
to $E.$ Therefore, the probability-preserving function 
$f$ is associated to the process of 
rolling a die and considering the likelihood of rolling an odd or an even roll
as opposed to rolling a 1, 2, 3, 4, 5, or 6. 
This describes a deterministic process since we know with certainty
that if we see odd, the die must be 1, 3, or 5, and analogously for even.
\ex

\begin{notation}
\label{not:finprob}
Let $\FinProb$ be the category whose objects are finite probability spaces $(X,p)$
and whose set of morphisms $\big\{(X,p)\to(Y,q)\big\}$ 
are probability-preserving functions from $(X,p)$ to $(Y,q).$
Similarly, let $\FinSetFun$ be the category whose objects are finite sets 
and whose morphisms are functions. 
Composition in both categories is defined to be the ordinary 
composition of functions. 
\end{notation}

The category $\FinProb$ was recently used by Baez, Fritz, and 
Leinster to provide a new categorical characterization of entropy
in the case of finite sets \cite{BFL11}. Extensions of this result
were obtained more recently for relative entropy on finite sets
by Baez and Fritz \cite{BaFr14} and relative entropy on Polish spaces
by Gagn\'e and Panangaden \cite{GaPa17}. 
In all cases, categories of stochastic maps were crucial in 
these characterizations. 
While functions take in one input and spit out one output, 
stochastic maps 
take in one input and have a ``spread'' of outputs
described by some probability distribution. 

\bd
\label{defn:prandstochasticmaps}
Let $X$ and $Y$ be two finite sets. Let $\mathrm{Pr}(Y)$
denote the set of probability measures on $Y.$
A
\emph{\uline{stochastic map}} from $X$ to $Y$ is a function
\be
\begin{split}
X &\xrightarrow{f} \mathrm{Pr}(Y) \\
x &\mapsto f(x)
\end{split}
\ee
whose evaluation on elements in $Y$ is written as 
\be
Y \ni y \xmapsto{f(x)} f_{yx}\in \R_{\ge 0}.
\ee
\ed

Note that by definition of $f(x)$ being a 
probability measure on $Y,$ this means that
\be
\sum_{y \in Y} f_{yx} = 1
\ee
for all $x \in X.$ 
The numbers $\{f_{yx}\}$ labelled by $x\in X$ and $y\in Y$ 
form what is often called a stochastic matrix. 
We denote stochastic maps using squiggly arrows $X \stoch Y$
following the convention of Baez and Fritz \cite{BaFr14}
to distinguish stochastic maps from ordinary functions. 
By abuse of notation, we often use the same notation $f$ for both
$f:X\to\mathrm{Pr}(Y)$ and $f:X\stoch Y$ and we may often refer
to $X\to\mathrm{Pr}(Y)$ as the stochastic map (the use of straight and curly
arrows avoids potential confusion). 

\bx
\label{ex:functionsarestochastic}
Let $X$ and $Y$ be finite sets and let $f:X\to Y$ be a function. 
Associated to $f$ is a canonical stochastic map $f:X\stoch Y$
defined by
\be
X\times Y\ni(x,y)\mapsto{f}_{yx}:=\de_{y f(x)}.
\ee
Here the function $\de:Y\times Y\to\R$ is the Kronecker-delta function
defined by 
\be
\label{eq:Kroneckerdelta}
\de_{y y'}:=
\begin{cases}
1&\mbox{ if $y=y'$}\\
0&\mbox{ otherwise}\\
\end{cases}
.
\ee
In other words, ${f}$ associates to $x$ the probability measure
that is $1$ on the element $f(x)$ and $0$ elsewhere, i.e.
if we have a function, we know with certainty where a point will go. 
\ex

\bx
Recall Example \ref{ex:die} with rolling a die. If instead of 
interpreting the set $Y$ as the set of odd or even die, suppose
that the $O$ and $E$ just stand for boxes labelled by $O$ or $E.$ 
A probability-preserving stochastic map $X\stoch Y$ can be interpreted
as saying that for any given roll of a die, that particular side has a 
probability of being placed in box $O$ and a complementary probability
of being placed in box $E.$ 
\ex

\bx
\label{ex:randomwalk}
A random walk on a circle is described by the following stochastic map. 
Let $N\in\N$ be the number of points on a circle, labeled 
$x_0,x_1,\dots,x_{N-1}$ with $x_{N}\equiv x_{0}$ and so on
so that $X:=\{x_0,x_1,\dots,x_{N-1}\}.$ The random walk on $X$ is
described by the stochastic map $f:X\stoch X$ defined by 
\be
f_{ji}:=
\begin{cases}
\frac{1}{2}&\mbox{ if }j\equiv i\pm1\!\!\mod N\\
0&\mbox{ otherwise}
\end{cases}
\ee
where the shorthand notation $f_{ji}$ is used instead of $f_{x_j x_i}.$ 
The reason this is called the random walk is because starting
at any point, you have a $\frac{1}{2}$ probability of moving
one step in either direction. A random walk on a line or plane 
can also be described though the set of points would be infinite, 
which we are not considering at the moment. 
\ex

\bx
\label{ex:stochasticfrominitial}
With such a definition in place, notice that the set of all stochastic maps
$\{\bullet\}\stoch X$ from a single element set $\{\bullet\}$ to a finite set $X$ is
isomorphic (bijective) 
to the set $\mathrm{Pr}(X)$ of all probability measures on $X.$
\ex

Given finite sets $X, Y,$ and $Z$
together with stochastic maps $f : X \stoch Y$ and $g : Y \stoch Z,$ 
one would like to iterate these maps. To do so, the probability
of $x\in X$ evolving to $z\in Z$ is given by summing over all 
possible intermediate elements $y\in Y$ weighted by their corresponding
probabilities. 
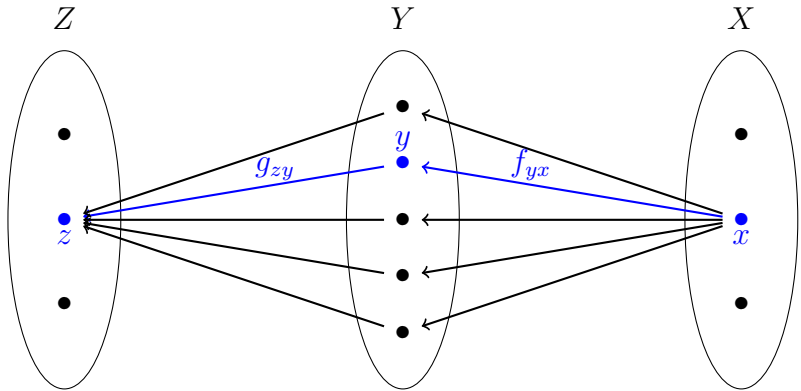
\begin{center}
\begin{tikzpicture}[scale=0.75]
\draw (6,0) ellipse (1cm and 3cm) node[above,yshift=2.4cm]{$X$};
\node at (6,1.5) {$\bullet$};
\node at (6,-1.5) {$\bullet$};
\node[blue] (x) at (6,0) {$\bullet$};
\node[blue] at (6,-0.3) {$x$};
\draw (0,0) ellipse (1cm and 3cm) node[above,yshift=2.4cm]{$Y$};
\node (0) at (0,0) {$\bullet$};
\node[blue] (y) at (0,1) {$\bullet$};
\node (2) at (0,2) {$\bullet$};
\node[blue] at (0,1.4) {$y$};
\node (n1) at (0,-1) {$\bullet$};
\node (n2) at (0,-2) {$\bullet$};
\draw (-6,0) ellipse (1cm and 3cm) node[above,yshift=2.4cm]{$Z$};
\node[blue] (z) at (-6,0) {$\bullet$};
\node[blue] at (-6,-0.3) {$z$};
\node at (-6,1.5) {$\bullet$};
\node at (-6,-1.5) {$\bullet$};
\draw[thick,->] (x) -- (0);
\draw[thick,->] (x) -- (2);
\draw[thick,->] (x) -- (n1);
\draw[thick,->] (x) -- (n2);
\draw[thick,->,blue] (x) -- node[above,xshift=-0.55cm]{$f_{yx}$} (y);
\draw[thick,->,blue] (y) -- node[above,xshift=0.55cm]{$g_{zy}$} (z);
\draw[thick,->] (n2) -- (z);
\draw[thick,->] (n1) -- (z);
\draw[thick,->] (0) -- (z);
\draw[thick,->] (2) -- (z);
\end{tikzpicture}
\end{center}
This motivates the following definition. 

\bd
Let $X, Y,$ and $Z$ be finite sets.
The \emph{\uline{composition}} of the
stochastic map $f : X \stoch Y$ followed by the stochastic map $g : Y \stoch Z,$ 
written as $g \circ f : X \stoch Z,$ is the function
$g \circ f : X \to \mathrm{Pr}(Z)$
defined by sending $x\in X$ to the probability
measure on $Z$ defined by
\be
\begin{split}
Z \ni z \xmapsto{(g\circ f)(x)} (g\circ f)_{zx}
:= \sum_{y\in Y}g_{zy}f_{yx}.
\end{split}
\ee
\ed

Note that $g\circ f$ is indeed a stochastic map because 
\be
\sum_{z \in Z} (g\circ f)_{zx}
\overset{\text{def}}{=}\sum_{z \in Z}\sum_{y\in Y}g_{zy}f_{yx}
=\sum_{y\in Y}\underbrace{\sum_{z \in Z}g_{zy}}_{1}f_{yx}
=\sum_{y\in Y}f_{yx}
=1.
\ee

\begin{notation}
\label{not:finstochandfinsetstoch}
Let $\FinSetStoch$ be the category whose objects are
finite sets and whose morphisms are stochastic maps. 
We leave checking the axioms of a category to the reader. 
\end{notation}

Notice that $\FinSetFun$ can be viewed as a subcategory of $\FinSetStoch$ by 
Example \ref{ex:functionsarestochastic} and we therefore
denote the associated functor by $\de:\FinSetFun\to\FinSetStoch.$ 

\br
The definition of a stochastic map is closely related to the notion
of a stochastic process. 
A stochastic map $f:X\stoch X$ from a finite set $X$
to itself can be used to describe the possible evolution
of a state. In the special case of a known
initial condition, i.e. an element $x_0\in X,$
one can consider the sequence of probability distributions 
(the first element of which is just the Kronecker delta
distribution at $x_{0}$) given by 
$\big(x_0, f(x_0),f^2(x_0),f^3(x_0), f^4(x_0), \dots\big).$ 
Such a sequence is called a Markov chain.
If each iteration of $f$ is interpreted as a time-step, 
then the probability distribution $f^n(x_0)$ is to be interpreted
in the following way. The probability of finding the initial condition
$x_0$ to evolve to the position $x_n\in X$ after $n$ time steps
is $(f^n)_{x_n x_0}.$ 
More generally, the
set of events can change under a given process (as we have
described above). An evolution of some initial states
can therefore be described by a sequence of stochastic maps
\be
X_0\stoch X_1\stoch X_2\stoch X_3\stoch\cdots,
\ee
where each of the $X_{i}$ is a set of possible events. 
A point in $X_{0},$ described by a function $\{\bullet\}\to X_{0},$ 
evolves just as a probability measure $\{\bullet\}\stoch X_{0},$
namely by precomposing:
\be
\{\bullet\}\stoch X_0\stoch X_1\stoch X_2\stoch X_3\stoch\cdots. 
\ee
\er

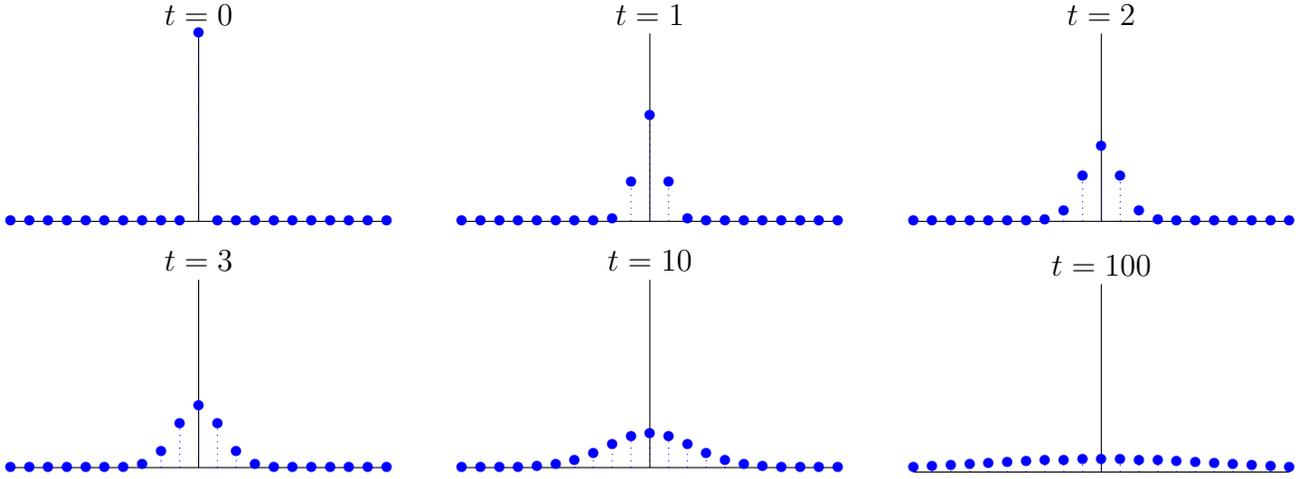
\begin{figure}
\centering
\begin{tikzpicture}[scale=2.5]
    \draw[-] (-1,0) -- (1,0);
    \draw[-] (0,0) -- (0,1.00);
    \node at (0,1.1) {$t=0$};
    \draw[dotted,blue] (0,0) -- (0,1) node {\footnotesize$\bullet$};
    \foreach \k in {-10,-9,...,-1}
    {
    \draw[dotted,blue] ({\k/10},0) -- ({\k/10},0) node {\footnotesize$\bullet$};
    }
    \foreach \k in {1,2,...,10}
    {
    \draw[dotted,blue] ({\k/10},0) -- ({\k/10},0) node {\footnotesize$\bullet$};
    }
\end{tikzpicture}
\quad
\begin{tikzpicture}[scale=2.5]
    \draw[-] (-1,0) -- (1,0);
    \draw[-] (0,0) -- (0,1.00);
    \node at (0,1.1) {$t=1$};
    \foreach \k in {-10,-9,...,10}
    {
    \draw[dotted,blue] ({\k/10},0) -- ({\k/10},{(1/1.77264)*exp(-\k*\k)}) node {\footnotesize$\bullet$};
    }
\end{tikzpicture}
\quad
\begin{tikzpicture}[scale=2.5]
    \draw[-] (-1,0) -- (1,0);
    \draw[-] (0,0) -- (0,1.00);
    \node at (0,1.1) {$t=2$};
    \foreach \k in {-10,-9,...,10}
    {
    \draw[dotted,blue] ({\k/10},0) -- ({\k/10},{(1/2.50663)*exp(-\k*\k/2)}) node {\footnotesize$\bullet$};
    }
\end{tikzpicture}
%
\begin{tikzpicture}[scale=2.5]
    \draw[-] (-1,0) -- (1,0);
    \draw[-] (0,0) -- (0,1.00);
    \node at (0,1.1) {$t=3$};
    \foreach \k in {-10,-9,...,10}
    {
    \draw[dotted,blue] ({\k/10},0) -- ({\k/10},{(1/3.06998)*exp(-\k*\k/3)}) node {\footnotesize$\bullet$};
    }
\end{tikzpicture}
\quad
\begin{tikzpicture}[scale=2.5]
    \draw[-] (-1,0) -- (1,0);
    \draw[-] (0,0) -- (0,1.00);
    \node at (0,1.1) {$t=10$};
    \foreach \k in {-10,-9,...,10}
    {
    \draw[dotted,blue] ({\k/10},0) -- ({\k/10},{(1/5.60498)*exp(-\k*\k/10)}) node {\footnotesize$\bullet$};
    }
\end{tikzpicture}
\quad
\begin{tikzpicture}[scale=2.5]
    \draw[-] (-1,0) -- (1,0);
    \draw[-] (0,0) -- (0,1.00);
    \node at (0,1.1) {$t=100$};
    \foreach \k in {-10,-9,...,10}
    {
    \draw[dotted,blue] ({\k/10},0) -- ({\k/10},{(1/15.2921)*exp(-\k*\k/100)}) node {\footnotesize$\bullet$};
    }
\end{tikzpicture}
\caption{A sharp Gaussian representing a definite initial condition at the origin 
tends towards a uniform distribution
after several iterations of the heat flow described by the matrix 
(\ref{eq:heatflowmatrix}).
$t=n$ with $n\in\{0\}\cup\N$ stands for the $n$-th time step. 
}
\label{fig:heatfinite}
\end{figure}

\bx
Let us model one-dimensional diffusion via a stochastic map on
a finite set, thought of as an approximation to diffusion on the interval $[-1,1]$
with periodic boundary conditions. Namely, let $X$ be the finite set given by 
\be
X:=\left\{\frac{k}{10}\;:\;k\in\{-10,-9,\dots,9\}\right\},
\ee
where the point $1$ is to be interpreted as being identified with $-1.$ 
According to the heat diffusion equation, after some time step,
a particle at a point has the following probability distribution.
Approximately, it has a $0.56$ likelihood to stay in the same place, 
$0.21$ likelihood to move one unit over in either direction, and a 
$0.01$ likelihood to move two units over in either direction. 
The stochastic matrix associated to this process is given by
\be
\label{eq:heatflowmatrix}
\begin{bmatrix}
0.56&0.21&0.01&0&\cdots&0.01&0.21\\
0.21&0.56&0.21&0.01&\cdots&0&0.01\\
0.01&0.21&0.56&0.21&\cdots&0&0\\
0&0.01&0.21&0.56&\cdots&0&0\\
\vdots&\vdots&\vdots&\vdots&\ddots&\vdots&\vdots\\
0.01&0&0&0&\cdots&0.56&0.21\\
0.21&0.01&0&0&\cdots&0.21&0.56
\end{bmatrix}
\ee
Under a few iterations of this process (matrix multiplication), the distribution (described by a unit vector) becomes
the uniform distribution (as one would expect---systems tend
to equilibriate). This is depicted in 
Figure \ref{fig:heatfinite}. 
\ex

Note that the formalism described above a-priori only describes discrete stochastic 
processes. One can adapt these ideas to study continuous processes
though we do not discuss this here except for an example from quantum
mechanics in Section \ref{sec:examplesinQM}. The reader is referred to
Giry's work, which is not only a foundational paper on categorical 
probability theory but also discusses this in some detail \cite{Gi82}. 

\subsection{$C^*$-algebras and states}

There are several definitions that must be provided
so that this work can be somewhat self-contained. 
A more thorough account of many definitions and results can be found in
Dixmier's book \cite{Di77} with additional references provided throughout.

\bd
\label{defn:algebra}
An (associative and unital) \emph{\uline{algebra}} 
is a vector space $\mA$ (over $\C$) together with 
\begin{enumerate}[(a)]
\setlength{\itemsep}{0pt}
\item
a binary multiplication operation
$\mA\times\mA\to\mA,$
\item
an element $1_{\mA}\in\mA\setminus\{0\}.$ 
\end{enumerate}
These data must satisfy the following conditions:
\begin{enumerate}[i.]
\setlength{\itemsep}{0pt}
\item
$(ab)c=a(bc)$ (associativity of multiplication),
\item
$a(b+c)=ab+ac$ and $(a+b)c=ac+bc$
(distributivity of multiplication over vector addition),
\item
$(\l a)b=a(\l b)=\l(ab)$, 
(distributivity of multiplication over scalar multiplication),
and 
\item
$a1_{\mA}=1_{\mA}a=a$
($1_{\mA}$ is a unit for the multiplication), 
\end{enumerate}
for all $a,b,c\in\mA$ and $\l\in\C.$ 
\ed

\bd
\label{defn:normedalgebra}
An algebra $\mA$ together with a norm
$\lVert \ \cdot \ \rVert:\mA\to\R_{\ge0}$ for which
$\mA$ with this norm is a normed vector space and
\be
\lVert ab\rVert\le\lVert a\rVert\lVert b\rVert\qquad\forall\;a,b\in\mA
\ee
is a \emph{\uline{normed algebra}}. 
$\mA$ is a \emph{\uline{Banach algebra}} iff 
all Cauchy sequences converge with respect to this norm.
\ed

\bd
An algebra $\mA$ together with an anti-automorphism 
${}^{*}:\mA\to\mA,$ i.e.
\be
(\l a+b)^{*}=\l^{*}a^{*}+b^{*}
\qquad\forall\;\l\in\C,\;a,b\in\mA
\ee
and
\be
(ab)^*=b^*a^*\qquad\forall\;a,b\in\mA,
\ee
is an \emph{\uline{involutive algebra}}.
$*$ is known as an \emph{\uline{involution}} on $\mA.$ 
\ed

\bd
Let $\mA$ be both an involutive algebra and a normed algebra.
If, in addition, $\mA$ satisfies 
\be
\lVert a^*\rVert=\lVert a\rVert\qquad\forall\;a\in\mA
\ee
then $\mA$ is called an \emph{\uline{normed involutive algebra}}.
If, in addition, all Cauchy sequences in $\mA$ converge, 
then $\mA$ is called a \emph{\uline{involutive Banach algebra}}. 
\ed

\bd
\label{defn:C*algebra}
A \emph{\uline{$C^*$-algebra}} is an involutive Banach algebra $\mA$ for which 
\be
\lVert a^*a\rVert=\lVert a\rVert^2 \qquad\forall\;a\in\mA.
\ee
A $C^*$-algebra is \emph{\uline{commutative}} iff 
$ab=ba$ for all $a,b\in\mA.$ 
\ed

\br
If $\mA$ is a normed algebra with an involution for which 
$\lVert a\rVert^2\le\lVert a^*a\rVert$ for all $a\in\mA,$
then $\mA$ is an involutive algebra. In fact, if $\mA$
is a Banach algebra with involution satisfying this condition, 
then $\mA$ is a $C^*$-algebra \cite{Di77}. 
\er

All of the data of these algebraic structures are
summarized in the following table (with no reference
to their conditions). 

\begin{center}
\begin{tabular}{|c|c|c|c|c|c|c|}
\hline
\begin{tabular}{c}
data/\\
structure
\end{tabular}
&
\begin{tabular}{c}
vector\\
space (VS)
\end{tabular}
&algebra
&\begin{tabular}{c}
normed\\
VS
\end{tabular}
&\begin{tabular}{c}
Banach\\
algebra
\end{tabular}
&
\begin{tabular}{c}
involutive\\
algebra (IA)
\end{tabular}
&
\begin{tabular}{c}
normed IA/\\
$C^*$-algebra
\end{tabular}
\\
\hline
$0$
&$\checkmark$&$\checkmark$&$\checkmark$&$\checkmark$&$\checkmark$&$\checkmark$\\
\hline
$+$
&$\checkmark$&$\checkmark$&$\checkmark$&$\checkmark$&$\checkmark$&$\checkmark$\\
\hline
scalar mult.
&$\checkmark$&$\checkmark$&$\checkmark$&$\checkmark$&$\checkmark$&$\checkmark$\\
\hline
$\times$
&&$\checkmark$&&$\checkmark$&$\checkmark$&$\checkmark$\\
\hline
$\lVert\;\cdot\;\rVert$
&&&$\checkmark$&$\checkmark$&&$\checkmark$\\
\hline
$*$
&&&&&$\checkmark$&$\checkmark$\\
\hline
$1$
&&$\checkmark$&&$\checkmark$&$\checkmark$&$\checkmark$\\
\hline
\end{tabular}
\end{center}

\bd
\label{defn:c*algebramorphism}
Let $\mA$ and $\mB$ be two $C^*$-algebras. 
A \emph{\uline{$*$-homorphism of $C^*$-algebras}}
from $\mB$ to $\mA$ is a bounded (i.e. continuous) linear map
$f: \mB\to\mA$  
such that 
\begin{enumerate}[i.]
\setlength{\itemsep}{0pt}
\item
$f(b^*)=f(b)^*$,
\item
$f(b_{1}b_{2})=f(b_1)f(b_2),$ and
\item
$f(1_{\mB})=1_{\mA}$ 
\end{enumerate}
for all $b,b_1,b_2\in\mB.$
\ed

\br
The word bounded in the definition of a $*$-homomorphism
of $C^*$-algebras is redundant. 
It is a theorem that such a map satisfying
conditions i., ii., and iii. is automatically bounded \cite{Di77}.
Thus, the reader may safely ignore that condition. 
\er

There are a few basic examples of $C^*$-algebras that we will focus on in 
this section.
\bx
$\C$ is a $C^*$-algebra with its usual structure associated with
complex numbers.
The involution $*$ is the complex conjugate. 
This is an example of a commutative $C^*$-algebra.
\ex

This previous example automatically induces a huge class
of examples. 

\bx
\label{ex:CXalgebra}
Let $X$ be a finite set and let $\C^{X}$ denote the set of all
functions $X\to\C.$ There is a unique algebra structure and involution
on $\C^{X}$ so that the evaluation functions $\ev_{x}:\C^{X}\to\C,$
defined by sending $\vf:X\to\C$ to $\vf(x),$
are $*$-homomorphisms for all $x\in X$
(these are just the pointwise algebraic structures).%
\footnote{This way of phrasing the natural algebraic structure on $\C^{X}$ 
makes precise in what sense it is the ``obvious'' algebraic structure.
Note that $X$ does not have to be finite for this statement. 
However, one must be careful about the topology when $X$ is 
infinite---demanding the smallest topology for which $\ev_{x}$
are all continuous would result in the topology of pointwise convergence, 
which is not the one we will use later. 
}
The norm on $\C^{X}$ is given by 
\be
\C^{X}\ni f\mapsto \lVert f\rVert:=\sup_{x\in X}\{|f(x)|\}.
\ee
\ex

The previous examples were commutative $C^*$-algebras. 
The following example is the quintessential $C^*$-algebra
of quantum mechanics of a finite system (such as a qubit).

\bx
\label{ex:nbybmatrices}
Fix $n\in\N.$ The set $M_{n}(\C)$ 
of $n\times n$ complex-valued matrices is a $C^*$-algebra with 
product given by matrix multiplication and $*$ is the transpose
complex conjugate. The norm can be taken to be
either the operator norm or the standard Euclidean norm 
(the topology will be the same so it does not matter). 
\ex

It is a theorem that every finite-dimensional $C^*$-algebra
is ($*$-isomorphic to) a direct sum of these (Theorem III.1.1 in \cite{Da96}).

\bx
More generally, for any
Hilbert space $\Hi,$ the set of bounded operators $\mB(\Hi)$ on
$\Hi$ is a $C^*$-algebra via operator composition, $*$ given by the adjoint, 
and the norm given by the operator norm. 
\ex

\begin{notation}
\label{defn:CAlg}
Let $\CAlg,$ $\cCAlg,$ and $\fdcCAlg$ denote the categories whose 
objects are $C^*$-algebras, commutative $C^*$-algebras, 
and finite-dimensional commutative $C^*$-algebras, respectively, 
and whose the morphisms are all $*$-homomorphisms of $C^*$-algebras. 
\end{notation}

\bd
\label{defn:state}
Given a $C^*$-algebra $\mA,$ 
a \emph{\uline{state}} on $\mA$ is a
bounded (i.e. continuous)%
\footnote{The word bounded/continuous here is redundant. 
It is a theorem that such a map on a $C^*$-algebra satisfying these
conditions is automatically bounded.}
linear function ${\w:\mA\to\C}$ such that
$\w(1_{\mA})=1$ and $\w(a^*a)\ge0$ for all $a\in\mA.$
\ed

The following theorem relates states to density matrices and motivates
why density matrices are the states of quantum mechanics 
as opposed to vectors (really, rays) in a Hilbert space \cite{Ha13}.

\bt
\label{thm:statesdensitymatrices}
There is a one-to-one correspondence between 
states on $M_{n}(\C)$ and matrices $\rho\in M_{n}(\C)$
satisfying $\tr(\rho)=1$ and $\rho^*=\rho$
(such matrices are known as \emph{\uline{density matrices}}). 
More precisely, for any density matrix $\rho,$ the function 
$\tr(\rho\;\cdot\;):M_{n}(\C)\to\C,$ sending $A\in M_{n}(\C)$ to
$\tr(\rho A),$ is a state. Conversely, for
any state $\w:M_{n}(\C)\to\C,$ there exists a unique density matrix
$\rho$ such that $\w=\tr(\rho\;\cdot\;).$ 
\et

This theorem extends to states on $\mB(\Hi)$ with $\Hi$ a separable
Hilbert space
and states that satisfy an additional continuity assumption.%
\footnote{
Such states go by the name of \emph{\uline{normal states}} for those readers 
who would like to look up more information. However, normal states
require a certain notion of continuity not available on arbitrary 
$C^*$-algebras, and one requires the notion of a von Neumann algebra, 
a special type of $C^*$-algebra that is closed with respect to a certain topology
\cite{ReSu07}.
$\mB(\Hi)$ is an example of a von Neumann algebra. 
}
In the context of quantum mechanics, $\mB(\Hi)$ is the set of
physical observables including non self-adjoint operators, 
such as ladder operators, so that it remains closed under 
multiplication by complex numbers. 
A state on $\mB(\Hi)$ assigns to every observable a number, 
interpreted as the expectation value of that observable
for that state.
States on $C^*$-algebras and their relationship to the more standard
approach of quantum mechanics in terms of Hilbert spaces is 
described in a categorical framework in \cite{Pa16}. 

\subsection{Two types of morphisms of $C^*$-algebras}

There is another notion of a morphism that is important
for $C^*$-algebras and utilizes less of its structure. 

\bd
Let $\mA$ be a $C^*$-algebra. An element $a\in\mA$ is 
\emph{\uline{positive}} iff there exists an element $b\in\mA$ such that
$a=b^*b.$ Positivity of elements induces a partial ordering on $\mA$
and the notation $a\ge0$ will be used to indicate that $a$ is positive
($b\ge a$ iff $b-a\ge0$).
\ed

\bd
Let $\mA$ and $\mB$ be two $C^*$-algebras. 
A \emph{\uline{positive (unital) map}} from $\mB$ to $\mA$ is a 
bounded linear function $\varphi:\mB\to\mA$ such that 
$\vf(1_{\mB})=1_{\mA}$ and $\vf(b^*b)\ge0$ for all $b\in\mB.$ 
\ed

A state is an example of a positive map. 
Every $*$-homomorphism of $C^*$-algebras is a positive map.
The converse, however, is far from true.
Examples are given in Section \ref{sec:examplesinQM}. 
Non-unital positive maps seem to be the weakest structure 
needed to make sense of operations and processes 
on quantum mechanical systems though one may also relax the
unit-preserving condition \cite{Kr83}.

\begin{notation}
Let $\mA$ be an involutive algebra and let $n\in\N.$ 
$M_{n}(\mA)$ is the set of of $n\times n$ matrices with
coefficients in $\mA$ with product given by matrix multiplication
and $*$ given by the transpose and applying $*$ in each entry. 
\end{notation}

\bt
Let $\mA$ be a $C^*$-algebra.
Then there exists a unique norm $\lVert\;\cdot\;\rVert$ on $M_{n}(\mA)$ 
such that 
\be
\lVert AB\rVert\le\lVert A\rVert \lVert B\rVert\qquad\forall\;A,B\in M_{n}(\mA)
\ee
and
\be
\lVert A^*A\rVert=\lVert A\rVert^2\qquad\forall\;A\in M_{n}(\mA).
\ee
Furthermore, $M_{n}(\mA)$ with this structure is a $C^*$-algebra.
\et

\bprf
See Chapter 1 of Paulsen \cite{Pa03} 
or Chapters 4 and 6 of Warner \cite{Wa10}.
\eprf

\br
There are $*$-isomorphisms
\be
M_{n}(\mA)\cong M_{n}(\C)\otimes\mA\cong \mA\otimes M_{n}(\C).
\ee
Therefore, all of these involutive algebras have a unique norm 
giving them the structure of $C^*$-algebras. 
The theory of tensor products of $C^*$-algebras is notoriously complicated.
Since we will work with commutative $C^*$-algebras for the most part, we
will avoid such complications. 
\er

\bd
\label{defn:CPU}
Let $\mA$ and $\mB$ be two $C^*$-algebras. 
A \emph{\uline{completely positive (unital) map}} from $\mB$ to $\mA$ is a 
function $\varphi:\mB\to\mA$ such that 
$\vf\otimes\id_{n}:\mB\otimes M_{n}(\C)\to\mA\otimes M_{n}(\C)$
is a positive map for all $n\in\N.$ The abbreviation ``CP map'' will often be used
for such maps.
\ed

\br
One motivation for this definition comes from
quantum mechanics \cite{Li76}, \cite{Kr83}. 
The concept of tensoring with the algebra
$M_{n}(\C)$ and the positive map with the identity 
means to adjoin another finite system that does not interact
with the original system on $\mB$ during the operation. The resulting map 
$\vf\otimes\id_{n}:\mB\otimes M_{n}(\C)\to\mA\otimes M_{n}(\C)$
should therefore still be positive. There are examples of positive
maps that are not completely positive so requiring this condition to be
satisfied is an additional constraint \cite{Pa03}, \cite{Kr83}.
Unital maps of $C^*$-algebras correspond to the trace-preserving
condition that is often referred to in texts on quantum information
\cite{NiCh10}. See Example \ref{ex:unitaryevolution} for an illustration. 
All maps in this work are unital so we will not use the
more standard notation ``CPU'' \cite{NiCh10}.
\er

The composition of CP maps is CP \cite{Pa03}. 

\bt
\label{thm:positivebetweencommutative}
Let $\mA$ be a $C^*$-algebra. 
Every positive map $\mA\to\C$ is CP.
In particular, a state on $\mA$ is CP.
Furthermore, every positive map between commutative $C^*$-algebras
is CP.
\et

\bprf
See Theorems 3.9 and 3.11 in \cite{Pa03}.
\eprf

\begin{notation}
\label{not:starhomoandCP}
Henceforth, a $*$-homomorphism from a $C^*$-algebra $\mB$ to 
another $C^*$-algebra $\mA$ will be denoted by a 
straight arrow $\mB\to\mA.$ A CP map, 
on the other hand, will be denoted with a curvy arrow
$\mB\stoch\mA.$ 
\end{notation}

We now collect all the categories of $C^*$-algebras that will be needed
in the prequel. 

\begin{notation}
\label{defn:CPCAlg}
Let $\CPCAlg,$ $\cCAlgPos,$ and $\fdcCAlgPos$ denote the categories whose 
objects are $C^*$-algebras, commutative $C^*$-algebras, 
and finite-dimensional commutative $C^*$-algebras, respectively, 
and whose the morphisms are all CP maps of $C^*$-algebras
(note that in the latter two cases, CP maps are the same as 
positive maps). 
\end{notation}

On a given physical system corresponding to some $C^*$-algebra, 
a state provides expectation values for observables. 
One can define a category of $C^*$-algebras equipped with a state
together with the various flavors of morphisms analogous to the case
of finite probability. However, this is tangential to our main goal
and will follow as an immediate consequence anyway. 
Examples are provided in Section \ref{sec:examplesinQM}, but first
we connect algebra with probability theory.

\subsection{From probability theory to algebra}
\label{sec:fromfiniteprobtoalg}

Now, we will construct a functor 
$C:\FinSetStoch\to\fdcCAlgPos$ and show that it restricts to a well-defined 
functor $C:\FinSetFun\to\fdcCAlg.$ The same notation for these functors will
be used as it will be clear from context which one is being referred to. 
The properties of these functors will be discussed 
in Section \ref{sec:fdfullyfaithful}. 
Associated to each finite set $Y,$ one has the $C^*$-algebra
${\C}^{Y}$ discussed in Example \ref{ex:CXalgebra}. 
Note that any element $v\in\C^{Y}$ can be uniquely expressed as 
\be
v=\sum_{y\in Y}v_{y}e_{y},
\ee
where $e_{y}:Y\to\C$ is defined by 
\be
Y\ni y'\mapsto e_{y}(y'):=\de_{yy'}
\ee
and $v_{y}\in\C$ is given by $v(y)=:v_{y}.$ 
Hence, if $X$ is another finite set, 
any CP map $\vf:\C^{Y}\stoch\C^{X}$
is determined by its value on the set of elements $\{e_{y}\}_{y\in Y}$ by linearity. 
Namely, 
\be
\vf(v)=\sum_{y\in Y}v_{y}\vf(e_{y}).
\ee
Hence, to every stochastic map $f:X\stoch Y,$ associate the linear function
$C^{f}:\C^{Y}\stoch\C^{X}$ 
uniquely determined by 
\be
\label{eq:curlyEf}
C^{f}(e_{y}):=\sum_{x\in X}f_{yx}e_{x}
\ee
for all $y\in Y.$ 
These assignments define the functor $C.$ 

\bt
\label{thm:FinProbtoCP}
The assignment%
\footnote{The $\op$ superscript is only meant to keep track that the
directionality of the arrows get reversed. This is shown in 
(\ref{eq:curlyEf}), where the domain of $C^{f}$ is $\C^{Y}$ instead of
$\C^{X}.$ The directionality change is also shown in (\ref{eq:Ccofunctoriality}). 
}
\be
\begin{split}
\FinSetStoch^{\op}&\xrightarrow{C}\fdcCAlgPos\\
X&\mapsto\C^{X}\\
\Big(f:X\stoch Y\Big)&\mapsto\Big(C^{f}:\C^{Y}\stoch\C^{X}\Big),
\end{split}
\ee
where $C^{f}$ is defined in (\ref{eq:curlyEf}), 
is a well-defined functor. 
\et

\bprf
The proof will be broken into three steps.
It will be shown that 
$C^{f}$ is CP,
the identity stochastic map goes to the identity CP map, 
and the composition of stochastic maps gets sent
to the composition of CP maps. 
\begin{enumerate}[i.]
\item
By Theorem \ref{thm:positivebetweencommutative}, 
it suffices to show that $C^{f}$ is positive. 
Every positive element in $\C^{Y}$ is of the form 
\be
\sum_{y\in Y}a_{y}e_{y},\qquad\text{with }\;a_{y}\ge0\quad\forall\;y\in Y.
\ee
Hence, by linearity, it suffices to show that $C^{f}(e_{y}),$
given by (\ref{eq:curlyEf}), 
is positive for every $y\in Y.$ By definition of $f$ being a stochastic
map, $f_{yx}\ge0$ for all $y\in Y$ and $x\in X.$ Hence, 
$C^{f}(e_{y})$ is a positive element in $\C^{X}.$
Furthermore, $C^{f}$ is unital because
\be
C^{f}(1_{Y})=C^{f}\left(\sum_{y\in Y}e_{y}\right)
=\sum_{y\in Y}\big(C^{f}(e_{y})\big)
=\sum_{y\in Y}\sum_{x\in X}f_{yx}e_{x}
=\sum_{x\in X}\underbrace{\left(\sum_{y\in Y}f_{yx}\right)}_{1}e_{x}
=1_{X}.
\ee

\item
The identity function $\id:X\to X$ has the corresponding stochastic map
given by the probability measure $\id_{x'x}=\de_{x'x}$ for all $x,x'\in X.$ 
Therefore, 
\be
C^{\id}(e_{x'})=\sum_{x\in X}\id_{x'x}e_{x'}
=\sum_{x\in X}\de_{x'x}e_{x}
=e_{x'}
\ee
for all $x'\in X.$ Thus, $C^{\id}$ is the identity CP map.

\item
Now suppose that $f:X\stoch Y$ and $g:Y\stoch Z$ are two
composable stochastic maps between finite sets.
This induces CP maps 
$C^{f}:\C^{Y}\stoch\C^{X}$
and $C^{g}:\C^{Z}\stoch\C^{Y}.$
Therefore, 
\be
\label{eq:Ccofunctoriality}
\begin{split}
C^{f}\big(C^{g}(e_{z})\big)
&=C^{f}\left(\sum_{y\in Y}g_{zy}e_{y}\right)\\
&=\sum_{y\in Y}g_{zy}C^{f}(e_{y})\\
&=\sum_{y\in Y}g_{zy}\sum_{x\in X}f_{yx}e_{x}\\
&=\sum_{x\in X}\underbrace{\left(\sum_{y\in Y}g_{zy}f_{yx}\right)}_{(g\circ f)_{zx}}e_{x}\\
&=C^{g\circ f}(e_{z})
\end{split}
\ee
for all $z\in Z.$ 
\end{enumerate} 
\eprf

\bt
\label{thm:restrictEtofinprob}
The functor $C:\FinSetStoch^{\op}\to\fdcCAlgPos$ from Theorem 
\ref{thm:FinProbtoCP} restricts to
a well-defined functor $C:\FinSetFun^{\op}\to\fdcCAlg.$ 
\et

\bprf
Using the notation of Theorem \ref{thm:FinProbtoCP} and its proof, 
all that is left to show is that $C^{f}:\C^{Y}\stoch\C^{X}$ is a 
$*$-homomorphism when $f$ is a function. By linearity, it suffices to prove that 
$C^{f}(e_{y}e_{y'})=C^{f}(e_{y})C^{f}(e_{y'})$
for all $y,y'\in Y.$ This follows from%
\footnote{Diagrams like this are to be read from top to bottom going along
either direction to recreate the thought-process involved in proving the claim.
The reason is because at some point near the bottom, a step which does not
seem obvious in one direction is perhaps more intuitive from the other 
direction (namely, $\de_{yf(x)}^{2}=\de_{yf(x)}$).}
\be
\xy0;/r.50pc/:
(6.84040286651,18.7938524157)*+{\ds C^{f}(e_y)C^{f}(e_{y'})}="1";
(17.3205080757,10)*+{\ds\sum_{x,x'\in X}\de_{yf(x)}\de_{y'f(x')}e_{x}e_{x'}}="2";
(19.6961550602,-3.47296355334)*+{\ds\sum_{x,x'\in X}\de_{yf(x)}\de_{y'f(x')}\de_{xx'}e_{x'}}="3";
(12.8557521937,-15.3208888624)*+{\;\;\ds\sum_{x\in X}\de_{yf(x)}\de_{y'f(x)}e_{x}}="4";
(0,-20)*+{\ds\de_{yy'}\sum_{x\in X}\de_{yf(x)}^2e_{x}}="5";
(-12.8557521937,-15.3208888624)*+{\ds\de_{yy'}\sum_{x\in X}\de_{yf(x)}e_{x}\;\;}="6";
(-19.6961550602,-3.47296355334)*+{\ds\de_{yy'}C^{f}(e_{y})}="7";
(-17.3205080757,10)*+{\ds C^{f}(\de_{yy'}e_{y})}="8";
(-6.84040286651,18.7938524157)*+{\ds C^{f}(e_{y}e_{y'})}="9";
{\ar@{=}@/^0.65pc/"1";"2"};
{\ar@{->}@<1.0ex>@/^0.50pc/(12.04040286651,17.2938524157);"2"};
{\ar@{=}@/^0.65pc/"2";"3"};
{\ar@{->}@<2.0ex>@/^0.65pc/"2";"3"};
{\ar@{=}@/^0.65pc/"3";"4"};
{\ar@{->}@<2.0ex>@/^0.65pc/"3";"4"};
{\ar@{=}@/^0.65pc/"4";"5"};
{\ar@{->}@<2.0ex>@/^0.65pc/"4";"5"};
{\ar@{=}@/^0.65pc/"5";"6"};
{\ar@{->}@<2.0ex>@/^0.65pc/"5";"6"};
{\ar@{=}@/^0.65pc/"6";"7"};
{\ar@{->}@<-2.0ex>@/_0.65pc/"7";"6"};
{\ar@{=}@/^0.65pc/"7";"8"};
{\ar@{->}@<-2.0ex>@/_0.65pc/"8";"7"};
{\ar@{=}@/^0.65pc/"8";"9"};
{\ar@{->}@<-2.0ex>@/_0.65pc/"9";"8"};
\endxy
\ee
\eprf

\subsection{From algebra to probability theory}
\label{sec:fdfullyfaithful}

We will show that the functors of Theorems
\ref{thm:FinProbtoCP} and \ref{thm:restrictEtofinprob} are fully faithful. 

\bd
\label{defn:fullyfaithful}
A functor $F:\mD\to\mathcal{E}$ from a category $\mD$ to a category $\mathcal{E}$
is \emph{\uline{fully faithful}} iff for every pair of objects $d$ and $d'$ in $\mD,$
the function $F:\mD(d,d')\to\mathcal{E}\big(F(d),F(d')\big)$ is a bijection
(\emph{\uline{full}} means this function is surjective and \emph{\uline{faithful}} 
means it is injective). Here, $\mD(d,d')$ denotes the set of morphisms
from $d$ to $d'$ in the category $\mD.$ 
\ed
In other words, if there is a fully faithful functor $\mD\to\mathcal{E}$ that
is also injective on objects, this allows one to study
$\mD$ in a broader context without altering the relations between its objects. 
Being faithful means that no information is lost when viewing 
$\mD$ as inside of $\mathcal{E}$ and being full means that no extra information is added to $\mD$ 
when viewing $\mD$ inside of $\mathcal{E}.$ In particular, $\mD$ when viewed on its own
has the same set (if it even is a set) of isomorphism classes of objects
as when viewing $\mD$ inside of $\mathcal{E}.$

\bt
\label{thm:fdfullyfaithful}
The functors $C:\FinSetStoch^{\op}\to\fdcCAlgPos$ and 
$C:\FinSetFun^{\op}\to\fdcCAlg$ 
from Theorems \ref{thm:FinProbtoCP} and \ref{thm:restrictEtofinprob}
are fully faithful and injective on objects.
\et

\bprf
Injectivity on objects is immediate from the definition of $C.$ 
Fix any two sets $X$ and $Y$ 
and a CP map $\vf:\C^{Y}\stoch\C^{X}.$
The goal is to show
there exists a unique stochastic map $f:X\stoch Y$ 
such that ${C}^{f}=\vf.$ 
Furthermore, when 
$\vf$ is a $*$-homomorphism, then the associated stochastic map
$f$ is a function. These two claims are proven in two steps.
\begin{enumerate}[i.]
\item
First note that for any $y\in Y,$ $\vf(e_{y})$ is an element of $\C^{X}$
and, because $\{e_{x}\}_{x\in X}$ is a basis of $\C^{X},$ 
there exist unique numbers $\{f_{yx}\in\C\}_{x\in X}$ such that
\be
\vf(e_{y})=\sum_{x\in X}f_{yx}e_{x}. 
\ee
Positivity of $\vf$ demands that 
\be
f_{yx}\ge0
\ee 
for all $x\in X$ (and $y\in Y$). 
Since $\vf$ is linear and unital, 
\be
\hspace{-1mm}
\sum_{x\in X}e_{x}=1_{X}=\vf(1_{Y})=\vf\left(\sum_{y\in Y}e_{y}\right)
=\sum_{y\in Y}\vf(e_{y})
=\sum_{y\in Y}\sum_{x\in X}f_{yx}e_{x}
=\sum_{x\in X}\left(\sum_{y\in Y}f_{yx}\right)e_{x}.
\ee
Since $\{e_{x}\}_{x\in X}$ is a linearly independent set of vectors in $\C^{X},$ 
this implies
\be
\label{eq:sumfprobability1}
\sum_{y\in Y}f_{yx}=1.
\ee
This shows that for every $x\in X,$ the assignment
\be
Y\ni y\mapsto f_{yx}
\ee
is a probability measure on $Y.$ In other words, the numbers
$\{f_{yx}\}_{y\in Y,x\in X}$ define a stochastic map $f:X\stoch Y.$ 
It is readily checked that ${C}^{f}=\vf.$ This proves that
${C}$ is full. 

To see that $C$ is faithful, suppose that
$f':X\stoch Y$ is another stochastic map such that 
${C}^{f}={C}^{f'}.$ In particular, this implies that
for any $y\in Y,$
\be
\xy0;/r.20pc/:
(-14.1421356237,14.1421356237)*+{{C}^{f}(e_{y})}="2";
(-14.1421356237,-14.1421356237)*+{\ds\sum_{x\in X}f_{yx}e_{x}}="1";
(14.1421356237,14.1421356237)*+{{C}^{f'}(e_{y})}="3";
(14.1421356237,-14.1421356237)*+{\ds\sum_{x\in X}f'_{yx}e_{x}}="4";
{\ar@{=}@/^0.95pc/"1";"2"};
{\ar@{=}@/^0.95pc/"2";"3"};
{\ar@{=}@/^0.95pc/"3";"4"};
\endxy
\ee
By linear independence of the set $\{e_{x}\}_{x\in X},$ this implies
\be
f_{yx}=f'_{yx}\qquad\forall\;x\in X.
\ee
Since this calculation was independent of $y,$ this also proves the equality
for all $y\in Y.$ Hence, $f=f'$ proving that ${C}$ is faithful. 
This concludes the proof that ${C}$ is fully faithful.

\item
Using the same notation as in the previous step, it suffices to show that
if $\vf:\C^{Y}\to\C^{X}$ is actually a $*$-homomorphism, then the 
corresponding stochastic map
$f:X\stoch Y$ comes from a function $f:X\to Y.$ 
First notice that for any pair of \emph{distinct} 
elements $y,y'\in Y,$ because $\vf$ is a homomorphism, 
\be
\xy0;/r.35pc/:
(15.6366296494,-12.4697960372)*+{\ds\sum_{x\in X}\sum_{x'\in X}f_{yx}f_{y'x'}\de_{xx'}e_{x'}}="7";
(0,-20)*+{\ds\sum_{x\in X}f_{yx}f_{y'x}e_{x}}="6";
(-15.6366296494,-12.4697960372)*+{0}="5";
(-19.4985582436,4.45041867913)*+{\vf(0)}="4";
(-8.67767478235,18.019377358)*+{\vf(e_{y}e_{y'})}="3";
(8.67767478235,18.019377358)*+{\vf(e_{y})\vf(e_{y'})}="2";
(19.4985582436,4.45041867913)*+{\ds\sum_{x\in X}\sum_{x'\in X}f_{yx}f_{y'x'}e_{x}e_{x'}}="1";
{\ar@{=}@/_0.75pc/"1";"2"};
{\ar@{=}@/_0.75pc/"2";"3"};
{\ar@{=}@/_0.75pc/"3";"4"};
{\ar@{=}@/_0.75pc/"4";"5"};
{\ar@{=}@/_0.75pc/"6";"7"};
{\ar@{=}@/_0.75pc/"7";"1"};
\endxy
\ee
proving that 
\be
\label{eq:fhomomorphismdistinct}
f_{yx}f_{y'x}=0\qquad\forall\;x\in X.
\ee
Putting this result aside for the time being,  
recall the result (\ref{eq:sumfprobability1}) which, together with the fact that 
$f_{yx}\ge0$ for all $y\in Y$ and $x\in X,$ implies that for each $x\in X,$
there exists a $y\in Y$ such that $f_{yx}>0$ (for otherwise, the sum
over all $y$ would not be $1$). For any other $y'\in Y\setminus\{y\},$ 
the result (\ref{eq:fhomomorphismdistinct}) then implies
\be
f_{y'x}=0\qquad\forall\;y'\in Y\setminus\{y\}.
\ee
Using (\ref{eq:sumfprobability1}) again implies 
\be
f_{yx}=1.
\ee
In other words, for every $x\in X,$ there exists a unique $y\in Y$ 
such that $f_{yx}=1$ and $f_{y'x}=0$ for all $y'\in Y\setminus\{y\}.$ 
Set $f(x)$ to be this unique element $y,$ i.e. $f(x):=y.$ 
This assignment then defines a function $f:X\to Y.$
\end{enumerate}
\eprf

It is a fact that every $n$-dimensional commutative $C^*$-algebra is
$*$-isomorphic to $\C^{X}$ for some set $X$ with cardinality $n$
(this will follow from results in Section \ref{sec:GelfandNaimark}). 
This, together with
Theorem \ref{thm:fdfullyfaithful}, says that $\FinSetStoch^{\op}$ is equivalent 
to the category $\fdcCAlgPos$ consisting of finite-dimensional
commutative $C^*$-algebras and positive maps. An explicit inverse will 
be constructed in greater generality in Section \ref{sec:stochasticGelfandspectrumfunctor}.

\subsection{Some quantum mechanics}
\label{sec:examplesinQM}

Completely positive maps represent general processes that
are allowed in quantum mechanics. We will provide
three different examples illustrating the versatility of the 
$C^*$-algebraic approach. 
The first two examples are actually examples of $*$-homomorphisms.
We thank Stefan Andronache and Marcelo Nomura for discussions 
leading to these examples. 

\bx
\label{ex:unitaryevolution}
Let $\Hi$ be a finite-dimensional%
\footnote{We use finite-dimensional Hilbert spaces 
to avoid technicalities involving domains of operators.}
 Hilbert space and $H$ a self-adjoint
operator, thought of as a Hamiltonian. For each $t\in\R,$ the operator
$\exp\left(-\frac{itH}{\hbar}\right),$ where $\hbar$ is the reduced Planck constant,
is unitary and describes unitary time evolution in quantum mechanics. 
Its action on observables is given by the adjoint action
(the Heisenberg picture)
\be
\begin{split}
\mB(\Hi)&\to\mB(\Hi)\\
A&\mapsto \exp\left(\frac{itH}{\hbar}\right)A\exp\left(-\frac{itH}{\hbar}\right).
\end{split}
\ee
Infinitesimally, this gives the differential equation that describes
the time evolution of an observable
\be
\frac{dA}{dt}=\frac{1}{i\hbar}[A,H].
\ee
Now, given an initial state $\tr(\rho\;\cdot\;):\mB(\Hi)\stoch\C$
for some density matrix $\rho$
(see Theorem \ref{thm:statesdensitymatrices} for notation) this pulls back
to a new state
$\tr(\rho\;\cdot\;)\circ\mathrm{Ad}_{\exp\left(-\frac{itH}{\hbar}\right)}$
under this evolution 
\be
\label{eq:adjointevolution}
\xy0;/r.25pc/:
(-15,7.5)*+{\mB(\Hi)}="B";
(15,7.5)*+{\mB(\Hi)}="A";
(0,-7.5)*+{\C}="Complex";
{\ar^{\mathrm{Ad}_{\exp\left(-\frac{itH}{\hbar}\right)}}"B";"A"};
{\ar@{~>}^{\tr(\rho\;\cdot\;)}"A";"Complex"};
{\ar@{~>}_{\tr(\rho\;\cdot\;)\circ\mathrm{Ad}_{\exp\left(-\frac{itH}{\hbar}\right)}}"B";"Complex"};
\endxy
\ee
and is given by the assignment
\be
\begin{split}
\mB(\Hi)&\to\C\\
A&\mapsto\tr\left(\rho\exp\left(\frac{itH}{\hbar}\right)A\exp\left(-\frac{itH}{\hbar}\right)\right)=\tr\left(\exp\left(-\frac{itH}{\hbar}\right)\rho\exp\left(\frac{itH}{\hbar}\right)A\right)
\end{split}
\ee
by cyclicity of the trace. By the uniqueness of density matrices representing
states (Theorem \ref{thm:statesdensitymatrices}), 
this proves that the time evolution of $\rho$ is given by the
differential equation 
\be
\frac{d\rho}{dt}=\frac{1}{i\hbar}[H,\rho].
\ee
This is the quantum Liouville equation. Notice how it has the opposite sign
of the evolution for observables. Also notice that if $\rho$ is a pure state, 
this reduces to Schr\"odinger's equation and reproduces the Schr\"odinger
picture. Thus, the $C^*$-algebraic approach naturally incorporates both
perspectives. It is important to notice that the straight horizontal arrow in 
(\ref{eq:adjointevolution}) indicates \emph{deterministic} evolution! 
The ``non-determinacy'' of this example 
only appears in the states and their expectation values. 
\ex

\bx
\label{ex:measurement}
Measurement is an example of a CP map, in fact a 
$*$-homomorphism.  
We think of measuring a quantum mechanical system 
in terms of numbers. For example, in the Stern-Gerlach experiment,
we measure a particle (such as a silver atom) to be spin up or spin
down depending on its position on a screen after moving through
a magnetic field \cite{Sa93}. A state of a quantum-mechanical system is
therefore reduced to a probability measure on the set
of the eigenvalues of the observable being measured. 

The general (finite-dimensional) situation is as follows.%
\footnote{For infinite systems, one must use projection-valued
measures \cite{NiCh10}.}
Let $\Hi$ be a finite-dimensional Hilbert space, let 
$A$ be a self-adjoint operator on $\Hi,$ and let $\w$ be a state
on $\mB(\Hi).$  
Let $\s(A)$ denote the spectrum of $A$ (in this case,
$\s(A)$ is just the set of eigenvalues of $A$). 
$A$ is an observable and $\s(A)$ are the possible values
that an observer will see when trying to measure the observable $A.$
It is a fact of life that an observer can \emph{only} measure elements
of $\s(A)$ (this is the origin of the ``discreteness'' of quantum mechanics). 
In this situation, what is the probability measure that
an observer expects to see on $\s(A)$? 
The resulting probability measure is obtained by pulling back $\w$ 
along the ``measurement'' map 
\be
\begin{split}
m:\C^{\s(A)}&\to\mB(\Hi)\\
e_{\l}&\mapsto P_{\l},
\end{split}
\ee
where $P_{\l}$ is the projection operator onto the eigenspace associated
with the eigenvalue $\l\in\s(A)$ of $A.$ 
This produces a state on $\C^{\s(A)}$ via pullback
\be
\xy0;/r.25pc/:
(-12.5,7.5)*+{\C^{\s(A)}}="B";
(12.5,7.5)*+{\mB(\Hi)}="A";
(0,-7.5)*+{\C}="Complex";
{\ar^{m}"B";"A"};
{\ar@{~>}_{\w\circ m}"B";"Complex"};
{\ar@{~>}^{\w}"A";"Complex"};
\endxy
\ee
By our analysis from Section \ref{sec:fdfullyfaithful}, $p:=\w\circ m$ 
precisely provides us with a probability measure
on $\s(A).$ This probability measure is interpreted as the
probability distribution of measuring the corresponding eigenvalues 
of $A$ with respect to the state $\w.$ 
Note that this example does not have an evolution of the form 
$\mB(\Hi)\to\mB(\mathcal{K})$ between observables on Hilbert spaces
$\Hi$ and $\mathcal{K}.$  
Instead, it is an ``evolution'' from the quantum phase space $\mB(\Hi)$ 
to the classical one $\s(A).$ 
Hence, the $C^*$-algebraic formulation allows one to treat classical
and quantum systems in the same category. 
\ex

In both of the previous situations, the evolution was described by
a $*$-homomorphism. In the following example, we provide 
completely positive evolution that is not a $*$-homomorphism. 
Such evolutions can be used to model noise and measurement
among other things \cite{NiCh10}.

\bx
\label{ex:bitflipchannel}
Let $\big\{|\!\up\>,|\!\dn\>\big\}$ denote the standard Euclidean basis in $\C^2.$ 
Let $|\!\up\>\<\up\!|$ denote the projection operator onto 
$\mathrm{span}\big(|\!\up\>\big)$ and similarly for $|\!\dn\>.$ 
The states on $\mB(\C^2)$ given by $\tr\big(|\!\up\>\<\up\!|\;\cdot\;\big)$
and  $\tr\big(|\!\dn\>\<\dn\!|\;\cdot\;\big)$ can be interpreted as a qubit 
being either in the spin up or spin down state with respect to an 
appropriately chosen coordinate frame. Consider an operation on the 
space of all density matrices whose action on these two states is given by 
\be
\label{eq:bitflipchannel}
\begin{split}
|\!\up\>\<\up\!|&\mapsto p |\!\up\>\<\up\!|+(1-p)|\!\dn\>\<\dn\!|\\
|\!\dn\>\<\dn\!|&\mapsto (1-p)|\!\up\>\<\up\!|+p|\!\dn\>\<\dn\!|
\end{split}
\ee
for some fixed $p\in[0,1].$ This is to be interpreted as an operation that 
has probability $1-p$ of flipping a spin up state to a spin down state
and has probability $1-p$ of flipping a spin down state to a spin up state. 
Does there exist a completely positive map $f:\mB(\C^2)\stoch\mB(\C^2)$ 
that pulls back the pure states above on the left 
in (\ref{eq:bitflipchannel}) to the states on the right in (\ref{eq:bitflipchannel})? 
The answer is yes. Let 
\be
E:=\sqrt{p}\begin{bmatrix}1&0\\0&1\end{bmatrix}
\aand
F:=\sqrt{1-p}\begin{bmatrix}0&1\\1&0\end{bmatrix}
\ee
and define $f$ to be given by 
\be
\mB(\C^2)\ni A\mapsto E^{*}AE+F^{*}AF. 
\ee
One can check that $f$ is completely positive and is not a 
$*$-homomorphism. The induced action on the density matrices
can be obtained by looking at the pullback of an initial density matrix $\rho$ 
\be
\mB(\C^2)\ni A\mapsto\tr\big(\rho f(A)\big)=\tr(\rho E^{*}AE+\rho F^{*}AF)
=\tr\big((E\rho E^{*}+F\rho F^*)A\big).
\ee
Hence, the resulting density matrix under such an operation is given by
$E\rho E^{*}+F\rho F^*.$ One can check that 
(\ref{eq:bitflipchannel}) holds under this transformation. 
So again, the $C^*$-algebraic framework
is able to reproduce familiar operations such as this one, known as the 
bit flip channel \cite{NiCh10}, all in the same category of
physical processes. 
\ex

A generalization of Example \ref{ex:unitaryevolution}
to CP maps gives rise to Lindblad's equation describing
\emph{non-deterministic} dynamics in (open) quantum systems
\cite{Li76}.

\section{An equivalence between spaces and algebraic structures}
\label{sec:GelfandNaimark}

On their own, the study of topological spaces and continuous maps
is quite different from the study of algebra and homomorphisms.
Topological spaces consist
of sets equipped with a subset of the power set (known as open sets) 
and continuous maps
are functions whose pre-image function takes open sets to open sets. 
Algebras are sets equipped with binary operations satisfying certain
laws such as associativity and maps between algebras are those 
that respect this algebraic structure. 
In the plain context of set theory, it simply does not make sense to compare
the two areas of mathematics. However, viewing these two areas 
\emph{as a whole}, which is achieved by viewing them as categories, 
it makes sense to compare them. 
In fact, under additional suitable restrictions on the respective sides, 
the two categories are \emph{equivalent}. 
There are several other versions of such equivalences between
topological categories and algebraic ones \cite{Ka41}, \cite{Ka51}. 
We will discuss one such equivalence, the equivalence
between the category of commutative (unital) $C^*$-algebras and $*$-homomorphisms
with the category of compact Hausdorff topological spaces and continuous
maps following the expositions in Folland  \cite{Fo94}, Rudin \cite{Ru91},
and Tao's online notes \cite{Ta09_245B11}.
The references already have well-written proofs of the main results.
Our purpose here is to frame these results in the appropriate
categorical context and establish certain tools that will be used in 
Section \ref{sec:stochasticGelfandNaimark} when we relax the class
of morphisms to include non-deterministic processes.

\subsection{From spaces to algebras}
\label{sec:spacestoalgebras}

In Section \ref{sec:fromfiniteprobtoalg}, we showed how given a finite set $X,$
the set of functions $\C^{X}$ can be made into a $C^*$-algebra
in a natural manner. This can also be done if $X$ is replaced
with a compact Hausdorff space and the function space is chosen
to be the subspace of \emph{continuous} functions from $X$ to $\C.$ 

\bx
\label{ex:ConcHaus}
Let $X$ be a compact Hausdorff space. 
Then the space of complex-valued continuous 
functions, denoted by $C(X),$ with the 
algebraic structure from $\C^{X}$ 
and with norm given by the uniform norm, i.e. 
\be
C(X)\ni \vf\mapsto \lVert \vf \rVert := \sup_{x\in X} |\vf(x)|,
\ee
is a commutative $C^*$-algebra. 
Cauchy completeness follows from the continuous (uniform) limit theorem
(see Theorem 7.12 in \cite{Ru76} or Theorem 21.6 in \cite{Mu00}). 
The constant function $1_{X}$ is the unit for $C(X).$ 
Note that although compactness is needed for $C(X)$ to be a 
$C^*$-algebra with this norm, it is not necessary that $X$ be Hausdorff.
Proving that $C(X)$ is a $C^*$-algebra under the only assumption that
$X$ is compact is a good exercise.%
\footnote{Sketch of proof: start with a Cauchy sequence
in $C(X),$ use Cauchy completeness in $\C$ to obtain a candidate function
pointwise, 
show that the resulting function is bounded, then show that the sequence
of functions converges to it uniformly, and finally apply the uniform limit theorem
to guarantee it is continuous.}
If $X$ is not Hausdorff, there exist points that cannot be separated
by continuous complex valued functions so that continuous functions
cannot distinguish them (see Remark \ref{rmk:whyHausdorff}
and Theorem \ref{thm:cHStoch} for further details).
\ex

\bn
\label{prop:Conmorphisms}
Let $X$ and $Y$ be compact Hausdorff spaces and let 
$f:X\to Y$ be a continuous function. The function $C^{f}:C(Y)\to C(X)$
defined by 
\be
\label{eq:Conmorphisms}
C(Y)\ni\b\mapsto C^{f}(\vf):=\vf\circ f
\ee
is a $*$-homomorphism. 
\en

\bprf
First note that $C^{f}(\vf)$ is continuous since it is the composition
of continuous functions. 
To see that $C^{f}$ respects the 
product, let $\psi,\vf\in C(Y)$ and let $x\in X.$ Then
\be
\begin{split}
\big(C^{f}(\psi\vf)\big)(x)&=(\psi\vf)\big(f(x)\big)\qquad\text{by def'n of $C^f$}\\
&=\Big(\psi\big(f(x)\big)\Big)\Big(\vf\big(f(x)\big)\Big)\qquad\text{by def'n of $\psi\vf$}\\
&=\big((\psi\circ f)(x)\big)\big((\vf\circ f)(x)\big)\qquad\text{by def'n of $\circ$}\\
&=\big((\psi\circ f)(\vf\circ f)\big)(x)\qquad\text{by def'n of $(\psi\circ f)(\vf\circ f)$}\\
&=\big(C^{f}(\psi)C^{f}(\vf)\big)(x)\qquad\text{by def'n of $C^f$}
\end{split}
\ee
proving that $C^{f}(\psi\vf)=C^{f}(\psi)C^{f}(\vf).$ Linearity, $*$-preservation,
and preservation of units follow from analogous calculations.
\eprf

\begin{notation}
\label{not:cHaus}
Let $\cH$ denote the category whose objects are compact Hausdorff
spaces and whose morphisms are continuous maps. 
\end{notation}

\bn
\label{prop:cHaustocCAlg}
The assignment
\be
\label{eq:cHaustocCAlg}
\begin{split}
\cH^{\op}&\xrightarrow{C}\cCAlg\\
X&\mapsto C(X)\\
\Big(X\xrightarrow{f}Y\Big)&\mapsto
\Big(C(X)\xleftarrow{C^{f}}C(Y)\Big)
\end{split}
\ee
from Example \ref{ex:ConcHaus} and Proposition \ref{prop:Conmorphisms}
defines a functor.
\en

\bprf
From the definition of $C,$ it is manifest that $C$ is a functor:
the identity $\id_{X}:X\to X$ gets sent to $C^{\id_{X}}=\id_{C(X)}$
and the composition $X\xrightarrow{f}Y\xrightarrow{g}Z$ gets sent 
to the function
\be
C(Z)\ni\vf\mapsto C^{g\circ f}(\vf)=\vf\circ(g\circ f)=(\vf\circ g)\circ f
=C^{f}\big(\vf\circ g\big)
=C^{f}\big(C^{g}(\vf)\big)
=(C^{f}\circ C^{g})(\vf)
\ee
since composition of functions is associative. 
\eprf

\subsection{Some topological preliminaries}
\label{sec:topologicalpreliminaries}

Given a Banach algebra $\mA,$ a certain subspace, denoted by $\mA^{*},$ of
the \emph{linear} dual space $\mA^{\vee}$
has a natural topology (the operator norm topology) 
coming from the norms on $\mA$ and $\C.$
However, this topology has too many open sets for some purposes. 
Another topology on $\mA^{*}$ that has fewer open sets is the weak* topology.
Fewer open sets allows more sets to be compact and more compact
sets means that sequences and nets have a better chance of
converging (see Tao's notes for an illuminating discussion \cite{Ta09_245B11}). 
There are several ways to describe this topology. The simplest for our purposes 
is in terms of convergent nets, though it will
also be useful to describe it in terms of topological
vector spaces whose topology is induced by a family of seminorms.
This latter perspective will also become more important in Section
\ref{sec:stochasticGelfandNaimark}.

\bd
A \emph{\uline{seminorm}} on a
vector space $V$ is a function 
$\lVert\ \cdot \ \rVert:V\to\R$ satisfying
\begin{enumerate}[i.]
\setlength{\itemsep}{0pt}
\item
$\lVert v\rVert\ge0$ for all $v\in V,$
\item
$\lVert \l v\rVert=|\l|\lVert v\rVert$ for all 
$\l\in\C,v\in V,$ and
\item
$\lVert v+w\rVert\le\lVert v\rVert+\lVert w\rVert$
for all $v,w\in V.$ 
\end{enumerate}
\ed

The only difference between a norm and a seminorm is that
a norm satisfies the additional condition that $\lVert v\rVert=0$
implies $v=0.$ 

\bx
Let $V$ be a vector space and $(W,\lVert\;\cdot\;\rVert)$ be a normed
vector space. For any linear map $\vf:V\to W,$ the function
$\Vert\;\cdot\;\rVert_{\vf}:V\to\R$ defined by sending $v\in V$
to $\big\lVert\vf(v)\big\rVert$ is a seminorm on $V.$ 
A particular case of interest is when $W:=\C.$ 
\ex

A single norm provides a natural topology via
open balls of varying radii at different points. 
While a single seminorm also provides a natural topology
in a similar way, it is rarely Hausdorff and does not enjoy
many desirable properties (see Example \ref{ex:openballs}). 
However, it does become a topological
vector space. 

\bd
A \emph{\uline{topological vector space}} consists of a vector
space $V$ equipped with a topology in which 
the addition of vectors $+:V\times V\to V$ and 
scalar multiplication $\C\times V\to V$ are both continuous. 
\ed

\bd
\label{defn:basefortopology}
A \emph{\uline{base}} on a set $X$ is
a collection $\mathfrak{B}\subseteq\mathcal{P}(X)$ of
subsets of $X$ satisfying
\begin{enumerate}[i.]
\setlength{\itemsep}{0pt}
\item
for each $x\in X,$ there exists a
$B\in\mathfrak{B}$ such that $x\in B,$
\item
for every pair $B_{1},B_{2}\in\mathfrak{B}$ and
every $x\in B_{1}\cap B_{2},$ there exists
a $B\in\mathfrak{B}$ such that 
$x\in B\subseteq B_{1}\cap B_{2}.$
\end{enumerate}
A base $\mathfrak{B}$ as above
can be used to construct a topology
on $X$ by setting $\t_{\mathfrak{B}}$ to be the smallest
topology containing $\mathfrak{B}.$ Explicitly, its elements
are given by unions of all possible elements in $\mathfrak{B},$
\be
\t_{\mathfrak{B}}=\left\{ \bigcup_{B\in\mathfrak{C}} B \;:\; \mathfrak{C}\subseteq\mathfrak{B}\right\}.
\ee
This topology is called the
\emph{\uline{topology generated}} by $\mathfrak{B}.$
\ed

\bx
\label{ex:openballs}
Let $\lVert\;\cdot\;\rVert$ be a seminorm on a vector space
$V.$ Then the set
\be
\Big\{B_{r}(v)\subseteq V\;:\;v\in V\;\text{ and }\; r>0\Big\}, 
\ee
where
\be
B_{r}(v):=\big\{w\in V\;:\;\lVert v-w\rVert<r\big\},
\ee
forms a base on $V.$ With respect to the topology generated by this base, $V$ is a 
topological vector space. Furthermore, it is Hausdorff if and only if
$\lVert\;\cdot\;\rVert$ is a norm (see Tao \cite{Ta09_245B11}).
$B_{r}(v)$ is called the \emph{\uline{open ball}} of radius $r\ge0$ 
centered at $v\in V.$
\ex

A sufficiently robust 
\emph{family} of seminorms typically has a more manageable topology
that naturally arises frequently enough to merit study. 

\begin{notation}
\label{not:nset}
For $n\in\N,$ set $\ov n:=\{1,\dots,n\}.$ Let $X$ be a set and 
$\mathfrak{S}\subseteq\mathcal{P}(X)$ a subset of the power set.%
\footnote{$\mathfrak{S}$ is a Fraktur ``S.''}
The value of a function $S:\ov n\to\mathfrak{S}$ at $i\in\ov n$ is
denoted by $S_{i}$ instead of $S(i).$ 
\end{notation}

\bd
\label{defn:subbase}
A \emph{\uline{subbase}} on a set $X$ is
a collection $\mathfrak{S}\subseteq\mathcal{P}(X)$ of
subsets of $X$ such that the set
\be
\mathfrak{B}:=
\bigcup_{n\in\{0\}\cup\N}
\left\{\bigcap_{i=1}^{n}S_{i}\;:\;\ov n\xrightarrow{S}\mathfrak{S}\right\}
\ee
of all finite intersections of elements in $\mathfrak{S}$ is a base
for $X.$ Here, $\ds\bigcap_{i=1}^{0}S_{i}$ is the empty intersection,
which is taken to be $X$ itself. 
The \emph{\uline{topology generated}} by $\mathfrak{S}$ 
is the topology generated by $\mathfrak{B}.$ 
\ed

It is a fact that any collection of
subsets of a set forms a subbase for some topology on that set 
and that topology is the smallest topology containing those subsets 
(see Section 5 of \cite{Wi04}).

\bn
\label{prop:familyoftopologies}
Let $V$ be a vector space with $\{\lVert\;\cdot\;\rVert_{\a}\}_{\a\in A}$ a
family of seminorms on $V$ indexed by some set $A.$ 
Let $\t_{\a}$ denote the topology 
associated to $\lVert\;\cdot\;\rVert_{\a}.$ 
Then $\ds\mathfrak{S}:=\bigcup_{\a\in A}\t_{\a}$
is a subbase on $V$ and 
$V$ with the topology generated by $\mathfrak{S}$ 
is a topological vector space. Furthermore, $V$ is 
Hausdorff if and only if for each $v\in V\setminus\{0\},$ there
exists an $\a\in A$ such that $\lVert v\rVert_{\a}\ne 0.$
Such a family of seminorms is said to \emph{\uline{separate points}} of $V$ 
or is said to be a \emph{\uline{separating family}}. 
\en

\bprf
See Theorem 5.14 and Proposition 5.16 in Folland \cite{Fo07}.
\eprf

\br
The use of seminorms is not necessary in the conclusion that $V$
is a topological vector space. All that is required
is that each of the topologies $\t_{\a}$ provide $V$ with the structure
of a topological vector space \cite{Ta09_245B11}. 
\er

Let $V$ be a vector space and let $V^{\vee}$ denote its linear dual
(linear functions from $V$ to $\C$)
In many applications where $V$ has a norm, 
$V^{\vee}$ is too large to deal with. Instead, one uses the norm
to restrict to a certain subspace of linear functionals. 

\bn
\label{prop:topologicaldual}
Let $V$ be a normed vector space with norm $\lVert\;\cdot\;\rVert.$
For $\chi\in V^{\vee},$ set (abusing notation slightly)
\be
\lVert\chi\rVert:=\inf\Big\{M\ge0\;:\;\big|\chi(v)\big|\le M\lVert v\rVert\quad\forall\;v\in V\Big\}, 
\ee 
which can be infinite. 
Set $V^*\subseteq V^{\vee}$ to be the subset given by 
\be
V^*:=\big\{\chi\in V^{\vee}\;:\;\lVert\chi\rVert<\infty\big\}. 
\ee
Then $\lVert\;\cdot\;\rVert$ is a norm on $V^*$ and $V^*$ equipped with 
this norm is a Banach space, known
as the \emph{\uline{topological dual space}} of $V.$ 
Elements of $V^*$ are called 
\emph{\uline{bounded/continuous linear functionals}} on $V.$ 
\en

\bprf
This is a general fact about operator norms (see Proposition 5.4 in \cite{Fo07}). 
\eprf

\bn
\label{prop:seminormsondualfromoriginal}
Let $V$ be a vector space. For each $v\in V,$ define 
$\lVert\;\cdot\;\rVert_{v}:V^{\vee}\to\C$ by
\be
V^{\vee}\ni\chi\mapsto\lVert\chi\rVert_{v} :=|\chi(v)|.
\ee
Then $\lVert\;\cdot\;\rVert_{v}$ is a seminorm on $V^{\vee}.$ 
Furthermore, for each $\chi\in V^{\vee}\setminus\{0\},$ there exists a $v\in V$
such that $\chi(v)\ne0,$ i.e. the family of seminorms 
$\{\lVert\;\cdot\;\rVert_{v}\}_{v\in V}$ is separating. 
The same statement holds for $V^*\subseteq V^{\vee}$ with the same family.
\en

\bprf
This immediately follows from the definitions. 
\eprf

\bd
Let $V$ be a vector space and let 
$\{\lVert\;\cdot\;\rVert_{v}\}_{v\in V}$ be the family
of seminorms as in Proposition \ref{prop:seminormsondualfromoriginal}.
The Hausdorff topology generated by this family
as described in Proposition \ref{prop:familyoftopologies}
is known as the \emph{\uline{weak* topology}} on $V^{*}.$
\ed

\bn
\label{prop:weakstartopologybase}
Let $V$ be a normed vector space.
For each $\chi\in V^{*}, \e>0, n\in\N,$ and
$\a_{n}:\ov n\to V,$ set
\be
B_{\e}^{\a_{n}}(\chi)
:=\bigcap_{j\in\ov n}\Big\{\xi\in V^{*} \;:\;
\lVert\chi-\xi\rVert_{\a_{n}(j)}<\e\Big\}.
\ee
Then 
\be
\mathfrak{B}:=\Big\{ B_{\e}^{\a_{n}}(\chi)\subseteq V^{\vee} \;:\;
\e>0,\;n\in\N,\;\a_{n}:\ov n\to V,\;\chi\in V^{*}\Big\}
\ee
is a base for the weak* topology on $V^{*}.$
\en

\bprf
This is a general fact for the topology generated by a
separating family of seminorms---see Section 2.4 in \cite{Pe89}.
\eprf

In other words, a base for the weak* topology consists of
finite intersections of balls coming from the different seminorms. 

\bn
\label{prop:weak*topologyconvergence}
Let $V$ be a vector space. A net%
\footnote{A net is a function whose domain is a directed set 
(see Section 11 of \cite{Wi04}).}
$\chi:\Theta\to V^{*}$
in $V^{*}$ converges to an element $\lim\chi\in V^{*}$ 
in the weak* topology if and only if 
\be
\lim_{\theta\in\Theta}\big(\chi_{\q}(v)\big)=(\lim\chi)(v)\qquad\forall\;v\in V.
\ee
\en

\bprf
See Tao \cite{Ta09_245B11} or Section 5.4 of Folland \cite{Fo07}.
\eprf

Nets are needed here because not all compact Hausdorff spaces
are separable (see Section \ref{sec:probmeasoncHaus}). 

\bn
Let $V$ be a vector space and let $V^{*}$ denote its topological dual space.
The weak* topology on $V^{*}$ is the weakest topology (meaning that it has
the fewest open sets) such that the functions 
\be
\begin{split}
V^{*}&\xrightarrow{\mathrm{ev}_{v}}\C\\
\chi&\mapsto
\chi(v)
\end{split}
\ee
are continuous for all $v\in V.$
\en

\bprf
See Section 3.14 of \cite{Ru91}. 
\eprf

\bt
[Banach-Alaoglu Theorem]
Let $V$ be a Banach space and let $V^{*}$ denote its topological dual space.
The closed unit ball (with respect to the norm on $V^*$)
in $V^{*}$ is compact with respect to the weak* topology. 
\et

\bprf
See Theorem 5.18 in \cite{Fo07}. 
\eprf

This theorem is surprising because the closed (and bounded) unit ball is compact
with respect to the norm topology if and only if the normed vector space 
is finite-dimensional \cite{Ta09_245B11}.

\subsection{The spectrum of a commutative $C^*$-algebra}

It is perhaps not so surprising that one can obtain a 
commutative $C^*$-algebra from a space as in Section \ref{sec:spacestoalgebras}.
What is more surprising is that there is a way to go back
from $C^*$-algebras to topological spaces.
Although $C^*$-algebras have a norm (and hence a topology), 
keep in mind that a $C^*$-algebra is a \emph{linear} space so 
obtaining a non-trivial topology from such an object will require some work. 

\bd
Let $\mathcal{A}$ be a commutative Banach algebra. 
A \emph{\uline{character}} 
on $\mathcal{A}$ is a continuous nonzero homomorphism 
$\chi : \mathcal{A} \to \C$ of Banach algebras, i.e. 
$\chi$ is linear, satisfies $\chi(ab)=\chi(a)\chi(b)$ for all $a,b\in\mA,$
and satisfies the condition that there exists an $a\in\mA$ 
for which $\chi(a)\ne0.$
\ed

Characters form a subset of $\mA^{*},$
the topological dual space of $\mA.$ This dual space
has a natural topology coming from
the norms on both $\mA$ and $\C,$ 
but the topology of interest for us is the
weak* topology on $\mA^{*}.$

\bd
\label{def:spectrum}
Let $\mathcal{A}$ be a commutative Banach algebra. The 
\emph{\uline{spectrum of $\mathcal{A}$}} is the set 
\be
\s(\mathcal{A}):=
\left\{\mathcal{A}\xrightarrow{\chi}\C \;:\;\chi\text{ is a character} \right \}. 
\ee
equipped with the subspace topology coming from $\mathcal{A}^{*}$  
via the weak* topology.
\ed

\br
The relationship between the spectrum of a commutative $C^*$-algebra
and the spectrum of a certain operator on a Hilbert space is 
described in Proposition 1.15 in \cite{Fo94}. Briefly, if $A$ is bounded and 
normal, then the operator-theoretic spectrum of $A$ is canonically 
homeomorphic with the 
spectrum of the commutative unital $*$-subalgebra generated by $A.$
This is the closure of the set of polynomials in $A$ and $A^*.$ By definition 
of being normal, $A$ and $A^*$ commute so that the order does not matter
in which the $A$'s and $A^*$'s appear in such polynomials guaranteeing that
this $C^*$-algebra is commutative. 
\er

\bn
\label{prop:Folland1.10}
Let $\mA$ be a Banach algebra and let $\chi\in\s(\mA)$ be a character.
Then
\begin{enumerate}[i.]
\setlength{\itemsep}{0pt}
\item
$\chi(1_{\mA})=1,$
\item
$\chi(a)\ne0$ for all invertible elements $a\in\mA,$ and
\item
$|\chi(a)|\le\lVert a\rVert$ for all $a\in\mA.$ 
\end{enumerate}
\en

\bprf
See Proposition 1.10 in \cite{Fo94}.
\eprf

\bn
Let $\mathcal{A}$ be a commutative Banach algebra. 
Then the spectrum
$\s (\mathcal{A})$ is a nonempty compact Hausdorff space. 
\en

\bprf
See Theorem 11.9 in \cite{Ru91}.
$\s(\mA)$ is nonempty because $\mA$ is unital. 
\eprf

\bn
\label{prop:spectrumofmap}
Let $\mB\xrightarrow{f}\mA$ be $*$-homomorphism of 
commutative Banach algebras. Then the function
\be
\label{eq:spectrumofmorphism}
\begin{split}
\s(\mA)&\xrightarrow{\s^{f}}\s(\mB)\\
\chi&\mapsto\chi\circ f
\end{split}
\ee
is continuous. 
\en

\bprf
First, note that since
\be
(\chi\circ f)(xy)=\chi\big(f(x)f(y)\big)
=\chi\big(f(x)\big)\chi\big(f(y)\big)
=(\chi\circ f)(x)(\chi\circ f)(y)
\ee
for all $x,y\in\mB,$ $\chi\circ f$ is indeed an element of $\s(\mB).$ 
We will give two proofs of continuity, 
one using open sets and the other using nets. 
\begin{enumerate}[i.]
\item
It suffices to show that the inverse
image of a base element gets sent
to an open set. By Proposition \ref{prop:weakstartopologybase}
and by definition of the subspace topology, a base for the
topology on $\s(\mB)$ is given by
\be
\s(\mB)\cap B_{\e}^{\a_{n}}(\xi)
=\left\{\chi\in\s(\mB)\;:\; \lVert\chi-\xi\rVert_{\a_{n}(j)}<\e\quad\forall\;
j\in\ov n\right\}
\ee
over all $\xi\in\mB^{*},\a_{n}:\ov n\to\mB,\e>0,n\in\N.$
$\s^{f}$ is continuous if%
\footnote{
Let $f:X\to X'$ be a function
and let $\mathfrak{B}$ and $\mathfrak{B}'$
be bases for topologies on $X$ and
$X',$ respectively. $f$ is continuous
if for every $x\in X,$
and for any $B'\in\mathfrak{B}'$ with
$f(x)\in B',$
there exists a $B\in\mathfrak{B}$ with $x\in B$ such that
$f(B)\subseteq B'$ (this is a simple exercise in the definitions). 
}
for every $\chi\in\s(\mA)$ and every
basic set%
\footnote{Technically, we should have said that
for every $\chi\in\s(\mA)$ and every
basic set $B$ containing $\chi\circ f,$
there exists an open set $U$ containing
$\chi$ whose image is contained in $B.$
However, due to the seminorms, it suffices
to take the basic set $B$ containing
$\chi\circ f$ to be centered at $\chi\circ f$ 
by choosing it to be sufficiently small. 
}
$\s(\mB)\cap B_{\e}^{\a_{n}}(\chi\circ f)$
there exists an open set containing
$\chi$ whose image is contained in this set.
In fact, the basic set $B_{\e}^{f\circ\a_{n}}(\chi)$
accomplishes this goal because
\be
\lVert\chi-\xi\rVert_{f\circ\a_{n}(j)}
=\left|\chi\Big(f\big(\a_{n}(j)\big)\Big)
-\xi\Big(f\big(\a_{n}(j)\big)\Big)\right|
=\lVert\chi\circ f-\xi\circ f\rVert_{\a_{n}(j)},
\ee
which shows that 
\be
\begin{split}
\s^{f}\Big(\s(\mA)\cap B_{\e}^{f\circ\a_{n}}(\chi)\Big)
&=\Big\{\xi\circ f\in\s(\mB)\;:\;
\lVert\chi-\xi\rVert_{f\circ\a_{n}(j)}<\e\quad\forall\;j\in\ov n\Big\}\\
&=\Big\{\xi\circ f\in\s(\mB)\;:\;
\lVert\chi\circ f-\xi\circ f\rVert_{\a_{n}(j)}<\e\quad\forall\;j\in\ov n\Big\}\\
&\subseteq\s(\mB)\cap B_{\e}^{\a_{n}}(\chi\circ f).
\end{split}
\ee
Thus, $\s^{f}$ is continuous.
\item
Let $\chi:\Theta\to\s(\mA)$ be a net converging to $\lim\chi$
(in the weak* topology).  
The goal is to show that the net $\s^{f}\circ\chi:\Theta\to\s(\mB)$
converges to $\s^{f}(\lim\chi)\equiv(\lim\chi)\circ f.$ 
For every $b\in\mB,$
\be
\lim_{\q\in\Theta}\big(\s^{f}(\chi_{\q})(b)\big)=\lim_{\q\in\Theta}\Big(\chi_{\q}\big(f(b)\big)\Big)
=(\lim\chi)\big(f(b)\big)=\big(\s^{f}(\lim\chi)\big)(b),
\ee
which establishes the required weak* convergence and thus
shows that $\s^{f}$ is continuous. 
\end{enumerate}
\eprf

As the reader may have noticed, the proof was much shorter using nets. 
The technique of using nets to prove continuity of various functions 
will be heavily used in Section \ref{sec:stochasticGelfandNaimark}. 

\bn
\label{prop:cCAlgtocHaus}
The assignment
\be
\begin{split}
\cCAlg^{\op}&\xrightarrow{\s}\cH\\
\mA&\mapsto\s(\mA)\\
\Big(\mB\xrightarrow{f}\mA\Big)&\mapsto
\Big(\s(\mB)\xleftarrow{\;\s^{f}}\s(\mA)\Big)
\end{split}
\ee
from Definition \ref{def:spectrum} and
Proposition \ref{prop:spectrumofmap}
defines a functor.
\en

\bprf
The proof is completely analogous to the proof of
Proposition \ref{prop:cHaustocCAlg}.
\eprf

\subsection{The Gelfand transform}

Up to this point, we have constructed functors 
\be
\xymatrix{
\cCAlg^{\op}
\ar@<1.0ex>[rr]^(0.525){\s} 
& & \ar@<1.0ex>[ll]^(0.475){C} \cH
}
.
\ee
In general, the diagrams 
\be
\xy0;/r.25pc/:
(-15,7.5)*+{\cH}="3";
(15,7.5)*+{\cH}="1";
(0,-7.5)*+{\cCAlg^{\op}}="2";
{\ar"1";"2"^{C}};
{\ar"2";"3"^{\s}};
{\ar"1";"3"_{\id}};
\endxy
\aand
\xy0;/r.25pc/:
(-15,7.5)*+{\cCAlg}="3";
(15,7.5)*+{\cCAlg}="1";
(0,-7.5)*+{\cH^{\op}}="2";
{\ar"1";"2"^{\s}};
{\ar"2";"3"^{C}};
{\ar"1";"3"_{\id}};
\endxy
\ee
do \emph{not} commute so that $\s$ and $C$ are not inverses
of each other. However, they are \emph{close}.
In the present section, we will construct natural  
isomorphisms
\be
\xy0;/r.25pc/:
(-15,7.5)*+{\cH}="3";
(15,7.5)*+{\cH}="1";
(0,-7.5)*+{\cCAlg^{\op}}="2";
{\ar"1";"2"^{C}};
{\ar"2";"3"^{\s}};
{\ar"1";"3"_{\id}};
{\ar@{=>}(0,6.5);"2"^(0.37){h}};
\endxy
\aand
\xy0;/r.25pc/:
(-15,7.5)*+{\cCAlg}="3";
(15,7.5)*+{\cCAlg}="1";
(0,-7.5)*+{\cH^{\op}}="2";
{\ar"1";"2"^{\s}};
{\ar"2";"3"^{C}};
{\ar"1";"3"_{\id}};
{\ar@{=>}(0,6.5);"2"^(0.37){\G}};
\endxy
\ee
indicating a precise sense in which $\s$ and $C$ are close to being inverses
of each other. 
A natural isomorphism between functors is analogous to a homotopy
between continuous functions or a unitary intertwiner between representations 
and what it amounts to more precisely will
be described presently. This will show that the categories 
$\cCAlg^{\op}$ and $\cH$ are \emph{equivalent}.

\begin{notation}
\label{def:Gelfand}
Let $\mathcal{A}$ be a 
commutative Banach algebra and set 
\be
\label{eq:Gelfandtransform}
\begin{split}
\mathcal{A}&\xrightarrow{\G_{\mathcal{A}}} C\big(\s(\mathcal{A})\big) \\
a &\mapsto\Big(\chi\mapsto\big(\G_{\mathcal{A}}(a)\big)(\chi)
:=\chi(a)\Big).
\end{split}
\ee
$\G_{\mathcal{A}}$ is called the 
\emph{\uline{Gelfand transform on $\mathcal{A}$}}.
\end{notation}

\bn
\label{prop:Gelfandfacts}
The following facts are true regarding the Gelfand transform.
\begin{enumerate}[i.]
\setlength{\itemsep}{0pt}
\item
For every commutative 
Banach algebra $\mA,$
the Gelfand transform $\G_{\mA}:\mA\to C_{c}(\s(\mA))$
is a 
homomorphism of Banach algebras (preserving units).
\item
If $\mA$ is a commutative $C^*$-algebra,
then $\G_{\mA}$ is a $*$-isomorphism (in particular, it is isometric).
\end{enumerate}
\en

\bprf
See
Theorem 1.13 parts (a) and (d)
and Theorem 1.20 in \cite{Fo94}
for the first and second claims, respectively. 
\eprf

\br
If $\mA$ is just an involutive Banach algebra, then 
$\G_{\mA}(a^*)$ is not always equal to 
$\big(\G_{\mA}(a)\big)^*$ for all $a\in\mA.$ Hence, $\G_{\mA}$
is not necessarily a $*$-homomorphism
(see Proposition 1.14 in \cite{Fo94}). 
\er

\bn
\label{prop:Gelfandtransform}
The assignment
\be
\begin{split}
\cCAlg^{\op}_{0}&\xrightarrow{\G}\cCAlg^{\op}_{1}\\
\mA&\mapsto\Big(\mA
\xrightarrow{\G_{\mA}}C\big(\s(\mA)\big)\Big)
\end{split}
\ee
from Definition \ref{def:Gelfand}
defines a natural isomorphism 
$\G:\id_{\cCAlg}\Rightarrow C\circ\s.$ 
\en

\bprf
Associated to any $*$-homomorphism $\mB\xrightarrow{f}\mA$
of $C^*$-algebras is the diagram
\be
\xy0;/r.25pc/:
(-12.5,7.5)*+{\mB}="A1";
(-12.5,-7.5)*+{\mA}="A2";
(12.5,7.5)*+{C\big(\s(\mB)\big)}="CsA1";
(12.5,-7.5)*+{C\big(\s(\mA)\big)}="CsA2";
{\ar"A1";"A2"_{f}};
{\ar"CsA1";"CsA2"^{C^{\s(f)}}};
{\ar"A1";"CsA1"^(0.35){\G_{\mB}}};
{\ar"A2";"CsA2"_(0.35){\G_{\mA}}};
\endxy
.
\ee
To see that it commutes, let $b\in\mB.$
Going along the top and right arrows gives
$
C^{\s(f)}\big(\G_{\mA}(b)\big)
=\G_{\mA}(b)\circ\s^{f}
$
by definition of $C$ as a functor
(see Proposition \ref{prop:cHaustocCAlg}).
Going along the left and bottom arrows
gives
$
\G_{\mA}\big(f(b)\big).
$
To see that these two elements of
$C\big(\s(\mA)\big)$ 
are equal, let $\chi\in\s(\mA).$
Then
\be
\xy0;/r.20pc/:
(11.7557,16.1803)*+{\qquad\qquad
\G_{\mA}\big(f(b)\big)\big(\chi\big)}="1";
(19.0211,-6.1803)*+{\chi\big(f(b)\big)}="2";
(0,-20)*+{(\chi\circ f)(b)}="3";
(-19.0211,-6.1803)*+{\G_{\mB}(b)\big(\chi\circ f\big)}="4";
(-11.7557,16.1803)*+{
\G_{\mB}(b)\Big(\s^{f}\big(\chi\big)\Big)
\qquad\qquad}="5";
{\ar@{=}@/^0.7pc/"1";"2"^(0.55){(\ref{eq:Gelfandtransform})}};
{\ar@{=}@/^0.7pc/"2";"3"};
{\ar@{=}@/^0.7pc/"3";"4"^(0.70){(\ref{eq:Gelfandtransform})}};
{\ar@{=}@/^0.7pc/"4";"5"^(0.45){(\ref{eq:spectrumofmorphism})}};
\endxy
,
\ee
which proves that $\G$ is a natural transformation. It is a natural 
{\emph{iso}morphism} by ii. of Proposition \ref{prop:Gelfandfacts}.
\eprf

\bn
\label{prop:evhomeo}
Let $X$ be a compact Hausdorff space.
The assignment
\be
\label{eq:evaluation}
\begin{split}
h_{X} : X &\to \s\big(C(X)\big) \\
x &\mapsto\Big(C(X)\ni f\mapsto\big(h_{X}(x)\big)\big(f\big)
:=f(x)\Big)
\end{split}
\ee
is a well-defined homeomorphism.
\en

\bprf
See Theorem 1.16 of \cite{Fo94}.
and Example 2 after
Theorem 1.30 in \cite{Fo94}.
\eprf

\br
\label{rmk:whyHausdorff}
If we had only assumed that $X$ was compact and not Hausdorff, then 
$h_{X}$ would not be one-to-one. To see this, let $x,x'\in X$ be two distinct
points that cannot be separated by open sets. For any continuous function
$\vf:X\to\C,$ the value of $\vf$ on these two points is forced to be the same, 
$\vf(x)=\vf(x'),$ i.e. $h_{X}(x)=h_{X}(x').$
\er

\bn
\label{prop:hnaturaltransform}
The assignment
\be
\begin{split}
\cH_{0}&\xrightarrow{h}\cH_{1}\\
X&\mapsto\Big(X
\xrightarrow{h_{X}}\s\big(C(X)\big)\Big)
\end{split}
\ee
from Proposition \ref{prop:evhomeo}
defines a natural isomorphism 
$h:\id_{\cH}\Rightarrow\s\circ C.$ 
\en

\bprf
Associated to any morphism $X\xrightarrow{f}Y$
of compact Hausdorff spaces is the diagram
\be
\xy0;/r.25pc/:
(-12.5,7.5)*+{X}="X";
(-12.5,-7.5)*+{Y}="Y";
(12.5,7.5)*+{\s\big(C(X)\big)}="scX";
(12.5,-7.5)*+{\s\big(C(Y)\big)}="scY";
{\ar"X";"Y"_{f}};
{\ar"scX";"scY"^{\s^{C(f)}}};
{\ar"X";"scX"^(0.35){h_{X}}};
{\ar"Y";"scY"_(0.35){h_{Y}}};
\endxy
.
\ee
To see that this diagram commutes, let $x\in X.$ 
Going along the top and right arrows gives
$
\s^{C(f)}\big(h_{X}(x)\big)
=h_{X}(x)\circ C^{f}
$
by definition of $\s$ as a functor 
(see Proposition \ref{prop:spectrumofmap}).
Going along the left and bottom arrows gives
$
h_{Y}\big(f(x)\big).
$
To see that these two elements of
$\s\big(C(Y)\big),$ i.e.
non-zero characters $C(Y)\to\C,$ are equal,
let $\varphi:Y\to\C$ be a continuous function on $Y.$
Then 
\be
\xy0;/r.20pc/:
(11.7557,16.1803)*+{\qquad\qquad
\Big(h_{Y}\big(f(x)\big)\Big)\Big(\varphi\Big)}="1";
(19.0211,-6.1803)*+{\varphi\big(f(x)\big)}="2";
(0,-20)*+{(\varphi\circ f)(x)}="3";
(-19.0211,-6.1803)*+{
\big(h_{X}(x)\big)\big(\varphi\circ f\big)}="4";
(-11.7557,16.1803)*+{
\big(h_{X}(x)\big)\big(C^{f}(\varphi)\big)
\qquad\qquad}="5";
{\ar@{=}@/^0.7pc/"1";"2"^(0.55){(\ref{eq:evaluation})}};
{\ar@{=}@/^0.7pc/"2";"3"};
{\ar@{=}@/^0.7pc/"3";"4"^(0.70){(\ref{eq:evaluation})}};
{\ar@{=}@/^0.7pc/"4";"5"^(0.45){(\ref{eq:Conmorphisms})}};
\endxy
,
\ee
which proves that $h$ is a 
natural transformation. It is a natural 
\emph{iso}morphism by Proposition 
\ref{prop:evhomeo}.
\eprf

\subsection{The commutative Gelfand-Naimark equivalence}

The following theorem is a categorical phrasing of a theorem
due to Gelfand and Naimark (Lemma 1 in \cite{GN43}). 
It describes in what sense
the functors $C$ and $\s$ are inverses of each other. 
The natural transformations $\G$ and $h$ also satisfy a universal 
property that indicates in what sense they are inverses of each other. 

\bt
Using the notation from above, 
\be
\left(\cCAlg^{\op}\xrightarrow{\s}\cH,
\cH\xrightarrow{C}\cCAlg^{\op},
\id\xRightarrow{\G}C\circ\s,
\s\circ C\xRightarrow{h^{-1}}\id
\right)
\ee
is an adjoint equivalence of categories. 
\et

\bprf
Because of Propositions \ref{prop:Gelfandtransform}
and \ref{prop:hnaturaltransform}, 
it suffices to check the zig-zag identities. 
The first one is given by 
\be
\xy0;/r.32pc/:
(-30,0)*+{\cH}="1";
(-10,0)*+{\cCAlg^{\op}}="2";
(10,0)*+{\cH}="3";
(30,0)*+{\cCAlg^{\op}}="4";
{\ar"1";"2"|-(0.45){C}};
{\ar"2";"3"|-{\s}};
{\ar"3";"4"|-(0.45){C}};
{\ar@{-->}@/^2.25pc/"1";"2"^{C}};
{\ar@/_2.75pc/"1";"3"_{\id_{\cH}}};
{\ar@/^2.75pc/"2";"4"^{\id_{\cCAlg^{\op}}}};
{\ar@{-->}@/_2.25pc/"3";"4"_{C}};
{\ar@{==>}(-20,6.0);(-20,1.0)^{\id_{C}}};
{\ar@{=>}"2";"2"+(0,-7.5)^{h^{-1}}};
{\ar@{=>}"3"+(0,7.5);"3"^(0.35){\G}};
{\ar@{==>}(20,-1.0);(20,-6.0)^{\id_{C}}};
\endxy
\ =\ 
\xy0;/r.32pc/:
(-10,0)*+{\cH}="1";
(10,0)*+{\cCAlg^{\op}}="4";
{\ar@/_1.5pc/"1";"4"_{C}};
{\ar@/^1.5pc/"1";"4"^{C}};
{\ar@{=>}(0,4.0);(0,-4.0)|-{\id_{C}}};
\endxy
\ee
which translates to commutativity of the diagram%
\footnote{The arrows are in the direction required
by contravariance of the functors.}
\be
\xy0;/r.25pc/:
(-15,-7.5)*+{C(X)}="3";
(0,7.5)*+{C\Big(\s\big(C(X)\big)\Big)}="2";
(15,-7.5)*+{C(X)}="1";
{\ar"1";"3"^{\id_{C(X)}}};
{\ar"1";"2"_(0.35){\G_{C(X)}}};
{\ar"3";"2"^(0.35){C(h_{X}^{-1})}};
\endxy
\ee
for every compact Hausdorff space $X.$ 
Because all morphisms here are invertible,
it is equivalent to show that
the diagram 
\be
\xy0;/r.25pc/:
(-15,-7.5)*+{C(X)}="3";
(0,7.5)*+{C\Big(\s\big(C(X)\big)\Big)}="2";
(15,-7.5)*+{C(X)}="1";
{\ar"1";"3"^{\id_{C(X)}}};
{\ar"1";"2"_(0.35){\G_{C(X)}}};
{\ar"2";"3"_(0.60){C(h_{X})}};
\endxy
\ee
commutes, which would provide a sense in which $h$ is the left inverse of $\G.$ 
Thus, let $\varphi\in C(X).$ 
Applying the composition along the top
two arrows to this element gives
\be
C^{h_{X}}\big(\G_{C(X)}(\varphi)\big)
=\G_{C(X)}(\varphi)\circ h_{X}
\ee
by Proposition \ref{prop:Conmorphisms}. This is a map
$X\to\s\big(C(X)\big)\to\C$ and is therefore
determined pointwise so let $x\in X$ and apply
this map to it. The result is
\be
\big(\G_{C(X)}(\varphi)\big)\big(h_{X}(x)\big)
\overset{(\ref{eq:Gelfandtransform})}{=\joinrel=}
h_{X}(x)\big(\varphi\big)
\overset{(\ref{eq:evaluation})}{=\joinrel=}
\varphi(x),
\ee
which proves the first zig-zag identity.
The other zig-zag identity
\be
\xy0;/r.32pc/:
(-30,0)*+{\cCAlg}="1";
(-10,0)*+{\cH^{\op}}="2";
(10,0)*+{\cCAlg}="3";
(30,0)*+{\cH^{\op}}="4";
{\ar"1";"2"|-(0.45){\s}};
{\ar"2";"3"|-{C}};
{\ar"3";"4"|-(0.45){\s}};
{\ar@{-->}@/^2.25pc/"3";"4"^{\s}};
{\ar@/_2.75pc/"2";"4"_{\id_{\cH^{\op}}}};
{\ar@/^2.75pc/"1";"3"^{\id_{\cCAlg}}};
{\ar@{-->}@/_2.25pc/"1";"2"_{\s}};
{\ar@{==>}(20,6.0);(20,1.0)^{\id_{\s}}};
{\ar@{=>}"3";"3"+(0,-7.5)^{h^{-1}}};
{\ar@{=>}"2"+(0,7.5);"2"^(0.35){\G}};
{\ar@{==>}(-20,-1.0);(-20,-6.0)^{\id_{\s}}};
\endxy
\ =\ 
\xy0;/r.32pc/:
(-10,0)*+{\cCAlg}="1";
(10,0)*+{\cH^{\op}}="4";
{\ar@/_1.5pc/"1";"4"_{\s}};
{\ar@/^1.5pc/"1";"4"^{\s}};
{\ar@{=>}(0,4.0);(0,-4.0)|-{\id_{\s}}};
\endxy
\ee
follows by a similar calculation. For a commutative $C^*$-algebra $\mA,$ 
this identity says that the diagram 
\be
\xy0;/r.25pc/:
(-15,-7.5)*+{\s(\mA)}="3";
(0,7.5)*+{\s\Big(C\big(\s(\mA)\big)\Big)}="2";
(15,-7.5)*+{\s(\mA)}="1";
{\ar"1";"3"^{\id_{\s(\mA)}}};
{\ar"1";"2"_(0.35){h_{\s(\mA)}}};
{\ar"2";"3"_(0.60){\s(\G_{\mA})}};
\endxy
\ee
commutes, which provides a sense in which $h$ is the right inverse of $\G.$ 
\eprf

The upshot of the Gelfand-Naimark Theorem that we emphasize here
is that it provides one with the perspective that the study of topology 
(compact Hausdorff spaces) and their continuous functions
is equivalent to the study of \emph{commutative} 
$C^*$-algebras and their $*$-homomorphisms. Therefore, 
\emph{non}-commutative $C^*$-algebras can be interpreted as
non-commutative topology though a more satisfying relationship
to topological concepts is still an area of active research \cite{He17}. 
Furthermore, the $*$-homomorphisms between these
non-commutative spaces can be interpreted as deterministic processes. 
In the next section, we will describe how to fit in stochastic maps
(non-deterministic processes)
and motivate non-commutative probability theory.

\section{Abstract probability theory}
\label{sec:stochasticGelfandNaimark}

The table below collects several of the categories that have been used
along with a few new ones. 

\begin{center}
\begin{tabular}{|c|c|c|}
\hline
Category name&Objects&Morphisms\\
\hline
\hline
$\FinSetFun$&finite sets&functions\\
\hline
$\FinSetStoch$&finite sets&stochastic maps\\
\hline
$\cH$&compact Hausdorff spaces&continuous functions\\
\hline
$\cHStoch$&compact Hausdorff spaces&(continuous) stochastic maps\\
\hline
$\fdcCAlg$&finite-dimensional commutative $C^*$-algebras&$*$-homomorphisms\\
\hline
$\fdcCAlgPos$&finite-dimensional commutative $C^*$-algebras&positive maps\\
\hline
$\cCAlg$&commutative $C^*$-algebras&$*$-homomorphisms\\
\hline
$\cCAlgPos$&commutative $C^*$-algebras&positive maps\\
\hline
\end{tabular}
\end{center}

In this section, we define the category $\cHStoch$ and
generalize the (commutative) Gelfand-Naimark Theorem
to prove an adjoint equivalence of categories 
\be
\xymatrix{
\cHStoch^{\op}
\ar@<1.0ex>[rr]^(0.525){C} &&
\ar@<1.0ex>[ll]^(0.475){\s}  \cCAlgPos
}
.
\ee
During this procedure, we provide an explicit formula for the composition
in $\cHStoch$ in Proposition \ref{prop:compositionstochastic}
by showing that a stochastic map induces a canonical 
Markov kernel in Lemma \ref{lem:evaluationmeasurable}. 
This formula was shown to be valid in the category of 
stochastic maps on Polish spaces in \cite{Gi82} 
but seems to have not been verified in the 
case of compact Hausdorff spaces. 
One can nevertheless \emph{define} the composition
by passing to the algebra of continuous functions and using 
the Riesz-Markov-Kakutani Representation Theorem \cite{EFHN15}
or states \cite{FuJa13}
though these constructions are a bit formal.
In addition, we construct a stochastic version of the Gelfand spectrum functor 
that does not use the Gelfand transform 
nor the Riesz-Markov-Kakutani Theorem 
in Theorem \ref{thm:representingpullbackstatesasmeasures}.
Instead, Choquet theory \cite{Ph01} is used to send states
onto probability measures on the spectrum.
Several of the key ingredients in these constructions 
are standard \cite{Ka41}, \cite{Ci06}, \cite{Fo07}, \cite{Fu11}
though we illustrate how these results fit into a broader context
by showing how all of these categories fit into the cube 
from (\ref{eq:thecube}) in the introduction. 

Sections \ref{sec:chausmeasures} and \ref{sec:probmeasoncHaus} 
review some measure-theoretic preliminaries. 
The main results begin in
Section \ref{sec:cHausStoch} with the construction
of the category $\cHStoch$ and end in 
Section \ref{sec:stochasticGN} with the stochastic Gelfand-Naimark Theorem.
The reader with the analytic background may skip ahead directly to these sections.
We discuss the relationship between our work and that of others
in more detail in Section \ref{sec:relationshiptoCTPT}. 

\subsection{Compact Hausdorff spaces and complex measures}
\label{sec:chausmeasures}

Given a measurable space $(Y,\mathcal{M}_{Y}),$ where $Y$ is a set
and $\mathcal{M}_{Y}$ is a set of measurable subsets, 
known as a \emph{\uline{$\s$-algebra}},
one can define the notion of a positive measure and a 
complex measure. 
Our measurable
spaces will come from compact Hausdorff topological spaces 
equipped with the Borel subsets for their $\s$-algebra.  
As a result, with the additional assumption of finite total measure, 
the set of positive measures will be a 
subset of the set of complex measures. We review the definitions, 
which will be important in specifying a particular topology on the set
of all such measures (Proposition \ref{prop:seminormonmeasures}), 
following Chapters 1 and 6 of \cite{Ru87}. 
We will obtain from this procedure, a 
topology on the set of all probability measures on a compact Hausdorff space. 
For technical reasons, we will then look at a particular subset of all 
measures satisfying a condition compatible with respect to the topology of the
underlying space---such measures are called regular measures 
(also Radon measures). 
We will explain the technicalities and why they are needed in remarks. 
Henceforth, all topological spaces will be assumed compact and Hausdorff
and the set of measurable subsets will always be taken to be 
the set of Borel subsets.
A discussion of the pros and cons of using Baire sets is deferred to 
Section \ref{sec:otherapproaches}.
Most of this background material does not require the measurable
space to be compact Hausdorff though we phrase everything in this 
setting since all our applications are for such spaces. 

\bd
\label{defn:complexmeasure}
Let $Y$ be a compact Hausdorff space with associated Borel measurable
subsets $\mathcal{M}_{Y}.$ A \emph{\uline{complex measure}} on 
$(Y,\mathcal{M}_{Y})$ is a function $\mu:\mathcal{M}_{Y}\to\C$
that is \emph{\uline{countably additive}} in the sense that for 
every measurable set $E\in\mathcal{M}_{Y},$ 
\be
\mu(E)=\sum_{i\in I}\mu(E_{i})
\ee
for every at most countable collection
$\{E_{i}\}_{i\in I\subseteq\N}$ of measurable sets satisfying%
\footnote{$\varnothing$ denotes the empty set.}
\be
E_{i}\cap E_{j}=\varnothing
\ee
whenever $i\ne j$ and 
\be
\bigcup_{i\in I}E_{i}=E. 
\ee
Such a collection $\{E_{i}\}_{i\in I}$ is called a \emph{\uline{measurable partition}}
of $E.$ 
A \emph{\uline{positive measure}} is a complex measure whose
value on all measurable sets is non-negative. 
A \emph{\uline{probability measure}} is a positive measure
satisfying $\mu(Y)=1.$ The set of all complex measures on 
$Y$ will be denoted by $\C\mathrm{Meas}(Y).$ 
\ed

Occasionally, we may write $Y$ instead of $(Y,\mathcal{M}_{Y})$ since our set of
measurable sets will always be Borel. 
Note that every positive  
measure can be
used to define an integral using simple functions 
(see Chapter 1 of \cite{Ru87}).
Since the construction is used frequently for positive measurable functions, 
we state it here while also setting some notation. 
\bd
\label{defn:integral}
Let $Y$ be a compact Hausdorff space and let $E\in\mathcal{M}_{Y}.$ 
The function $\chi_{E}:Y\to\R$ defined by 
\be
Y\ni y\mapsto\chi_{E}(y):=\begin{cases}1&\mbox{ if } y\in E\\0&\mbox{ otherwise}\end{cases}
\ee
is called the \emph{\uline{characteristic function}} on $E.$ 
A \emph{\uline{simple function}} is a function $s:Y\to\R_{\ge0}$ of the form 
\be
s=\sum_{i}^{\text{finite}}s_{i}\chi_{E_{i}}
\ee
for some finite set of non-negative numbers $\{s_i\}$ 
and measurable subsets $\{E_{i}\}.$
Let $\mu$ be a positive measure on $Y.$ 
For a measurable function $f:Y\to\R_{\ge0},$ the 
\emph{\uline{integral of $f$ with respect to $\mu$}} is the number
\be
\int_{Y}f\;d\mu:=\sup_{0\le s\le f}\left\{\sum_{i}s_{i}\mu(E_{i})\;:\;s=\sum_{i}s_{i}\chi_{E_{i}}\text{ is simple}\right\}.
\ee
\ed
The integral of a complex measurable function can be defined by splitting the function
up into a sum of its real and imaginary parts each of which can be decomposed
into a sum of positive and negative parts. 
Similarly, every complex measure can be uniquely decomposed 
via Jordan decomposition into an appropriate sum and difference of 
positive measures so as to define the integral of a positive function 
with respect to a complex measure and then finally the integral
of any complex measurable function with respect to a complex
measure (see Section 6.18 of \cite{Ru87}---we will review a necessary
fact to make sense of this definition in Theorem 
\ref{thm:polardecompofmeasure} below). 

The following Lemma is useful for several of the facts that will follow.
It is unfortunately long, but is used so often that we felt it was important
to state it in its full form. 

\blem
\label{lem:contcompactsequenceapprox}
Let $Z$ be a compact topological space and let $\vf:Z\to\R_{\ge0}$ be
a bounded measurable function. 
Then there exists a sequence $s:\N\to \R^{Z}$ 
of simple functions satisfying
\be
0\le s_1\le s_2\le s_3\le \cdots\le \vf
\aand
\lim_{n\to\infty}s_{n}(z)=\vf(z)\quad\forall\;z\in Z. 
\ee
Furthermore,
the coefficients $\{s_{n,i_n}\}$ and measurable sets $\{E_{n,i_{n}}\}$
in the expression 
\be
s_{n}=\sum_{i_{n}}s_{n,i_{n}}\chi_{E_{n,i_n}}
\ee
can be chosen so that 
$i_{n}\in\{0,1,\dots,2^{n}-1\},$ 
\be
\label{eq:splittingEtotwoparts}
E_{n,i_{n}}=E_{n+1,2i_{n}}\cup E_{n+1,2i_{n}+1}
\aand
E_{n+1,2i_{n}}\cap E_{n+1,2i_{n}+1}=\varnothing, 
\ee
and 
\be
\label{eq:breakingsimpleintotwo}
s_{n,i_{n}}\le\min\big\{s_{n+1,2i_{n}},s_{n+1,2i_{n}+1}\big\}
\ee
for all $i_{n}\in\{0,1,\dots,2^{n}-1\}$ and $n\in\N.$ 
\elem

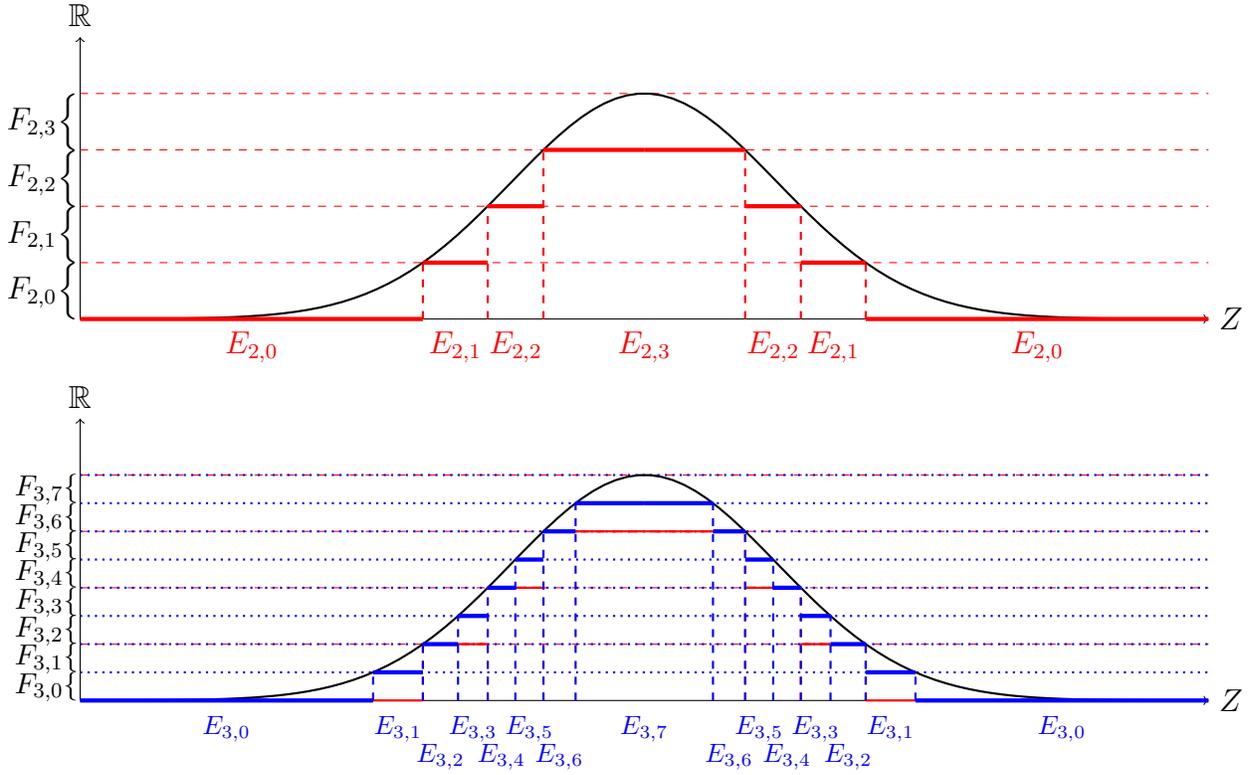
\begin{figure}
\centering
\begin{tikzpicture}[xscale=2.5,yscale=3.0]
    \draw[->] (-3,0) -- (3,0) node[right] {$Z$};
    \draw[->] (-3,0) -- (-3,1.25) node[above] {$\R$};
    \node at (-3.20,0.125) {$F_{2,0}\Big\{$};
    \node at (-3.20,0.375) {$F_{2,1}\Big\{$};
    \node at (-3.20,0.625) {$F_{2,2}\Big\{$};
    \node at (-3.20,0.875) {$F_{2,3}\Big\{$};
    \node[red] at ({(-3-sqrt(ln(1/0.25)))/2},-0.125) {$E_{2,0}$};
    \node[red] at ({(-sqrt(ln(1/0.25))-sqrt(ln(1/0.5)))/2},-0.125) {$E_{2,1}$};
    \node[red] at ({(-sqrt(ln(1/0.5))-sqrt(ln(1/0.75)))/2},-0.125) {$E_{2,2}$};
    \node[red] at (0,-0.125) {$E_{2,3}$};
    \node[red] at ({(sqrt(ln(1/0.5))+sqrt(ln(1/0.75)))/2},-0.125) {$E_{2,2}$};
    \node[red] at ({(sqrt(ln(1/0.25))+sqrt(ln(1/0.5)))/2},-0.125) {$E_{2,1}$};
    \node[red] at ({(3+sqrt(ln(1/0.25)))/2},-0.125) {$E_{2,0}$};
    \foreach \m in {0.25,0.5,0.75,1}
    {
    \draw[dashed,red] (-3,\m) -- (3,\m);
    }
    \draw[thick,domain=-3:3,samples=100] plot (\x,{exp(-\x*\x)});
    \draw[red,ultra thick] (-3,0) -- ({-sqrt(ln(1/.25))},0);
    \draw[red,ultra thick] ({sqrt(ln(1/.25))},0) -- (3,0);
    \foreach \n in {.25,.5,.75}
    {
    \draw[red,ultra thick] ({-sqrt(ln(1/\n))},{\n}) -- ({-sqrt(ln(1/(\n+.25)))},{\n});
    }
    \foreach \n in {.25,.5,.75}
    {
    \draw[red,ultra thick] ({sqrt(ln(1/\n))},{\n}) -- ({sqrt(ln(1/(\n+.25)))},{\n});
    }
    \foreach \n in {.25,.5,.75}
    {
    \draw[red,thick,dashed] ({-sqrt(ln(1/\n))},{0}) -- ({-sqrt(ln(1/\n))},{\n}) -- ({-sqrt(ln(1/(\n+.25)))},{\n});
    }
    \foreach \n in {.25,.5,.75}
    {
    \draw[red,thick,dashed] ({sqrt(ln(1/\n))},{0}) -- ({sqrt(ln(1/\n))},{\n}) -- ({sqrt(ln(1/(\n+.25)))},{\n});
    }
\end{tikzpicture}
%
%
\begin{tikzpicture}[xscale=2.5,yscale=3.0]
    \draw[->] (-3,0) -- (3,0) node[right] {$Z$};
    \draw[->] (-3,0) -- (-3,1.25) node[above] {$\R$};
    \node at (-3.20,0.0625) {\small$\;F_{3,0}\{$};
    \node at (-3.20,0.1875) {\small$\;F_{3,1}\{$};
    \node at (-3.20,0.3125) {\small$\;F_{3,2}\{$};
    \node at (-3.20,0.4375) {\small$\;F_{3,3}\{$};
    \node at (-3.20,0.5625) {\small$\;F_{3,4}\{$};
    \node at (-3.20,0.6875) {\small$\;F_{3,5}\{$};
    \node at (-3.20,0.8125) {\small$\;F_{3,6}\{$};
    \node at (-3.20,0.9375) {\small$\;F_{3,7}\{$};
    \node[blue] at ({(-3-sqrt(ln(1/0.125)))/2},-0.125) {\footnotesize$E_{3,0}$};
    \node[blue] at ({(-sqrt(ln(1/0.125))-sqrt(ln(1/0.25)))/2},-0.125) {\footnotesize$E_{3,1}$};
    \node[blue] at ({(-sqrt(ln(1/0.25))-sqrt(ln(1/0.375)))/2},-0.25) {\footnotesize$E_{3,2}$};
    \node[blue] at ({(-sqrt(ln(1/0.375))-sqrt(ln(1/0.5)))/2},-0.125) {\footnotesize$E_{3,3}$};
    \node[blue] at ({(-sqrt(ln(1/0.5))-sqrt(ln(1/0.625)))/2},-0.25) {\footnotesize$E_{3,4}$};
    \node[blue] at ({(-sqrt(ln(1/0.625))-sqrt(ln(1/0.75)))/2},-0.125) {\footnotesize$E_{3,5}$};
    \node[blue] at ({(-sqrt(ln(1/0.75))-sqrt(ln(1/0.875)))/2},-0.25) {\footnotesize$E_{3,6}$};
    \node[blue] at (0,-0.125) {\footnotesize$E_{3,7}$};
    \node[blue] at ({(sqrt(ln(1/0.75))+sqrt(ln(1/0.875)))/2},-0.25) {\footnotesize$E_{3,6}$};
    \node[blue] at ({(sqrt(ln(1/0.625))+sqrt(ln(1/0.75)))/2},-0.125) {\footnotesize$E_{3,5}$};
    \node[blue] at ({(sqrt(ln(1/0.5))+sqrt(ln(1/0.625)))/2},-0.25) {\footnotesize$E_{3,4}$};
    \node[blue] at ({(sqrt(ln(1/0.375))+sqrt(ln(1/0.5)))/2},-0.125) {\footnotesize$E_{3,3}$};
    \node[blue] at ({(sqrt(ln(1/0.25))+sqrt(ln(1/0.375)))/2},-0.25) {\footnotesize$E_{3,2}$};
    \node[blue] at ({(sqrt(ln(1/0.125))+sqrt(ln(1/0.25)))/2},-0.125) {\footnotesize$E_{3,1}$};
    \node[blue] at ({(3+sqrt(ln(1/0.125)))/2},-0.125) {\footnotesize$E_{3,0}$};
    \foreach \m in {0.25,0.5,0.75,1}
    {
    \draw[dashed,red] (-3,\m) -- (3,\m);
    }
    \foreach \k in {0.125,0.25,...,1}
    {
    \draw[dotted,thick,blue] (-3,\k) -- (3,\k);
    }
    \draw[thick,domain=-3:3,samples=100] plot (\x,{exp(-\x*\x)});
    \draw[red,thick] (-3,0) -- ({-sqrt(ln(1/.25))},0);
    \draw[red,thick] ({sqrt(ln(1/.25))},0) -- (3,0);
    \foreach \n in {.25,.5,.75}
    {
    \draw[red,thick] ({-sqrt(ln(1/\n))},{\n}) -- ({-sqrt(ln(1/(\n+.25)))},{\n});
    }
    \foreach \n in {.25,.5,.75}
    {
    \draw[red,thick] ({sqrt(ln(1/\n))},{\n}) -- ({sqrt(ln(1/(\n+.25)))},{\n});
    }
    \foreach \n in {.25,.5,.75}
    {
    \draw[red,thick,dashed] ({-sqrt(ln(1/\n))},{0}) -- ({-sqrt(ln(1/\n))},{\n}) -- ({-sqrt(ln(1/(\n+.25)))},{\n});
    }
    \foreach \n in {.25,.5,.75}
    {
    \draw[red,thick,dashed] ({sqrt(ln(1/\n))},{0}) -- ({sqrt(ln(1/\n))},{\n}) -- ({sqrt(ln(1/(\n+.25)))},{\n});
    }
    \draw[blue,ultra thick] (-3,0) -- ({-sqrt(ln(1/0.125))},0);
    \draw[blue,ultra thick] ({sqrt(ln(1/0.125))},0) -- (3,0);
    \foreach \n in {0.125,0.25,...,0.875}
    {
    \draw[blue,ultra thick] ({-sqrt(ln(1/\n))},{\n}) -- ({-sqrt(ln(1/(\n+0.125)))},{\n});
    }
    \foreach \n in {0.125,0.25,...,0.875}
    {
    \draw[blue,ultra thick] ({sqrt(ln(1/\n))},{\n}) -- ({sqrt(ln(1/(\n+0.125)))},{\n});
    }
    \foreach \n in {0.125,0.25,...,0.875}
    {
    \draw[blue,thick,dashed] ({-sqrt(ln(1/\n))},{0}) -- ({-sqrt(ln(1/\n))},{\n}) -- ({-sqrt(ln(1/(\n+0.125)))},{\n});
    }
    \foreach \n in {0.125,0.25,...,0.875}
    {
    \draw[blue,thick,dashed] ({sqrt(ln(1/\n))},{0}) -- ({sqrt(ln(1/\n))},{\n}) -- ({sqrt(ln(1/(\n+0.125)))},{\n});
    }
\end{tikzpicture}
\caption{A non-negative bounded measurable function is an increasing 
sequence of simple functions satisfying several convenient properties
as described in Lemma \ref{lem:contcompactsequenceapprox}.}
\label{fig:simplefunctions}
\end{figure}

\bprf
Since $\vf$ is bounded, there exists an $M\ge0$ such that
$\vf(Z)\subseteq[0,M].$ For each $n\in\N,$ set
\be
F_{n,i_{n}}:=
\begin{cases}
\left[\frac{i_{n}M}{2^{n}},\frac{(i_{n}+1)M}{2^{n}}\right)&\mbox{ for $i_{n}\in[0,2^{n}-1)\cap\Z$}\\
\left[\frac{i_{n}M}{2^{n}},\frac{(i_{n}+1)M}{2^{n}}\right]&\mbox{ for $i_{n}=2^{n}-1$},\\
\end{cases}
\ee
where $i_{n}\in\{0,1,\dots,2^{n}-1\}.$ 
Since $\vf$ is Borel measurable, the preimages 
$E_{n,i_{n}}:=\vf^{-1}\big(F_{n,i_{n}}\big)$
form a measurable partition of $Z$ for each $n.$ 
Set $s_{n}:Z\to\R_{\ge0}$ to be the simple function 
\be
s_{n}:=\sum_{i_n=0}^{2^{n}-1}\frac{i_{n}M}{2^{n}}\chi_{E_{n,i_{n}}}. 
\ee
An example of such a simple function together with the next function
in the sequence is
depicted in Figure \ref{fig:simplefunctions}.
The sequence $\N\ni n\mapsto s_n$ of simple functions 
satisfies the required properties because for every $n,$
the pointwise difference between $s_{n}$ and $\vf$ is at most $\frac{M}{2^{n}}.$ 
Furthermore, $s_{n}\le s_{n+1}$ for each $n\in\N$ 
because $\{F_{n+1,i_{n+1}}\}$ is a refinement of $\{F_{n,i_{n}}\}$ 
and hence likewise $\{E_{n+1,i_{n+1}}\}$ is a refinement of $\{E_{n,i_{n}}\}$ 
(this was the reason for the choice of $2^{n}$ as opposed to merely $n$ 
in the above formulas). 
The other properties are immediate from the construction. 
\eprf

\bt
[Monotone Convergence Theorem]
Let $Y$ be a compact Hausdorff space and let 
$\vf:\N\to[0,\infty]^{Y}$ be an increasing sequence of Borel measurable functions, 
i.e. $\vf_{n+1}(y)\ge\vf_{n}(y)$ for all $y\in Y$ and $n\in\N.$ 
Then $\lim\vf$ is Borel measurable and 
\be
\lim_{n\to\infty}\int_{Y}\vf_{n}\;d\mu=\int_{Y}\lim\vf\;d\mu
\ee
for every Borel measure $\mu$ on $Y.$ 
\et

\bprf
See Theorem 1.26 in \cite{Ru87}. 
\eprf

\bd
Let $Y$ be a compact Hausdorff space. 
Given a complex measure $\mu\in\C\mathrm{Meas}(Y),$ the 
\emph{\uline{total variation}} of $\mu$ is the function 
$|\mu|:\mathcal{M}_{Y}\to\R$ defined by 
\be
\mathcal{M}_{Y}\ni E\xmapsto{|\mu|}\sup\left\{\sum_{i\in I}\big|\mu(E_{i})\big|
\right\},
\ee
where the supremum is over all measurable partitions $\{E_{i}\}_{i\in I}$ of $E.$
\ed

\bn
Let $Y$ be a compact Hausdorff space and $\mu\in\C\mathrm{Meas}(Y).$ 
Then, the total variation of $\mu$ is a positive measure on $Y.$
Furthermore, the set of all complex measures on $Y$
has the structure of a normed complex vector space with norm given by 
\be
\C\mathrm{Meas}(Y)\ni\mu\mapsto|\mu|(Y).
\ee
\en

\bprf
See Chapter 6 of \cite{Ru87}. 
\eprf

\bt
[Polar decomposition for complex measures]
\label{thm:polardecompofmeasure}
Let $Y$ be a compact Hausdorff space and $\mu\in\C\mathrm{Meas}(Y).$ 
Then there exists a measurable function $h:Y\to\C$ such that
$|h(y)|=1$ for all $y\in Y$ and 
\be
\mu(E)=\int_{E}h\;d|\mu|\qquad\forall\;E\in\mathcal{M}_{Y}.
\ee
This is often written in shorthand form as 
$d\mu=h\;d|\mu|.$
\et

\bprf
See Theorem 6.12 in \cite{Ru87}. 
\eprf

This is an important consequence of the
Lebesgue-Radon-Nikodym Theorem.
It allows one to define the integral of complex-valued 
measurable functions with respect to complex measures, 
namely using the notation of 
Theorem \ref{thm:polardecompofmeasure}, the integral of 
a complex-valued measurable function $\vf:Y\to\C$ with respect to $\mu$ is
\be
\int_{Y}\vf\;d\mu:=\int_{Y}\vf h\;d|\mu|. 
\ee
This definition satisfies the usual properties of the integral 
(see Section 6.18 in \cite{Ru87}). 
This theorem can also be used to prove the following result, which we will need.

\blem
\label{lem:absvalueonintegral}
Let $Y$ be a compact Hausdorff space and $\mu\in\C\mathrm{Meas}(Y).$ 
Then 
\be
\left|\int_{Y}\vf\;d\mu\right|\le\int_{Y}|\vf|\;d|\mu|
\ee
for all $\vf\in C(Y).$ 
\elem

\bprf
Using a polar decomposition $d\mu=h\;d|\mu|,$ 
\be
\left|\int_{Y}\vf\;d\mu\right|=\left|\int_{Y}\vf h\;d|\mu|\right|
\le\int_{Y}|\vf||h|\;d|\mu|=\int_{Y}|\vf|\;d|\mu|
\ee
for all $\vf\in C(Y).$ 
\eprf

Being a normed space, $\C\mathrm{Meas}(Y)$ has a natural topology on it. 
However, as for many infinite-dimensional vector spaces, too few sequences
have convergent subsequences with respect to this topology. There are several
ways to implement topologies on $\C\mathrm{Meas}(Y).$ Rather than using
all measurable functions, we will use continuous functions to define a collection
of seminorms (recall the discussion in Section 
\ref{sec:topologicalpreliminaries}). 

\bn
\label{prop:seminormonmeasures}
Let $Y$ be a compact Hausdorff space. For every continuous function $\vf:Y\to\C,$
set $\lVert\;\cdot\;\rVert_{\vf}:\C\mathrm{Meas}(Y)\to\R$ to be the function given by
\be
\C\mathrm{Meas}(Y)\ni\mu\mapsto\int_{Y}|\vf|\;d|\mu|. 
\ee 
Then $\lVert\;\cdot\;\rVert_{\vf}$ is a seminorm on $\C\mathrm{Meas}(Y).$ 
\en

\bprf
This follows immediately from the definition of the integral and the fact that
the total variation $|\;\cdot\;|$ defines a norm on $\C\mathrm{Meas}(Y).$
Note that the function $\lVert\;\cdot\;\rVert_{\vf}$ makes sense even when 
$\vf$ is just a measurable function. 
\begin{enumerate}[i.]
\item
$\lVert\mu\rVert_{\vf}=\int_{Y}|\vf|\;d|\mu|\ge0$ 
for all complex measures $\mu\in\C\mathrm{Meas}(Y)$ 
follows immediately from 
the definition of the integral of a positive function with 
respect to a positive measure. 
\item
Fix $\l\in\C$ and $\mu\in\C\mathrm{Meas}(Y).$ 
By Lemma \ref{lem:contcompactsequenceapprox}, there exists a sequence
$s:\N\to\R^{Y}$ of simple functions on $Y$ of the form 
\be
s_{n}:=\sum_{i_{n}}^{\text{finite}}s_{n,i_{n}}\chi_{E_{n,i_n}},
\ee
satisfying 
\be
0\le s_{1}\le s_{2}\le s_{3}\le \cdots\le|\vf|
\ee
and 
\be
\lim_{n\to\infty} s_{n}(y)=|\vf(y)|\qquad\forall\;y\in Y.
\ee
By the definition of the integral 
\be
\lim_{n\to\infty}\int_{Y}s_{n}\;d|\l\mu|=\int_{Y}|\vf|\;d|\l\mu|
=\lVert\l\mu\rVert_{\vf}.
\ee
Meanwhile, by definition of simple functions, 
the left-hand-side of this expression equals
\be
\begin{split}
\lim_{n\to\infty}\int_{Y}s_{n}\;d|\l\mu|
&=\lim_{n\to\infty}\sum_{i_{n}=1}^{\text{finite}}|\l\mu|(E_{n,i_{n}})s_{n,i_{n}}\\
&=\lim_{n\to\infty}\sum_{i_{n}=1}^{\text{finite}}|\l||\mu|(E_{n,i_{n}})s_{n,i_{n}}\\
&=|\l|\lim_{n\to\infty}\sum_{i_{n}=1}^{\text{finite}}|\mu|(E_{n,i_{n}})s_{n,i_{n}}\\
&=|\l|\lim_{n\to\infty}\int_{Y}s_{n}\;d|\mu|\\
&=|\l|\int_{Y}|\vf|\;d|\mu|\\
&=|\l|\lVert\mu\rVert_{\vf}.
\end{split}
\ee
Hence, $\lVert\l\mu\rVert_{\vf}=|\l|\lVert\mu\rVert_{\vf}$ as needed. 
\item
Let $\mu,\nu\in\C\mathrm{Meas}(Y)$ and let $s$ be a sequence
as in the previous part. A similar argument as above shows that
\be
\lVert\mu+\nu\rVert_{\vf}=\int_{Y}|\vf|\;d|\mu+\nu|
\le\int_{Y}|\vf|\;d|\mu|+\int_{Y}|\vf|\;d|\nu|
=\lVert\mu\rVert_{\vf}+\lVert\nu\rVert_{\vf}.
\ee
\end{enumerate}
All of the conditions of a seminorm have therefore been checked. 
\eprf

Henceforth, we equip $\C\mathrm{Meas}(Y)$ with the topology 
generated by the family of seminorms 
$\big\{\lVert\;\cdot\;\rVert_{\vf}\big\}_{\vf\in C(Y)}.$
This topology is known as the \emph{\uline{vague topology}} on 
$\C\mathrm{Meas}(Y).$ It is helpful to know when nets converge
with respect to this topology so that we may use nets to prove
whether or not certain functions are continuous. 

\bn
Let $Y$ be a compact Hausdorff space.
A net $\mu:\Theta\to\C\mathrm{Meas}(Y)$ converges to $\lim\mu$ 
with respect to the vague topology if and only if 
\be
\lim_{\q\in\Theta}\int_{Y}\vf\;d\mu_{\q}=\int_{Y}\vf\;d\lim\mu
\ee
for all $\vf\in C(Y).$ 
\en

\bprf
{\color{white}{space}}

\noindent
($\Rightarrow$)
Suppose that $\mu:\Theta\to\C\mathrm{Meas}(Y)$ converges to $\lim\mu.$ 
Fix $\e>0.$ By assumption, for any continuous function $\vf\in C(Y),$ 
there exists a $\q_{0}\in\Theta$ such that 
\be
\lVert\lim\mu-\mu_{\q}\rVert_{\vf}<\e\qquad\forall\;\q\ge\q_{0}.
\ee
Then, by Lemma \ref{lem:absvalueonintegral}, 
\be
\begin{split}
\left|\int_{Y}\vf\;d\lim\mu-\int_{Y}\vf\;d\mu_{\q}\right|
&\le\left|\int_{Y}\vf\;d(\lim\mu-\mu_{\q})\right|\\
&\le\int_{Y}|\vf|\;d|\lim\mu-\mu_{\q}|\\
&=\lVert\lim\mu-\mu_{\q}\rVert_{\vf}\\
&<\e\qquad\forall\;\q\ge\q_{0}.
\end{split}
\ee

\noindent
($\Leftarrow$)
Suppose that $\ds\lim_{\q\in\Theta}\int_{Y}\vf\;d\mu_{\q}=\int_{Y}\vf\;d\lim\mu$ 
for all $\vf\in C(Y).$ Fix $\e>0,$ $n\in\N,$ and $\a:\ov{n}\to C(Y)$ whose value
at $i\in\ov{n}$ will be denoted by $\a_{i}$ (we will show that $\mu$ is eventually
in a neighborhood base of $\lim\mu$). By assumption, for each $i\in\ov{n},$ 
there exist $\q_{i}\in\Theta$ such that 
\be
\left|\int_{Y}\vf\;d\lim\mu-\int_{Y}\vf\;d\mu_{\q}\right|<\e\qquad\forall\;\q\ge\q_{i}.
\ee
Set 
\be
\q_{0}:=\max_{i\in\ov{n}}\{\q_{i}\}.
\ee
Then, $\mu_{\q}\in B^{\a}_{\e}(\lim\mu)$ for all $\q\ge\q_{0}.$ 
Hence $\mu$ converges to $\lim\mu$ with respect to the vague topology.  
\eprf

\subsection{Probability measures on compact Hausdorff spaces}
\label{sec:probmeasoncHaus}

Naively, one might use the set of \emph{all} probability measures on $Y$
and define suitably continuous stochastic maps between compact 
Hausdorff spaces using this definition. 
However, in order to obtain an equivalence of categories between
such spaces and the category of commutative $C^*$-algebras and positive maps, 
an additional restriction needs to be made on the set of probability measures
to ensure that the equivalence holds. 
This restriction is solely due to the fact that not all compact Hausdorff
spaces are separable, or equivalently, metrizable
(for instance, see counterexamples 24, 43, 107, and 111 from \cite{StSe95}). 
The measures we must restrict our attention to are regular (Radon) measures. 
Since we will mainly be concerned with probability measures, 
we will immediately assume throughout that all measures are finite. 

\bd
Let $Y$ be a compact Hausdorff space. A positive Borel measure 
$\mu:\mathcal{M}_{Y}\to\R_{\ge0}$ is \emph{\uline{regular}} iff
\be
\mu(E)=\inf\big\{\mu(U)\;:\; E\subseteq U\; \text{with $U$ open in $Y$} \big\}
\ee
and
\be
\mu(E)=\sup\big\{\mu(K)\;:\; K\subseteq E\; \text{with $K$ compact in $Y$}\big\}
\ee
for all $E\in\mathcal{M}_{Y}.$ 
A complex Borel measure 
$\mu:\mathcal{M}_{Y}\to\C$ is \emph{\uline{regular}} iff the 
positive and negative parts of the real and complex parts of $\mu$ 
are all regular. 
\ed

Because being regular is a property of a measure, the set of 
regular measures on $Y$ is a subset 
$\C\mathrm{Reg}(Y)\subseteq\C\mathrm{Meas}(Y)$ of all 
complex measures. Since $\C\mathrm{Meas}(Y)$ has a normed
vector space structure, one could ask if 
$\C\mathrm{Reg}(Y)$ inherits the properties of this structure as well. 

\bt
$\C\mathrm{Reg}(Y)$ with the total variation norm 
is a Banach space for all compact Hausdorff spaces $Y.$ 
\et

\bprf
See Exercise 3 in Chapter 6 of \cite{Ru87}. 
\eprf

\begin{notation}
\label{not:ProbMeas}
Because $\C\mathrm{Reg}(Y)\subseteq\C\mathrm{Meas}(Y),$
it also inherits the vague topology described after 
Proposition \ref{prop:seminormonmeasures}. 
Let $\mathrm{ProbMeas}(Y)$ denote the subset of $\C\mathrm{Reg}(Y)$
of regular measures that are also probability measures. 
Equip $\mathrm{ProbMeas}(Y)$ with the subspace
topology coming from the vague topology on $\C\mathrm{Meas}(Y).$
\end{notation}

Notice that when $Y$ is a finite set with the discrete topology, 
$\mathrm{ProbMeas}(Y)=\mathrm{Pr}(Y).$

\br
When $Y$ is just a \emph{locally} compact Hausdorff space, 
one must use Radon probability measures instead 
(see Theorem 7.2 in \cite{Fo07}).
When $Y$ is a compact Hausdorff space, it is a fact that every
Radon probability measure is regular (see Corollary 7.6 in \cite{Fo07}).
When $Y$ is, in addition, second countable, \emph{every} Borel measure
on $Y$ is regular (see Theorem 7.8 in \cite{Fo07}). 
Hence, many of the complications that follow can be ignored provided
that $Y$ is second countable (and therefore separable). 
In fact, many treatments of probability theory on topological spaces
assume the spaces are Polish. We do not need this requirement, which
is one of the interesting features of our findings. 
\er

We mention several facts regarding regular measures on compact
Hausdorff spaces that will be needed for many of our proofs. 

\blem
[Urysohn's Lemma]
\label{lem:Urysohn}
Let $Y$ be a compact Hausdorff space and let $K\subseteq U\subseteq Y$ with
$K$ closed (and hence compact) and $U$ open. 
Then there exists a continuous function $\vf\in C\big(Y,[0,1]\big)$ such that 
$\vf=1$ on $K$ and $\vf=0$ outside a compact subset of $U.$ 
\elem

\bprf
See Theorem 4.32 in \cite{Fo07}. 
\eprf

\bd
Let $Y$ be a topological space. A function $\vf:Y\to\R$ is
\emph{\uline{lower semi-continuous}} iff 
$f^{-1}\big((a,\infty)\big)$ is open for all $a\in\R$
and is 
\emph{\uline{upper semi-continuous}} iff 
$f^{-1}\big((-\infty,a)\big)$ is open for all $a\in\R.$
\ed

It is immediate from the definition that semi-continuous functions are
measurable. This is because sets of the form $(a,\infty)$ also generate
the Borel $\s$-algebra on $\R$ and because a function is measurable if and 
only if the preimages of generating sets are measurable.   

\bn
\label{prop:semicontinuous}
Let $Y$ be a topological space. 
\begin{enumerate}[(a)]
\setlength{\itemsep}{0pt}
\item
If $U\subseteq Y$ is open, then $\chi_{U},$ the characteristic
function on $U,$ is lower semi-continuous. 
\item
If $K\subseteq Y$ is closed, then $\chi_{K}$ is upper semi-continuous. 
\item
If $Y$ is compact and Hausdorff and $\chi:Y\to\R_{\ge0}$ is 
lower semi-continuous, then 
\be
\chi(y)=\sup_{\vf\in C(Y)}\big\{\vf(y)\;:\;0\le\vf\le\chi\big\}
\qquad\forall\;y\in Y.
\ee
\item
If $Y$ is compact and Hausdorff and $\chi:Y\to\R_{\ge0}$ is 
upper semi-continuous, then 
\be
\chi(y)=\inf_{\vf\in C(Y)}\big\{\vf(y)\;:\;0\le\chi\le\vf\big\}
\qquad\forall\;y\in Y.
\ee

\end{enumerate}
\en

\bprf
See Proposition 7.11 in  \cite{Fo07} (though the reader may also check that
this follows from the definitions and Urysohn's Lemma). 
\eprf

\bc
\label{cor:LSCintegralapprox}
Let $\mu$ be a positive regular measure on a compact Hausdorff space $Y$ 
and let $\chi:Y\to\R_{\ge0}$ be a lower semi-continuous function on $Y.$ Then
\be
\int_{Y}\chi\;d\mu=\sup_{\vf\in C(Y)}\left\{\int_{Y}\vf\;d\mu\;:\;0\le\vf\le\chi\right\}.
\ee
\ec

\bprf
See Corollary 7.13 in \cite{Fo07}. 
\eprf

\br
Corollary \ref{cor:LSCintegralapprox} follows from Lusin's theorem, which is 
a more quantitative version of Urysohn's Lemma for compact Hausdorff
spaces \emph{equipped with} a (finite) regular measure \cite{Fo07}. 
\er

All regular probability measures are regular positive measures, 
and therefore all of these results apply to regular probability measures
on a compact Hausdorff space. 

\subsection{Continuous stochastic maps}
\label{sec:cHausStoch}

We are finally ready to define the category $\cHStoch,$ the 
generalization of the category $\FinSetStoch$ to compact 
Hausdorff spaces and (continuous) stochastic maps. 
In the process of doing so, we prove several facts of independent interest
(see Lemma \ref{lem:evaluationofmeasurecontinuous}
and \ref{lem:evaluationmeasurable} for instance). 

\bd
\label{defn:continuousstochasticprocess}
Let $X$ and $Y$ be two compact Hausdorff spaces. 
A (continuous) \emph{\uline{stochastic map}} from $X$ 
to $Y$ is a continuous function $\mu:X\to\mathrm{ProbMeas}(Y)$
with respect to the vague topology and is denoted by $\mu:X\stoch Y.$ 
The value of $\mu$ at $x\in X$ is denoted by $\mu_{x}$
and is a regular probability measure on $Y.$ 
\ed

\bx
\label{ex:diracstochastic}
Let $f:X\to Y$ be a continuous \emph{function} of compact Hausdorff spaces. 
Then
\be
\begin{split}
X&\xrightarrow{\de_{f}}\mathrm{ProbMeas}(Y)\\
x&\mapsto\de_{f(x)},
\end{split}
\ee
where $\de_{f(x)}$ is the measure defined by
\be
\mathcal{M}_{Y}\ni E\mapsto \de_{f(x)}(E):=
\begin{cases}
1&\mbox{ if }f(x)\in E\\
0&\mbox{ if }f(x)\notin E,\\
\end{cases}
\ee
is a stochastic map from $X$ to $Y.$ 
$\de_{f(x)}$ is called the Dirac delta measure on $Y$ concentrated at the point $f(x).$
The fact that $\de_{f(x)}$ is a regular probability measure follows immediately 
from the definition.
Before proving that this is a stochastic map, notice that 
for any continuous function $\vf\in C(Y)$ and for any $x\in X,$ 
\be
\int_{Y}\vf\;d\de_{f(x)}=\vf\big(f(x)\big).
\ee
To see that $X\ni x\mapsto\de_{f(x)}$ is continuous, let 
$x:\Theta\to X$ be a net converging to $\lim x\in X.$ Then, 
for any continuous function $\vf\in C(Y),$ 
\be
\lim_{\q\in\Theta}\int_{Y}\vf\;d\de_{f(x_{\q})}
=\lim_{\q\in\Theta}\vf\big(f(x_\q)\big)
=\vf\big(f(\lim x)\big)
=\int_{Y}\vf\;d\de_{f(\lim x)}.
\ee
Since $\vf$ was arbitrary, this shows that 
$\ds\lim_{\q\in\Theta}\de_{f(x_\q)}=\de_{f(\lim x)}$ in the vague topology. 

In particular, the identity function from a space to itself induces
a stochastic map known as the \emph{\uline{identity stochastic map}}. 
Why it is called the identity will be justified after the definition of
composition of stochastic maps is provided. 
\ex

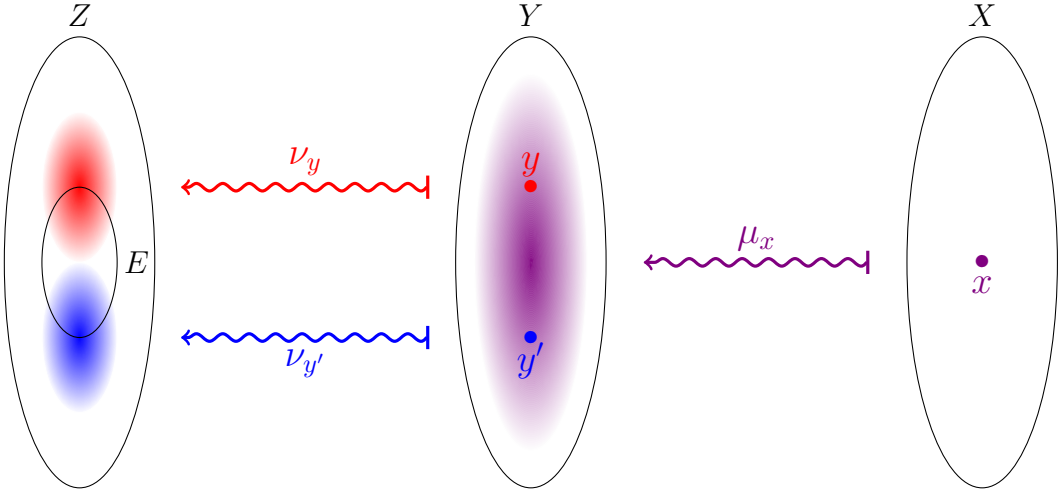
\begin{figure}
\centering
\begin{tikzpicture}[xscale=-1]
\draw (-6,0) ellipse (1cm and 3cm) node[above,yshift=3cm]{$X$};
\node[violet] at (-6,0) {$\bullet$};
\node[violet] at (-6,-0.3) {{\large$x$}};
\draw (0,0) ellipse (1cm and 3cm) node[above,yshift=3cm]{$Y$};
\filldraw[white,even odd rule, inner color=violet, outer color=white,fill opacity=0.5] (0,0) ellipse (0.75cm and 2.5cm);
\node[red] at (0,1) {$\bullet$};
\node[red] at (0,1.3) {{\large$y$}};
\node[blue] at (0,-1) {$\bullet$};
\node[blue] at (0,-1.3) {{\large$y'$}};
\filldraw[white,even odd rule, inner color=red, outer color=white,fill opacity=0.75] (6,1) ellipse (0.5cm and 1cm);
\filldraw[white,even odd rule, inner color=blue, outer color=white,fill opacity=0.75] (6,-1) ellipse (0.5cm and 1cm);
\draw (6,0) ellipse (1cm and 3cm) node[above,yshift=3cm]{$Z$};
\draw (6,0) ellipse (0.5cm and 1cm) node[right,xshift=0.45cm]{$E$};
\draw[very thick,|->,violet,decorate,decoration={snake, amplitude=0.5mm}] (-4.5,0) -- node[above]{{\large$\mu_{x}$}}(-1.5,0);
\draw[very thick,|->,red,decorate,decoration={snake, amplitude=0.5mm}] (1.35,1.0) -- node[above]{{\large$\nu_{y}$}}(4.65,1.0);
\draw[very thick,|->,blue,decorate,decoration={snake, amplitude=0.5mm}] (1.35,-1.0) -- node[below]{{\large$\nu_{y'}$}}(4.65,-1.0);
\end{tikzpicture}
\caption{
$x\in X$ evolves to the probability measure $\mu_{x}$ 
and all the points in $Y$
that have neighborhoods with non-vanishing measure themselves get
sent to probability measures on $Z.$ 
To calculate the overlap with a measurable set $E\subseteq Z,$
the measures are evaluated on $E$ and the results are integrated over
all of $Y.$ 
Only $\nu_{y}$ and $\nu_{y'}$ for two points $y,y'\in Y$ 
have been sketched on the left for ease of visualization.}
\label{fig:continuousstochasticcomposition}
\end{figure}

The composition of stochastic maps is quite subtle for compact Hausdorff
spaces due to the possibility that the underlying spaces need not be
metrizable. Before going into the details, we isolate the issue that occurs
by what we expect the composition to be heuristically based on our
understanding of the finite set case. In this regard, 
let $X,Y,$ and $Z$ be compact Hausdorff spaces
and let $\mu:X\stoch Y$ and $\nu:Y\stoch Z$ be two stochastic maps.
The composition of $\mu$ followed by $\nu$ should be some stochastic map 
$X\stoch Z$ which when applied to a Borel set $E\subseteq Z$
and a point $x\in X$ should give the probability of $x$ evolving 
through the intermediate possibilities in $Y$ that end up at the ``window'' $E$
as sketched in Figure \ref{fig:continuousstochasticcomposition}.

In other words, we sum up $\nu_{y}(E)$ over all $y$ using the measure
$\mu_{x}$ on $Y$ to obtain 
\be
\int_{Y}\nu_{y}(E)\;d\mu_{x}(y). 
\ee
However, this formula assumes that the function 
$Y\ni y\mapsto\nu_{y}(E)$ is measurable for all Borel sets $E\in\mathcal{M}_{Z}.$ 
The validity of this formula is well known if all the spaces in question
are Polish \cite{Gi82} but does not seem to be known for arbitrary
compact Hausdorff spaces. We will prove this result soon, 
but first we prove a crucial result that is interesting in its own right. 

\blem
\label{lem:evaluationofmeasurecontinuous}
Let $Y$ be a compact Hausdorff space 
and let $E\in\mathcal{M}_{Y}$ be a Borel set. Then the evaluation function 
\be
\begin{split}
\mathrm{ProbMeas}(Y)&\xrightarrow{\mathrm{ev}_{E}}\R\\
\mu&\mapsto\mu(E)
\end{split}
\ee
is Borel measurable with respect to the vague topology on 
$\mathrm{ProbMeas}(Y).$ 
In fact, $\mathrm{ev}_{E}$ is lower semi-continuous when $E$ is open
and upper semi-continuous when $E$ is closed. 
\elem

To prove this lemma, we recall 
a useful theorem, 
which we learned through Sengupta's lucid notes \cite{Se09}. 

\bt
[Dynkin's $\pi$-$\l$ Theorem]
\label{thm:Dynkin}
Let $X$ be a set and let $P$ and $L$ be collections of subsets of $X$
satisfying the following conditions. 
\begin{enumerate}[i.]
\setlength{\itemsep}{0pt}
\item
$P$ is closed under finite intersections (such a collection is called a
\emph{\uline{$\pi$-system}}).
\item
$L$ contains the empty set, is closed under taking complements (in $X$), 
and is closed under countable unions of pairwise disjoint elements (such
a collection is called a \emph{\uline{$\l$-system}}). 
\item
$P\subseteq L.$
\end{enumerate}
Then the $\s$-algebra generated by $P$ is contained in $L.$ 
\et

\bprf
See Theorem 3.2 in \cite{Bi95} or the notes \cite{Se09}. 
\eprf

We are indebted to Iddo Ben-Ari for the crucial insight that
Dynkin's $\pi$-$\l$ theorem can be used to give a nice proof
of Lemma \ref{lem:IddosLemma}, which lead to the proof of
Lemma \ref{lem:evaluationofmeasurecontinuous}. 

\bprf[Proof of Lemma \ref{lem:evaluationofmeasurecontinuous}]
The proof will be broken up into several steps.
First, it will be shown that the collection, $P,$ of all open subsets
of $Y$ forms a $\pi$-system. Then it will be shown that the set
\be
L:=
\left\{E\subseteq Y\;:\;\mathrm{ProbMeas}(Y)\xrightarrow{\mathrm{ev}_{E}}[0,1]\text{ is Borel measurable}\right\}
\ee
is a $\l$-system. Finally, it will be shown that $\mathrm{ev}_{U}$ is 
lower semi-continuous (and hence Borel measurable), 
for all open sets $U\subseteq Y$ showing that $P\subseteq L.$  
\begin{enumerate}[i.]
\setlength{\itemsep}{0pt}
\item
The fact that open sets are closed under finite intersections
is part of the definition of a topology. 

\item
For each $\mu\in\mathrm{ProbMeas}(Y),$ 
\be
\mathrm{ev}_{\varnothing}(\mu)=\mu(\varnothing)=0
\ee
so that $\mathrm{ev}_{\varnothing}$ is constant and hence measurable. 
For $E\in L,$ 
\be
\mathrm{ev}_{E^{c}}(\mu)=\mu(E^{c})=\mu(Y)-\mu(E)=(1-\mathrm{ev}_{E})(\mu)
\ee
is the difference of two measurable functions with respect to $\mu$ 
and is therefore measurable with respect to $\mu.$
For $\{E_{n}\}$ a countable collection of pairwise disjoint 
measurable sets in $L$ with $E:=\bigcup_{n}E_{n},$ 
\be
\mathrm{ev}_{E}(\mu)
=\mu\left(\bigcup_{n=1}^{\infty}E_{n}\right)
=\sum_{n=1}^{\infty}\mu(E_{n})
=\left(\sum_{n=1}^{\infty}\mathrm{ev}_{E_{n}}\right)(\mu)
\ee
by additivity of the measure $\mu.$ Hence, since this is a bounded countable
sum of measurable functions, the sum is measurable. 

\item
Let $U\subseteq Y$ be open and non-empty (if $U$ is empty, the result 
was already proved in part ii.). 
By Proposition \ref{prop:semicontinuous}, $\chi_{U}$ is lower semi-continuous.
Hence, for any $a\in\R,$ the preimage of $(a,\infty)$ under $\mathrm{ev}_{U}$
is 
\be
\begin{split}
\mathrm{ev}_{U}^{-1}\big((a,\infty)\big)&=\big\{\mu\in\mathrm{ProbMeas}(Y)\;:\;\mu(U)>a\big\}\\
&=\big\{\mu\in\mathrm{ProbMeas}(Y)\;:\;\int_{Y}\chi_{U}\;d\mu>a\big\}\\
&=\left\{\mu\in\mathrm{ProbMeas}(Y)\;:\;\sup_{\substack{\vf\in C(Y)\\0\le\vf\le\chi_{U}}}\left\{\int_{Y}\vf\;d\mu\right\}>a\right\}
\text{ by Corollary \ref{cor:LSCintegralapprox}}\\
&=\bigcup_{\substack{\vf\in C(Y)\\0\le\vf\le\chi_{U}}}\left\{\mu\in\mathrm{ProbMeas}(Y)\;:\;\int_{Y}\vf\;d\mu>a\right\}.
\end{split}
\ee
By definition of the vague topology, the set in curly brackets in the last equality
is open. Hence, this is a union of open sets and is therefore
an open subset of $\mathrm{ProbMeas}(Y).$ 
Therefore, this shows that $\mathrm{ev}_{U}$ is
lower semi-continuous for all open sets $U$ in $Y.$ 
An analogous argument using infima shows that $\mathrm{ev}_{K}$ 
is upper semi-continuous for all closed sets $K$ in $Y.$ 
\end{enumerate}
Thus, since all the conditions of Dynkin's $\pi$-$\l$ theorem are 
satisfied, the $\s$-algebra generated by $P$
which is the Borel $\s$-algebra, is contained in $L.$ 
\eprf

\br
\label{rmk:LSCversusC}
Although $\mathrm{ev}_{U}$ was shown to be lower semi-continuous
for all open sets $U,$ it is in general \emph{not} continuous. 
For an explicit example, 
enumerate the rationals in $[0,1]$ by some bijective function 
$q:\N\to\Q\cap[0,1]$ and set
\be
U:=\left(\bigcup_{n=1}^{\infty}\left(q_{n}-\frac{\e}{2^{n+1}},q_{n}+\frac{\e}{2^{n+1}}\right)\right)\cap[0,1]
\ee
for some $\e\in(0,1).$ Then, $U$ is open, 
contains all of the rationals in $[0,1],$
and has Lebesgue measure $\e<1$ so that it does not contain all of $[0,1].$ 
Fix $x_{\infty}\in[0,1]\setminus U.$ 
Since $\Q\cap[0,1]$ is dense in $[0,1],$ there exists a sequence
$x:\N\to\Q\cap[0,1]$ converging to $x_{\infty}$ in $[0,1].$ 
Then, the corresponding sequence of Dirac measures 
$n\mapsto\de_{x_{n}}$ converges to $\de_{x_{\infty}}$ in the vague topology
but
\be
\de_{x_{n}}(U)=1
\aand
\de_{x_{n}}\big([0,1]\setminus U\big)=0
\qquad\forall\;n\in\N
\ee
while
\be
\de_{x_{\infty}}(U)=0
\aand
\de_{x_{\infty}}\big([0,1]\setminus U\big)=1
\ee
showing that $\mathrm{ev}_{U}$ is not continuous at $\de_{x_{\infty}}.$ 
In fact, checking that $\mathrm{ev}_{U}$ is lower semi-continuous 
at $\de_{x_{\infty}}$ directly follows from the definition since
$\mathrm{ev}_{U}(\mu)\ge\mathrm{ev}_{U}(\de_{x_{\infty}})-\e$
reads
$\mu(U)\ge0-\e=-\e,$ which is true for all $\e>0$ and probability measures $\mu.$ 
However, one runs into trouble checking if $\mathrm{ev}_{U}$ is upper
semi-continuous for the reason illustrated above. 
\er
 
\blem
\label{lem:evaluationmeasurable}
Let $\mu:X\stoch Y$ be a stochastic map as in 
Definition \ref{defn:continuousstochasticprocess}. 
For every Borel set $E\in\mathcal{M}_{Y},$
the function 
\be
X\ni x\mapsto\mu_{x}(E)
\ee
is Borel measurable. In fact, if $E$ is open (closed), then
this function is lower (upper) semi-continuous. 
\elem

\bprf
This follows from Lemma \ref{lem:evaluationofmeasurecontinuous}
and the fact that this function is the composition
\be
X\xrightarrow{\mu}\mathrm{ProbMeas}(Y)\xrightarrow{\mathrm{ev}_{E}}\R
\ee
of Borel measurable functions. If $E$ is open (closed), then this function is the 
composition of a continuous function followed by a 
lower (upper) semi-continuous
function, which is still lower (upper) semi-continuous 
(this is immediate from the definitions).
\eprf

Note that this shows, in particular, that a stochastic map gives rise to a 
Markov kernel. The converse is not true as an arbitrary Markov kernel takes
no account of continuity. Hence, although the composition of Markov kernels is 
a Markov kernel,
the requirement of being a stochastic map imposes additional restrictions.

\bn
\label{prop:compositionstochastic}
Let $X,Y,$ and $Z$ be compact Hausdorff spaces
and let $\mu:X\stoch Y$ and $\nu:Y\stoch Z$ be two stochastic maps.
For each $x\in X,$ the assignment 
\be
\label{eq:compositionofstochasticmaps}
\mathcal{M}_{Z}\ni E\mapsto(\nu\circ\mu)_{x}(E):=\int_{Y}\nu_{y}(E)\;d\mu_{x}(y),
\ee
is a regular probability measure on $Z.$ Furthermore, the assignment
\be
\begin{split}
X&\to\mathrm{ProbMeas}(Z)\\
x&\mapsto(\nu\circ\mu)_{x},
\end{split}
\ee
is a stochastic map. This stochastic map is called the \emph{\uline{composition}}
of $\mu$ followed by $\nu$ and is denoted by $\nu\circ\mu.$
\en

\bprf
The proof will be split into three parts to prove the following claims: 
i. $(\nu\circ\mu)_{x}$ is a Borel probability measure on $Z,$
ii. it is regular, and 
iii. the assignment  $X\ni x\mapsto(\nu\circ\mu)_{x}$ is a stochastic map $X\stoch Z.$
\begin{enumerate}[i.]
\setlength{\itemsep}{0pt}
\item 
First, 
\be
(\nu\circ\mu)_{x}(Z)=\int_{Y}\underbrace{\nu_{y}(Z)}_{1}\;d\mu_{x}(y)
=\int_{Y}d\mu_{x}(y)
=\mu_{x}(Y)
=1.
\ee
Now, fix a Borel set $E\in\mathcal{M}_{Z}.$
Since $Y\ni y\mapsto\nu_{y}(E)$ is Borel measurable (by Lemma 
\ref{lem:evaluationmeasurable})
and since $\nu_{y}(E)$ is non-negative for all Borel sets $E$
in $Z,$ its integral with respect to a positive measure is non-negative. 
We now check that $(\nu\circ\mu)_{x}$ is additive on (at most)
countable measurable partitions.
Let $\{E_{i}\}_{i\in I}$ be a countable
measurable partition of $E.$ Then
\be
\begin{split}
(\nu\circ\mu)_{x}(E)&=\int_{Y}\nu_{y}(E)\;d\mu_{x}(y)\quad\text{ by definition of $(\nu\circ\mu)_{x}$}\\
&=\int_{Y}\sum_{i\in I}\nu_{y}(E_{i})\;d\mu_{x}(y)\quad\text{ since $\nu_{y}$ is a measure}\\
&=\sum_{i\in I}\int_{Y}\nu_{y}(E_{i})\;d\mu_{x}(y)\quad\text{ by the Monotone Convergence Theorem}\\
&=\sum_{i\in I}(\nu\circ\mu)_{x}(E_{i})\quad\text{ by definition of $(\nu\circ\mu)_{x}$},
\end{split}
\ee
In the third equality, 
$\nu_{y}(E_{i})$ is non-negative and measurable for all $y\in Y$
so the partial sums form an increasing sequence of measurable functions. 
Hence, the Monotone Convergence Theorem
allows one to exchange the sum operation with the integral. 
Therefore, $(\nu\circ\mu)_{x}$ is a Borel probability measure. 

\item
To prove regularity of the measure $(\nu\circ\mu)_{x},$
fix $E\in\mathcal{M}_{Z}$ and $\e>0.$ The goal is to find 
an open set $U\subseteq Z$ and a compact set $K\subseteq Z$ 
such that $K\subseteq E\subseteq U$ and 
\be
(\nu\circ\mu)_{x}(U)-(\nu\circ\mu)_{x}(E)<\e
\aand
(\nu\circ\mu)_{x}(E)-(\nu\circ\mu)_{x}(K)<\e.
\ee
Note, that it suffices to show that 
(by a single application of the triangle inequality)
\be
(\nu\circ\mu)_{x}(U)-(\nu\circ\mu)_{x}(K)<\e.
\ee
Fix $y\in Y.$ Since $\nu_{y}$ is regular, there exist 
an open set $U_{y}\subseteq Z$ and a closed set $K_{y}\subseteq Z$
such that $K_{y}\subseteq E\subseteq U_{y}$ and 
\be
\label{eq:regularityUK}
\nu_{y}(U_{y})-\nu_{y}(K_{y})<\frac{\e}{3}. 
\ee
Since $Y\ni y'\mapsto\nu_{y'}(U_{y})$ is lower semi-continuous 
and $Y\ni y'\mapsto\nu_{y'}(K_{y})$ is upper semi-continuous 
by Lemma \ref{lem:evaluationmeasurable}, there exists an open
neighborhood $V_{y}\subseteq Y$ containing $y$ such that 
\be
\label{eq:regularityUKlocal}
\nu_{y}(U_{y})-\nu_{y'}(U_{y})\le\frac{\e}{3}
\aand
\nu_{y'}(K_{y})-\nu_{y}(K_{y})\le\frac{\e}{3}
\qquad\forall\;y'\in V_{y}.
\ee
The collection of such open and closed sets can be obtained for every $y\in Y.$
Hence, $\{V_{y}\}_{y\in Y}$ forms an open cover of $Y.$ 
By compactness of $Y,$ there exists a finite collection
$\{y_1,\dots,y_n\}$ of points in $Y$ for which 
$\{V_{y_1}, \dots,V_{y_n}\}$ covers $Y.$ 
Set 
\be
\label{eq:finiteUandK}
U:=\bigcap_{i=1}^{n}U_{y_i}
\aand
K:=\bigcup_{i=1}^{n}K_{y_i}. 
\ee
Since these intersections and unions are finite, 
$U$ is open and $K$ is closed in $Z$.
Furthermore, they satisfy $K\subseteq E\subseteq U.$ 
Now, fix $y\in Y.$ There exists some $i\in\ov n$ such that 
$y\in V_{y_{i}}$ so that 
\be
\begin{split}
\big|\nu_{y}(U)-\nu_{y}(K)\big|&\le
\big|\nu_{y}(U_{y_{i}})-\nu_{y}(K_{y_{i}})\big|\quad\text{by (\ref{eq:finiteUandK}) and sub-additivity of $\nu_{y}$}\\
&=\big|\nu_{y}(U_{y_{i}})-\nu_{y_i}(U_{y_{i}})+\nu_{y_i}(U_{y_{i}})-\nu_{y_{i}}(K_{y_{i}})+\nu_{y_{i}}(K_{y_{i}})-\nu_{y}(K_{y_{i}})\big|\\
&\le\big|\nu_{y}(U_{y_{i}})-\nu_{y_i}(U_{y_{i}})\big|+\big|\nu_{y_i}(U_{y_{i}})-\nu_{y_{i}}(K_{y_{i}})\big|+\big|\nu_{y_{i}}(K_{y_{i}})-\nu_{y}(K_{y_{i}})\big|\\
&<\e\quad\text{by (\ref{eq:regularityUK}) and (\ref{eq:regularityUKlocal})}.
\end{split}
\ee
Since this inequality is true for every $y\in Y,$
\be
(\nu\circ\mu)_{x}(U)-(\nu\circ\mu)_{x}(K)
=\int_{Y}\Big(\nu_{y}(U)-\nu_{y}(K)\Big)\;d\mu_{x}(y)
<\int_{Y}\e\;d\mu_{x}
=\e
\ee
because $\mu_{x}$ is a probability measure. 
This proves that $(\nu\circ\mu)_{x}$ is regular for arbitrary $x\in X.$ 

\item
Now, we check that the function $X\ni x\mapsto(\nu\circ\mu)_{x}$ is continuous 
so that $\nu\circ\mu:X\stoch Z$ is a stochastic map.  
In this regard, fix an arbitrary $\vf\in C(Z)$ and a net $x:\Theta\to X$ converging
to $\lim x.$ The goal is to show that 
\be
\lim_{\q\in\Theta}\int_{Z}\vf\;d(\nu\circ\mu)_{x_{\q}}=
\int_{Z}\vf\;d(\nu\circ\mu)_{\lim x}.
\ee
By the Jordan decomposition theorem, it suffices to assume that $\vf$ is positive.
Since $\vf$ is continuous, by Lemma 
\ref{lem:contcompactsequenceapprox} there exists a sequence of 
simple functions $0\le s_1\le s_2\le \cdots\le\vf$ of the form 
\be
s_{n}=\sum_{i_{n}}s_{n,i_{n}}\chi_{E_{n,i_n}}
\ee
defined on Borel sets $\{E_{n,i_n}\}$ in $Z$ satisfying all of the
conclusions of that Lemma. Hence
\be
\label{eq:provingcompositioncontinuous}
\begin{split}
\lim_{\q\in\Theta}\int_{Z}\vf\;d(\nu\circ\mu)_{x_\q}&=
\lim_{\q\in\Theta}\lim_{n\to\infty}\sum_{i_n}s_{n,i_n}(\nu\circ\mu)_{x_\q}(E_{n,i_n})\quad\text{by definition of the integral}\\
&=\lim_{\q\in\Theta}\lim_{n\to\infty}\sum_{i_n}s_{n,i_n}\int_{Y}\nu_{y}(E_{n,i_n})\;d\mu_{x_\q}\quad\text{by definition of $\nu\circ\mu$}\\
&=\lim_{\q\in\Theta}\lim_{n\to\infty}\int_{Y}\sum_{i_n}s_{n,i_n}\nu_{y}(E_{n,i_n})\;d\mu_{x_\q}\quad\text{by linearity of the integral.}
\end{split}
\ee
Now, to verify that the limit can be interchanged with the integral, note that
by the properties of the simple functions constructed above 
and since $\nu_{y}$ is a measure, 
\be
\begin{split}
\sum_{i_n}s_{n,i_n}\nu_{y}(E_{n,i_{n}})&=
\sum_{i_n}s_{n,i_n}\nu_{y}(E_{n+1,2i_{n}})
+\sum_{i_n}s_{n,i_n}\nu_{y}(E_{n+1,2i_{n}+1})\quad\text{by (\ref{eq:splittingEtotwoparts})}\\
&=\sum_{\substack{i_{n+1}\\\text{even}}}s_{n,\frac{i_{n+1}}{2}}\nu_{y}\left(E_{n+1,i_{n+1}}\right)
+\sum_{\substack{i_{n+1}\\\text{odd}}}s_{n,\frac{i_{n+1}-1}{2}}\nu_{y}\left(E_{n+1,i_{n+1}}\right)\\
&\le\sum_{\substack{i_{n+1}\\\text{even}}}s_{n+1,i_{n+1}}\nu_{y}\left(E_{n+1,i_{n+1}}\right)
+\sum_{\substack{i_{n+1}\\\text{odd}}}s_{n,i_{n+1}-1}\nu_{y}\left(E_{n+1,i_{n+1}}\right)\quad\text{by (\ref{eq:breakingsimpleintotwo})}\\
&=\sum_{i_{n+1}}s_{n+1,i_{n+1}}\nu_{y}(E_{n+1,i_{n+1}})
\end{split}
\ee
for all $y\in Y.$ Therefore, the sequence of functions
\be
\N\times Y\ni(n,y)\mapsto\sum_{i_n}s_{n,i_n}\nu_{y}(E_{n,i_{n}})
\ee
on $Y$ is a monotonic increasing sequence of Borel 
measurable functions. Hence, continuing on from 
equation (\ref{eq:provingcompositioncontinuous}), 
\be
\label{eq:provingcompositioncontinuous2}
\begin{split}
\lim_{\q\in\Theta}\int_{Z}\vf\;d(\nu\circ\mu)_{x_\q}&=
\lim_{\q\in\Theta}\int_{Y}\lim_{n\to\infty}\sum_{i_n}s_{n,i_n}\nu_{y}(E_{n,i_n})\;d\mu_{x_\q}\\
&\qquad\quad\text{by the Monotone Convergence Theorem}\\
&=\lim_{\q\in\Theta}\int_{Y}\left(\int_{Z}\vf\;d\nu_{y}\right)\;d\mu_{x_\q}\quad\text{by definition of the integral.}
\end{split}
\ee
To continue with this calculation, note that by assumption, the assignment
$y\mapsto \nu_{y}$ is continuous with respect to the vague topology on
$\mathrm{ProbMeas}(Z).$ This means that
for any continuous function such as $\vf\in C(Z)$ above, 
$y\mapsto\int_{Z}\vf\;d\nu_{y}$ is a continuous real-valued function on $Y.$ 
Hence, following from (\ref{eq:provingcompositioncontinuous2}),
\be
\lim_{\q\in\Theta}\int_{Z}\vf\;d(\nu\circ\mu)_{x_\q}
=\int_{Y}\left(\int_{Z}\vf\;d\nu_{y}\right)\;d\mu_{\lim x}
\ee
since $\mu$ is continuous with respect to the vague topology on 
$\mathrm{ProbMeas}(Y).$ 
By applying an analogous procedure 
(except without taking any limit in $\Theta$), 
one similarly concludes that 
\be
\int_{Z}\vf\;d(\nu\circ\mu)_{\lim x}=\int_{Y}\left(\int_{Z}\vf\;d\nu_{y}\right)\;d\mu_{\lim x}.
\ee
Since these two results are equal, continuity of $\nu\circ\mu$ has been
exhibited. 
\end{enumerate}
Therefore, $\nu\circ\mu:X\stoch Z$ is a stochastic map. 
\eprf

From this fact, it immediately follows that the identity stochastic
map is indeed an identity for the composition of stochastic maps. 

\bn
\label{prop:cHausStochAssociative}
The composition of stochastic maps is associative. 
\en

\bprf
Let $\mu:X\stoch Y,\nu:Y\stoch Z,$ and $\z:Z\stoch W$ be stochastic maps
of compact Hausdorff spaces.
The goal is to show that 
$\big(\z\circ(\nu\circ\mu)\big)_{x}(E)=\big((\z\circ\nu)\circ\mu\big)_{x}(E)$
for all Borel sets $E\in\mathcal{M}_{W}$ and points $x\in X.$ 
By expanding out the definitions, it must be shown that
\be
\label{eq:associativitygoal}
\int_{Z}\z_{z}(E)\;d(\nu\circ\mu)_{x}(z)=\int_{Y}\left(\int_{Z}\z_{z}(E)\;d\nu_{y}(z)\right)d\mu_{x}(y).
\ee
Contrary to what one might naively think, Fubini's Theorem does not apply 
here. Instead, since $Z\ni z\mapsto\z_{z}(E)$ is Borel measurable
and bounded, 
by Lemma \ref{lem:contcompactsequenceapprox}, there exists a 
sequence $s:\N\to \R^{Z}$ of simple functions satisfying
\be
0\le s_1(z)\le s_2(z)\le s_3(z)\le\cdots\le\z_{z}(E)
\aand
\lim_{n\to\infty}s_{n}(z)=\z_{z}(E)\qquad\forall\;z\in Z. 
\ee
and all the other properties of that Lemma.
Hence, by a completely analogous 
calculation as in (\ref{eq:provingcompositioncontinuous})--(\ref{eq:provingcompositioncontinuous2}) but without a net, 
the equality in (\ref{eq:associativitygoal}) is obtained. 
This proves that composition of stochastic maps is associative. 
\eprf

In conclusion, we have proved the first part of the following theorem. 

\bt
\label{thm:cHStoch}
The collection of compact Hausdorff spaces together with
morphisms as stochastic maps forms a category, 
denoted by $\cHStoch.$ In addition, the Dirac measure
furnishes a faithful (but not full) 
functor $\de:\cH\to\cHStoch$ sending a space $X$ to $X$
and a continuous function $f:X\to Y$ to $\de_{f}:X\stoch Y$ 
as in Example \ref{ex:diracstochastic}. 
\et

\bprf
The faithfulness of $\de$ follows immediately from the Hausdorff condition.
More precisely, seeking to prove the contrapositive, 
suppose that $f,g:X\to Y$ were two
different functions. Then there exists an $x\in X$ such that $f(x)\ne g(x).$ 
Since $Y$ is Hausdorff, there exist open sets $U$ and $V$ in $Y$ such that
$x\in U, y\in V,$ and $U\cap V=\varnothing.$ Then $\de_{f(x)}(U)=1$ while
$\de_{g(x)}(U)=0.$ Hence $\de_{f}\ne\de_{g}.$ 
\eprf

\subsection{From continuous probability theory to algebra}
\label{sec:stochasticCfunctor}

We have already associated to every compact Hausdorff
space a commutative $C^*$-algebra by taking 
continuous functions on that space. 
To every stochastic map, we can
also associate a canonical positive map of $C^*$-algebras. 

\bn
\label{prop:curlyEcodomain}
Let $\mu:X\stoch Y$ be a stochastic map between compact
Hausdorff spaces. For every continuous function $\vf:Y\to\C,$ 
the function $C^{\mu}(\vf):X\to\C$ defined by 
\be
\label{eq:curlyEonstochasticmaps}
X\ni x\mapsto\int_{Y}\vf\;d\mu_{x}
\ee
is a continuous function on $X.$ 
Furthermore, $C^{\mu}:C(Y)\to C(X)$ is a (unital) positive map. 
\en

\bprf
Linearity and positivity of $C^{\mu}$ follow immediately from the definitions
and properties of a probability measure. 
To prove that $C^{\mu}(\vf)$ is a continuous function for any 
continuous function $\vf\in C(Y),$ let
$x:\Theta\to X$ be a net converging to $\lim x\in X.$ 
By definition of the vague topology on $\mathrm{ProbMeas}(Y)$
and because $\mu:X\to\mathrm{ProbMeas}(Y)$ is continuous, 
\be
\lim_{\q\in\Theta}\big(C^{\mu}(\vf)(x_{\q})\big)
=\lim_{\q\in\Theta}\int_{Y}\vf\;d\mu_{x_{\q}}
=\int_{Y}\vf\;d\mu_{\lim x}
=C^{\mu}(\vf)(\lim x)
\ee
showing that $C^{\mu}(\vf)$ is continuous. 
\eprf

\bt
\label{thm:stochasticcontinuousfunctor}
The assignment 
\be
\begin{split}
\cHStoch^{\op}&\xrightarrow{C}\cCAlgPos\\
X&\mapsto C(X)\\
(\mu:X\stoch Y)&\mapsto \big(C^{\mu}:C(Y)\stoch C(X)\big)
\end{split}
\ee
is a functor. Moreover, the diagram 
\be
\xy0;/r.25pc/:
(-17.5,7.5)*+{\cH^{\op}}="1";
(-17.5,-7.5)*+{\cHStoch^{\op}}="2";
(17.5,7.5)*+{\cCAlg}="3";
(17.5,-7.5)*+{\cCAlgPos}="4";
{\ar"1";"3"^{C}};
{\ar"2";"4"_{C}};
{\ar"1";"2"_{\de}};
{\ar"3";"4"};
\endxy
\ee
commutes (the functor on the right is the inclusion since every 
{$*$-homomorphism} is positive). 
\et

\bprf
The fact that the codomain of $C$ is as claimed
follows from Proposition \ref{prop:curlyEcodomain}. 
 
If $\id_{X}:X\stoch X$ is the identity stochastic map, then
\be
\big(C^{\id_{X}}(\vf)\big)(x)
=\int_{X}\vf\;d\de_{x}
=\vf(x)
\ee
for all $x\in X$ and for all $\vf\in C(X).$  Hence $C^{\id_{X}}(\vf)=\vf$
for all $\vf\in C(X).$ Thus $C$ preserves identities. 

To see that $C$ preserves the composition, let
$\mu:X\stoch Y$ and $\nu:Y\stoch Z$ be two stochastic maps. 
The goal is to show
that $C^{\nu\circ\mu}=C^{\mu}\circ C^{\nu}.$ 
The proof will be similar to the proof that composition in 
$\cHStoch$ is associative. To 
do this, fix $\vf\in C(Z).$ By the Jordan decomposition theorem, 
it suffices to assume that $\vf(Z)\subseteq\R_{\ge0}$ since
any measurable complex-valued function can be decomposed
into its real and imaginary parts and each of those into positive and negative
parts. 
Since $\vf$ is continuous and $Z$ is compact, 
let $s$ be a sequence of simple functions as in Lemma 
\ref{lem:contcompactsequenceapprox}. 
Then, 
\be
\begin{split}
\big(C^{\nu\circ\mu}(\vf)\big)(x)&=\int_{Z}\vf\;d(\nu\circ\mu)_{x}\quad\text{by definition of $C^{\nu\circ\mu}$}\\
&=\lim_{n\to\infty}\sum_{i_{n}}s_{n,i_{n}}(\nu\circ\mu)_{x}(E_{n,i_n})\quad\text{by definition of the Lebesgue integral}\\
&=\lim_{n\to\infty}\sum_{i_{n}}s_{n,i_{n}}\int_{Y}\nu_{y}(E_{n,i_n})\;d\mu_{x}(y)\quad\text{by definition of $(\nu\circ\mu)_{x}$}\\
&=\lim_{n\to\infty}\int_{Y}\sum_{i_{n}}s_{n,i_{n}}\nu_{y}(E_{n,i_n})\;d\mu_{x}(y)
\quad\text{by (finite) linearity of the integral}\\
&=\int_{Y}\lim_{n\to\infty}\sum_{i_{n}}s_{n,i_{n}}\nu_{y}(E_{n,i_n})\;d\mu_{x}(y)\quad\text{by Monotone Convergence Theorem}\\
&=\int_{Y}\left(\int_{Z}\vf\;d\nu_{y}\right)d\mu_{x}(y)\quad\text{by definition of the Lebesgue integral}\\
&=\int_{Y}\big(C^{\nu}(\vf)\big)(y)\;d\mu_{x}(y)\quad\text{by definition of $C^{\nu}$}\\
&=\Big(C^{\mu}\big(C^{\nu}(\vf)\big)\Big)(x)\quad\text{by definition of $C^{\mu}$}
\end{split}
\ee
for all $x\in X.$ The Monotone Convergence Theorem can be applied
in this calculation due to Lemma \ref{lem:contcompactsequenceapprox}
and the way the simple functions were constructed in that Lemma
(recall step iii. in the proof of Proposition \ref{prop:compositionstochastic}).
Since $\vf$ was arbitrary, this shows that 
$C^{\nu\circ\mu}=C^{\mu}\circ C^{\nu}.$

Finally, if $f:X\to Y$ is a continuous function, let $\vf\in C(Y)$ and $x\in X.$ 
Then 
\be
\big(C^f(\vf)\big)(x)=\int_{Y}\vf\;d\de_{f(x)}=\vf\big(f(x)\big)
\ee
showing that the diagram in the statement of the theorem commutes. 
\eprf

\subsection{The stochastic Gelfand spectrum functor}
\label{sec:stochasticGelfandspectrumfunctor}

In this section, we prove that the functor 
$C:\cHStoch^{\op}\to\cCAlgPos$ 
constructed in Section \ref{sec:stochasticCfunctor} is fully faithful. 
To do this, we review the 
Riesz-Markov-Kakutani Representation Theorem.
Then we use these results to construct an analogue of the spectrum 
functor for positive maps on arbitrary commutative $C^*$-algebras. 

\bt
[Riesz-Markov-Kakutani (RMK) Representation Theorem]
\label{thm:Riesz}
Let $Y$ be a compact Hausdorff space and 
let $I:C(Y)\stoch\C$ be a state on $C(Y).$ 
Then, there exists a unique regular probability measure
$\mu\in\mathrm{ProbMeas}(Y)$ such that 
\be
I(\vf)=\int_{Y}\vf\;d\mu\qquad\forall\;\vf\in C(Y).
\ee
Furthermore, this measure $\mu$ satisfies 
\be
\mu(U)=\sup\big\{ I(\vf)\;:\;\vf\in C(Y)
\text{ and }0\le\vf\le\chi_{U}\big\}
\ee
for all open sets $U\subseteq Y$ and 
\be
\mu(K)=\inf\big\{I(\vf)\;:\;\vf\in C(Y)
\text{ and }\vf\ge\chi_{K}\big\}
\ee
for all closed (and hence compact) sets $K\subseteq Y.$ 
\et

Note that by our definition, a state is a positive linear functional that
is unit-preserving. This is why the measure $\mu$ 
in the statement of the theorem is a probability measure. 
We will not prove this theorem in its entirety, 
but an explicit construction of the measure
is useful and provides a close comparison to the finite-dimensional case. 
Namely, when $X$ is finite, $C(X)=\C^{X}$ and it makes sense 
to assign $\mu(E):=I(\chi_{E})$ for any finite subset $E\subseteq X$ 
as we have done in Section \ref{sec:stochasticmatrices}. 
This is because the characteristic function $\chi_{E}$ is continuous on a finite set. 
Hence, heuristically, given a measurable set $E\in\mathcal{M}_{Y},$
one would expect $\mu(E)$ to be given by $I(\chi_{E}).$
However, $\chi_{E}$ is not continuous and so this expression simply does not make 
sense because $\chi_{E}$ is not in the domain of $I.$ 
Nevertheless, we know by regularity that $E$ can be approximated using 
open and closed sets. One can use this fact to first extend $I$ to be a positive linear
functional on the larger space of lower semi-continuous functions
as is done in \cite{Ta09_245B12}. 

\bprf
[Sketch of proof]
First, let $U\subseteq Y$ be open and let $K\subseteq Y$ be closed. 
Set 
\be
\mu(U):=\sup_{\substack{\vf\in C(Y)\\0\le\vf\le\chi_{U}}}\big\{I(\vf)\big\}
\aand
\mu(K):=\inf_{\substack{\vf\in C(Y)\\0\le\chi_{K}\le\vf}}\big\{I(\vf)\big\}.
\ee
For an arbitrary \emph{subset} $E$ of $Y,$ set 
\be
\mu^{+}(E)=\inf_{\substack{\text{$U$ open}\\E\subseteq U}}\big\{\mu(U)\big\}
\aand
\mu^{-}(E)=\sup_{\substack{\text{$K$ closed}\\K\subseteq E}}\big\{\mu(K)\big\}.
\ee
It is a fact that for all Borel subsets $E\in\mathcal{M}_{Y},$
\be
\mu^{+}(E)=\mu^{-}(E)
\ee
and we therefore denote this assignment by $\mu.$ It is also a fact that
$\mu$ is the unique regular probability measure on $Y$ satisfying the
condition in the statement of the theorem. 
See Theorem 7.2 in \cite{Fo07}, Theorem 8 in 
\cite{Ta09_245B12}, or Theorem 2.14 in \cite{Ru87} for proofs of these claims. 
\eprf

\br
In general, the assignments $\mu^{\pm}$ on arbitrary subsets $E$ of $Y$
given by 
\be
\mu^{+}(E):=\inf_{\substack{\text{$U$ open}\\E\subseteq U}}
\left\{\sup_{\substack{\vf\in C(Y)\\0\le\vf\le\chi_{U}}}\big\{I(\vf)\big\}
\right\}
\ee
and
\be
\mu^{-}(E):=\sup_{\substack{\text{$K$ closed}\\K\subseteq E}}
\left\{\inf_{\substack{\vf\in C(Y)\\0\le\chi_{K}\le\vf}}\big\{I(\vf)\big\}
\right\}
\ee
are not the same. These are known as \emph{outer} and \emph{inner}
measures associated to $I.$ The collection of subsets $E$ for which 
$\mu^{+}(E)=\mu^{-}(E)$ contains the Borel sets of $Y.$
In general, however, the resulting $\s$-algebra is larger than the
Borel $\s$-algebra (see Theorem 2.14 in \cite{Ru87}). Nevertheless, 
$\mu$ restricted to the Borel $\s$-algebra associated to $Y$
provides $Y$ with a regular probability measure. 
This is the measure constructed in the proof of Theorem \ref{thm:Riesz}. 
\er

\bn
$\cHStoch^{\op}\xrightarrow{C}\cCAlgPos$ is fully faithful.
More explicitly, 
let $X$ and $Y$ be compact Hausdorff spaces and let 
$f:C(Y)\stoch C(X)$ be a positive (unital) map.
Then there exists a unique stochastic map 
$\mu:X\stoch Y$ such that 
\be
\label{eq:positivemapstochasticprocess}
\big(f(\vf)\big)(x)=\int_{Y}\vf\;d\mu_{x}
\qquad\forall\;x\in X\qquad\forall\;\vf\in C(Y). 
\ee
\en

\bprf
First note that for every $\vf\in C(Y)$ and $x\in X,$ 
\be
\big(f(\vf)\big)(x)=(\ev_{x}\circ f)(\vf),
\ee
where $\ev_{x}:C(X)\to\C$ is the evaluation at $x$ function
(recall Example \ref{ex:CXalgebra}),
which is a {$*$-homomorphism} and therefore a positive map.
Hence, the composition
$C(Y)\;\xy0;/r.25pc/:\ar@{~>}(-3,0);(3,0)^{f}\endxy\;C(X)\xrightarrow{\ev_{x}}\C$
is a state on $C(Y).$ 
By the RMK Representation Theorem, 
there exists a unique regular probability measure
$\mu_{x}\in\mathrm{ProbMeas}(Y)$ satisfying 
(\ref{eq:positivemapstochasticprocess}). 
The only thing left to check is that the assignment
$X\ni x\mapsto\mu_{x}\in\mathrm{ProbMeas}(Y)$
is continuous. However, this immediately follows from the
fact that $\vf(f)$ is continuous since if 
$x:\Theta\to X$ is a net converging to $\lim x,$ 
then 
\be
\lim_{\q\in\Theta}\int_{Y}\vf\;d\mu_{x_\q}=\lim_{\q\in\Theta}\big(f(\vf)\big)(x_\q)
=\big(f(\vf)\big)(\lim x)
=\int_{Y}\vf\;d\mu_{\lim x}.
\ee
Hence, $C$ is fully faithful. 
\eprf

Thus, the RMK Representation Theorem precisely encodes the fact that 
the functor $C:\cHStoch^{\op}\to\cCAlgPos$ is fully faithful while  
the Gelfand-Naimark Theorem shows that this functor is essentially
surjective. Together, these results hold if and only if $C$ is part of
adjoint equivalence of categories. 
However, this result is somewhat formal. It is useful to have 
an explicit functor $\s:\cCAlgPos\to\cHStoch^{\op}$ that (1) requires no
reference to any underlying space structure and (2) agrees
with the usual spectrum functor from the ordinary (commutative) Gelfand-Naimark
correspondence when restricted to {$*$-homomorphisms} of $C^*$-algebras.

\begin{figure}
\centering
  \begin{tikzpicture}[scale=0.3333]
  \pgfmathsetmacro{\R}{5}
  \pgfmathsetmacro{\r}{2.887}
  \pgfmathsetmacro{\s}{30}
  \pgfmathsetmacro{\d}{2.887}
  \node at (0,0) {$\bullet$};
  \node at (1.9,0.6) {$\chi\circ f$};
  \node[blue] at (-6.6,-2.5) {$\sigma^{f}_{\chi}$};
  \foreach \n in {0,1,2}
  {
  \draw[blue,ultra thick] ({\d*cos(90+120*\n)},{\d*sin(90+120*\n)}) ++({\s+\n*120}:{\r}) arc ({\s+\n*120}:{\s+2*(90-\s)+\n*120}:{\r});
  \fill[blue,fill opacity=0.2] ({\d*cos(90+120*\n)},{\d*sin(90+120*\n)}) ++({\s+\n*120}:{\r}) arc ({\s+\n*120}:{\s+2*(90-\s)+\n*120}:{\r}) -- cycle;
  }
  \foreach \n in {0,1,2}
  {
  \draw ({\R*cos(120+\n*120)},{\R*sin(120+\n*120)}) -- ({\R*cos(60+(\n+1)*120)},{\R*sin(60+(\n+1)*120)});
  \fill[blue,fill opacity=0.2] 
  (0,0) --
  ({\R*cos(60+\n*120)},{\R*sin(60+\n*120)}) -- ({\R*cos(120+\n*120)},{\R*sin(120+\n*120)}) -- cycle;
  }
  \end{tikzpicture}
  \phantom{wxyz}
\caption{A character provides a Dirac point measure on the 
spectrum. A state is ``smeared'' onto the spectrum in terms of a 
probability measure representing the associated ``weight'' of that state
on the characters in the spectrum.}
\label{fig:stateconvexspectrum}
\end{figure}
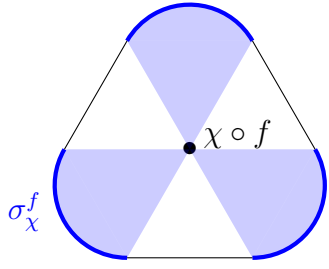

We have already constructed the functor $\s$ at the level of objects, 
namely, given any commutative $C^*$-algebra $\mathcal{A},$
take $\s(\mA)$ to be the spectrum of $\mA$ (see Definition \ref{def:spectrum}). 
Given a {$*$-homomorphism} of commutative $C^*$-algebras 
$f:\mB\to\mA,$ the associated map $\s(f)$ assigns to each
character $\chi\in\s(\mA)$ the linear functional $\chi\circ f:\mB\to\C.$
This linear functional is a character precisely because $f$ and $\chi$ are both
(non-vanishing) {$*$-homomorphisms.} Hence, when $f$ is simply a positive map, 
$\chi\circ f$ is no longer a character in general. Nevertheless, from these data, 
namely a positive map $f:\mB\stoch\mA$ and a character $\chi\in\s(\mA),$
we must construct a canonical regular probability measure $\s^{f}_{\chi}$
on $\s(\mB)$ in such a way so that this assignment 
$\s^{f}:\s(A)\to\mathrm{ProbMeas}\big(\s(\mB)\big)$ is a 
stochastic map. The idea should be that we smear $\chi\circ f$ onto $\s(\mB)$ 
in some way that captures the information of the state $\chi\circ f$
as depicted in Figure \ref{fig:stateconvexspectrum}. 

This can be done formally using the Gelfand transform
and the fact that $C$ is fully faithful. Namely, we can fill in the 
diagram 
\be
\xy0;/r.25pc/:
(-17.5,7.5)*+{\mB}="B";
(17.5,7.5)*+{\mA}="A";
(-17.5,-7.5)*+{C\big(\s(\mB)\big)}="CsB";
(17.5,-7.5)*+{C\big(\s(\mA)\big)}="CsA";
{\ar@{~>}^{f}"B";"A"};
{\ar"CsB";"B"^{\G_{\mB}^{-1}}};
{\ar"A";"CsA"^{\G_{\mA}}};
{\ar@{~~>}_{C^{\s(f)}}"CsB";"CsA"};
\endxy
\ee
on the bottom by using the inverse of the Gelfand transform, 
the positive map $f,$ and then 
the Gelfand transform once more. This defines a positive map 
$C^{\s(f)}:C\big(\s(\mB)\big)\stoch C\big(\s(\mA)\big)$
which we know is uniquely characterized by a 
stochastic map $\s^{f}:\s(\mA)\stoch\s(\mB)$ by the RMK theorem.  
We could simply define $\s(f)$ in this manner and move on. 
However, a more geometrically appealing 
formula for $\s(f)$ is available in terms of the convex structure on
the set of states with reference to neither the Gelfand transform
nor the RMK theorem. 
As the reader may have guessed from the suggestive drawing above, 
this merits a review of convex sets, 
their extreme points, the Krein-Milman theorem, and Choquet theory
\cite{Si11}, \cite{Ph01}. 

\bd
\label{defn:convex}
Let $V$ be a topological vector space. A subset $C\subseteq V$ 
is said to be \emph{\uline{convex}} iff 
\be
\l v+(1-\l)w\in C\qquad\forall\;\l\in[0,1]\text{ whenever }v,w\in C. 
\ee
An \emph{\uline{extreme point}} of $C$ is a point $u\in C$ 
such that if $u=\l v+(1-\l)w$ for some $v,w\in C$ and $\l\in[0,1],$
then $v=u$ and/or $w=u.$ 
The set of extreme points of $C$ is denoted by $\mathrm{ex}(C).$ 
\ed

\bx
\label{ex:statesconvexspectrumextreme}
Let $\mA$ be a commutative $C^*$-algebra. 
Then, $\mS(\mA),$ the set of all states, is a compact convex
subset of $\mA^{*}$ in the weak* topology. 
$\s(\mA),$ the spectrum of $\mA,$ is the set of extreme points of $\mS(\mA).$
This follows, for instance, from Theorem 11.33 in Rudin 
\cite{Ru91}.
\ex

\bd
\label{defn:convexhull}
Let $V$ be a topological vector space and let $X\subseteq V.$ 
The \emph{\uline{convex hull}} of $X$ is given by 
\be
\mathrm{ch}(X):=\left\{\sum_{i}^{\text{finite}}\l_{i}v_{i}\;:\;v_{i}\in X,\;\l_{i}\in[0,1],\;\sum_{i}\l_{i}=1\right\}, 
\ee
i.e. the set of all finite convex combinations of elements of $X.$ 
\ed

The following is an application of the Krein-Milman theorem suitable
for our purposes \cite{Si11}. 

\bt
Let $\mA$ be a commutative $C^*$-algebra. Then the set of all states
on $\mA$ is the closure of the convex hull of the set of characters on $\mA,$ i.e.
\be
\mS(\mA)=\overline{\mathrm{ch}\big(\s(\mA)\big)}.
\ee
\et

\bprf
This follows from the usual Krein-Milman theorem, the conclusions
of Example \ref{ex:statesconvexspectrumextreme}, and the fact
that $\mA^{*}$ is a locally convex topological vector space
(it is locally convex due to the separating family of seminorms
generating its topology). 
\eprf

This important result implies the existence of a 
probability measure associated to any state. 

\bt
\label{thm:regularprobmeasonspectrumfromstate}
Let $\mA$ be a commutative $C^*$-algebra and let $\w\in\mS(\mA)$ be a
state on $\mA.$ Then there exists a unique regular probability measure 
$\mu$ on $\s(\mA)$ such that
\be
\label{eq:regularprobmeasonspectrumfromstate}
\int\displaylimits_{\s(\mA)}\mathrm{ev}_{a}\;d\mu=\w(a)\qquad\forall\;a\in\mA. 
\ee
Here $\mathrm{ev}_{a}:\s(\mA)\to\C$ denotes the evaluation at the point $a$
sending a character $\chi\in\s(\mA)$ to $\chi(a).$ 
\et

The statement of the theorem can be obtained from the RMK 
Representation Theorem combined with the Gelfand-Naimark Theorem. 
However, we prefer to prove this theorem directly using Choquet theory,  
which provides a concrete and geometric construction of the measure.
We thank Benjamin Russo for discussions leading to this construction.

\bprf
We first establish existence and then uniqueness. 

By the Krein-Milman Theorem, 
$\w\in\overline{\mathrm{ch}\big(\s(\mA)\big)}.$ 
Hence, there exists a net
$\xi:\Theta\to\mathrm{ch}\big(\s(\mA)\big)$ such that 
$\lim\xi=\w$ in the weak* topology on $\mS(\mA).$ 
Since each $\xi_{\theta}$ is in the 
convex hull of $\s(\mA),$ there exists an $N_{\theta}\in\N$ 
and a convex decomposition of $\xi_{\theta}$ into characters, i.e. 
\be
\xi_{\theta}=\sum_{i_{\theta}=1}^{N_{\theta}}\l_{i_{\theta}}\chi_{i_{\theta}},
\ee
where $\chi_{i_{\theta}}\in\s(\mA)$ for all $i_{\theta}\in\{1,\dots,N_{\theta}\}.$ 
Set 
\be
\mu_{\theta}:=\sum_{i_{\theta}=1}^{N_{\theta}}\l_{i_{\theta}}\de_{\chi_{i_{\theta}}}
\ee
to be the convex combination of Dirac measures on these characters. 
This defines a net of regular probability measures on $\s(\mA).$ 
Since $\mathrm{ProbMeas}\big(\s(\mA)\big)$ is compact in the vague 
topology, 
there exists a subnet $(\mu_{\theta'})_{\theta'\in\Theta'\subseteq\Theta}$ 
of these measures converging in the vague topology
to some regular probability measure, which we denote by $\mu.$ 
This probability measure satisfies (\ref{eq:regularprobmeasonspectrumfromstate})
since
\be
\begin{split}
\int\displaylimits_{\s(\mA)}\mathrm{ev}_{a}\;d\mu&=\lim_{\theta'\in\Theta'}\int\displaylimits_{\s(\mA)}\mathrm{ev}_{a}\;d\mu_{\theta'}\quad\text{ by continuity of $\mathrm{ev}_{a}$ and vague convergence}\\
&=\lim_{\theta'\in\Theta'}\sum_{i_{\theta'}=1}^{N_{\theta'}}\l_{i_{\theta'}}\int\displaylimits_{\s(\mA)}\mathrm{ev}_{a}\;d\de_{\chi_{i_{\theta'}}}\quad\text{ by linearity of the measure and integral}\\
&=\lim_{\theta'\in\Theta'}\sum_{i_{\theta'}=1}^{N_{\theta'}}\l_{i_{\theta'}}\chi_{i_{\theta'}}(a)\quad\text{ by definition of the Dirac measure}\\
&=\lim_{\theta'\in\Theta'}\sum_{i_{\theta'}=1}^{N_{\theta'}}\xi_{\theta'}(a)\quad\text{ by definition of $\xi$}\\
&=\w(a)\quad\text{ by weak* convergence.}
\end{split}
\ee
This proves existence of such a measure. 

Uniqueness is more subtle and follows, for instance, 
from the Choquet-Meyer Theorem. More specifically, it follows from
Corollary 10.9 in \cite{Ph01}, whose notation we follow. 
We do not introduce all the necessary definitions to explain this in full detail,
but we outline the idea using the notation and terminology from \cite{Ph01}.
$\mA^{*},$ the topological dual space of $\mA$ with the weak* topology, is a 
locally convex vector space (this follows, for instance, from 
Proposition \ref{prop:weakstartopologybase} above 
and Theorem 1.37 in \cite{Ru91}). 
The set of all states $\mS(\mA)$ on a commutative (unital)
$C^*$-algebra is a compact convex subset of $\mA^{*}$ and is
a Choquet simplex (see Theorem 3.1.18 in \cite{Sa71} and note that states on a 
commutative $C^*$-algebra are automatically tracial).
Because $\mathrm{ex}\big(\mS(\mA)\big)=\s(\mA)$ is a closed set
(since it is compact), 
every maximal measure is supported on
$\s(\mA)$ (see the comments preceding the Bishop-de Leeuw Theorem
in Section 4 of \cite{Ph01}). To see that  
\be
\int\displaylimits_{\mS(\mA)}\vf\;d\mu=\int\displaylimits_{\mS(\mA)}\vf\;d\de_{\w}
\ee
($\mu\sim\de_{\w}$ in the notation of \cite{Ph01})
for all continuous real-valued affine functions $\vf$ on $\mS(\mA),$
notice that all functions of the form $\mathrm{ev}_{a}:\mA^{*}\to\C$ 
are continuous and linear and hence remain continuous when restricted to 
$\mS(\mA)\subseteq\mA^{*}.$ By construction of $\mu$ from the subnet
$(\mu_{\theta'})$ and the definition of the Dirac measure,
\be
\int\displaylimits_{\mS(\mA)}\mathrm{ev}_{a}\;d\mu=
\int\displaylimits_{\s(\mA)}\mathrm{ev}_{a}\;d\mu=
\w(a)=\int\displaylimits_{\mS(\mA)}\mathrm{ev}_{a}\;d\de_{\w}\qquad\forall\;a\in\mA.
\ee
Now, suppose that $\vf$ is an arbitrary continuous real-valued affine function
on $\mS(\mA).$ Then, $\vf$ restricted to $\s(\mA)$ is a continuous
function and therefore equals $\mathrm{ev}_{a}$ for some
$a\in\mA$ by Proposition \ref{prop:Gelfandfacts}.
If we show that $\vf=\mathrm{ev}_{a}$ on all of $\mS(\mA),$
then we will be done by using the previous calculation. 
Note that because $\mS(\mA)$ is a compact convex subset 
of a locally convex space, 
the supremum and infimum of $\vf$ are achieved on a face of $\mS(\mA).$
This is a consequence, for instance, of Proposition 8.26 in \cite{Si11} by noting
that the image of a non-empty compact convex set under such a continuous 
map is an interval in $\R$ and the supremum and infimum of an interval 
are at most two points. These endpoints are faces of the interval
and so their pre-images are faces of the simplex $\mS(\mA).$ 
Furthermore, the extreme points of a face are 
extreme points of the original convex subset 
(a simple calculation from the definitions shows this). Hence, $|\vf|$ achieves its
supremum on $\s(\mA)$ (and possibly on a larger domain as well). 
Hence, by the Hahn-Banach Theorem
(see Section 5.22 of \cite{Ru87}), there exists a unique extension 
of $\vf|_{\s(\mA)}=\mathrm{ev}_{a}$ to $\mS(\mA).$ But we already know that
$\mathrm{ev}_{a}$ is such an extension. Hence $\vf=\mathrm{ev}_{a}.$
Since the result was already proved for such functions, 
$\mu\sim\de_{\w}$ so uniqueness of the regular probability measure 
follows from the Choquet-Meyer theorem.
\eprf

\br
Uniqueness also follows from a theorem due to Bauer
since $\mS(\mA)$ is a Bauer simplex for $\mA$ a commutative
$C^*$-algebra (see Theorem II.4.1 in \cite{Al71}). 
\er

\br
\label{rmk:lackofcontinuity}
Note that because the evaluation map 
$\mathrm{ev}_{E}:\mathrm{ProbMeas}\big(\s(\mA)\big)\to\R$ is not necessarily
continuous for every Borel subset $E\subseteq\s(\mA),$ in general
\be
\mu(E)\ne\lim_{\q\in\Theta}\sum_{i_{\theta}\in I_{\q}}\l_{i_{\theta}}(E),
\ee
where $I_{\q}$ is the subset of $\{1,\dots,N_{\q}\}$ for which 
$\chi_{i_\q}\in E.$ 
A counterexample is given by the sequence 
\be
\N\ni n\mapsto\mu_{n}:=\frac{1}{n}\sum_{k=1}^{n}\de_{\frac{k}{n}}
\ee
of a uniform sum of Dirac measures on $[0,1].$ This sequence 
converges to the standard Lebesgue measure $\mu$ 
on $[0,1]$ in the vague topology.
However, 
$
\mu_{n}\big(\Q\cap[0,1]\big)=1
$
while
$\mu\big(\Q\cap[0,1]\big)=0.$ 

In the finite-dimensional setting, this does not occur, 
nets are not required, and $\w$ \emph{is}
uniquely expressed as a convex decomposition of extremal states. 
The measure of $E$ is the sum of the weights corresponding to
restricting this sum to the contributions coming only from $E.$ 
This is because the definition of a Choquet simplex is one that agrees
with the definition of a simplex in finite dimensional vector spaces
(see Proposition 10.10 in \cite{Ph01}). 
\er

\bt
\label{thm:representingpullbackstatesasmeasures}
Let $\mB$ and $\mA$ be commutative $C^*$-algebras and let 
$f:\mB\stoch\mA$ be a positive (unital) map. For each character
$\chi\in\s(\mA),$ there exists a unique regular probability measure
$\s^{f}_{\chi}\in\mathrm{ProbMeas}\big(\s(\mB)\big)$ 
such that 
\be
\label{eq:characterpullsbacktostate}
\int\displaylimits_{\s(\mB)}\mathrm{ev}_{b}\;d\s^{f}_{\chi}=\chi\big(f(b)\big)
\qquad\forall\;b\in\mB. 
\ee
Furthermore, the assignment 
\be
\begin{split}
\s(\mA)&\xrightarrow{\s(f)\equiv\s^{f}}\mathrm{ProbMeas}\big(\s(\mB)\big)\\
\chi&\xmapsto{\qquad\quad}\s^{f}_{\chi},
\end{split}
\ee
where $\s^{f}_{\chi}$ is the unique regular probability measure on $\s(\mB)$ 
associated to the state $\chi\circ f\in\mS(\mB)$ via Theorem  
\ref{thm:regularprobmeasonspectrumfromstate}, 
defines a stochastic map $\s(\mA)\stoch\s(\mB).$
\et

\bprf
Since $\chi\circ f$ is a state on $\mA,$ Theorem 
\ref{thm:regularprobmeasonspectrumfromstate} 
guarantees there exists a unique regular probability measure $\s^{f}_{\chi}$
satisfying (\ref{eq:characterpullsbacktostate}).
The only thing left to check is that this defines a stochastic map,
i.e. that the function $\s^{f}:\s(\mA)\to\mathrm{ProbMeas}\big(\s(\mB)\big)$
is continuous with respect to the vague topology. 
Let $\chi:\Theta\to\s(\mA)$ be a net of characters converging to $\lim\chi.$ 
The goal is to show that 
\be
\lim_{\q\in\Theta}\int\displaylimits_{\s(\mB)}\vf\;d\s^{f}_{\chi_{\q}}=\int\displaylimits_{\s(\mB)}\vf\;d\s^{f}_{\lim\chi}
\qquad\forall\;\vf\in C\big(\s(\mB)\big). 
\ee
The Gelfand transform $\G_{\mB}:\mB\to C\big(\s(\mB)\big)$ sending
$b\in\mB$ to the evaluation function $\mathrm{ev}_{b}$ is a
$C^*$-algebra isomorphism by Proposition \ref{prop:Gelfandfacts}.
Hence, every continuous function on the spectrum of $\mB$ is 
of this form. This fact simplifies the above goal to simply showing that 
\be
\lim_{\q\in\Theta}\int\displaylimits_{\s(\mB)}\mathrm{ev}_{b}\;d\s^{f}_{\chi_{\q}}
=\int\displaylimits_{\s(\mB)}\mathrm{ev}_{b}\;d\s^{f}_{\lim\chi}
\qquad\forall\;b\in\mB. 
\ee
This follows from the previous results since
\be
\begin{split}
\lim_{\q\in\Theta}\int\displaylimits_{\s(\mB)}\mathrm{ev}_{b}\;d\s^{f}_{\chi_{\q}}&=
\lim_{\q\in\Theta}\chi_{\q}\big(f(b)\big)\quad\text{by (\ref{eq:characterpullsbacktostate})}\\
&=(\lim\chi)\big(f(b)\big)\quad\text{ by weak* convergence of the net $\chi$}\\
&=\int\displaylimits_{\s(\mB)}\mathrm{ev}_{b}\;d\s^{f}_{\lim\chi}\quad\text{by (\ref{eq:characterpullsbacktostate})}
\end{split}
\ee
for all $b\in\mB.$ 
Therefore $\s^{f}$ is a stochastic map. 
\eprf

\bt
\label{thm:stochasticspectrum}
The assignment 
\be
\begin{split}
\cCAlgPos^{\op}&\xrightarrow{\s}\cHStoch\\
\mA&\mapsto\s(\mA)\\
\Big(\mB\;\xy0;/r.25pc/:\ar@{~>}(-3,0);(3,0)^{f}\endxy\;\mA\Big)&\mapsto\Big(\s(\mA)\;\xy0;/r.25pc/:\ar@{~>}(-6,0);(6,0)^{\s(f)\equiv\s^{f}}\endxy\;\s(\mB)\Big)\\
\end{split}
\ee
from Definition \ref{def:spectrum}
and Theorem \ref{thm:representingpullbackstatesasmeasures}
defines a functor, henceforth referred to as the 
\emph{\uline{stochastic spectrum functor}}. 
\et

\bprf
The identity positive map gets sent to the identity stochastic map
by uniqueness of the measure and the fact that the Dirac measure
is a regular measure satisfying the required conditions. 
Let $\mC,\mB,$ and $\mA$ be commutative $C^*$-algebras and let 
$f:\mB\stoch\mA$ and $g:\mC\stoch\mB$ be positive maps. 
Let $\s^{f}:\s(\mA)\to\mathrm{ProbMeas}\big(\s(\mB)\big),$
$\s^{g}:\s(\mB)\to\mathrm{ProbMeas}\big(\s(\mC)\big),$
and 
$\s^{f\circ g}:\s(\mA)\to\mathrm{ProbMeas}\big(\s(\mC)\big)$
denote the images of $f,g,$ and $f\circ g,$ under $\s,$ respectively. 
Let $\s^{g}\circ\s^{f}:\s(\mA)\to\mathrm{ProbMeas}\big(\s(\mC)\big)$ 
denote the composition $\s(\mA)\stoch\s(\mB)\stoch\s(\mC)$ of the stochastic
map $\s^{f}$ followed by $\s^{g}.$ The goal is to show 
$\s^{f\circ g}=\s^{g}\circ\s^{f}.$ Since $\s^{f\circ g}_{\chi}$
and $(\s^{g}\circ\s^{f})_{\chi}$ are regular probability measures for each
$\chi\in\s(\mA),$
by uniqueness in Theorem \ref{thm:representingpullbackstatesasmeasures}, 
it suffices to show 
\be
\int\displaylimits_{\s(\mC)}\mathrm{ev}_{c}\;d(\s^{g}\circ\s^{f})_{\chi}=
\chi\Big(f\big(g(c)\big)\Big)
\ee
for all $c\in\mC$ and for all $\chi\in\s(\mA).$ 
Therefore, first let $c\in\mC$ be positive so that $\mathrm{ev}_{c}\ge0,$
and fix $\chi\in\s(\mA),$ and let 
\be
0\le s_1\le s_2\le \cdots\le\mathrm{ev}_{c}
\ee
be a sequence $s:\N\to\R^{\s(\mC)}$ of simple functions converging to 
$\mathrm{ev}_{c}$ and satisfying all the conditions as in 
Lemma \ref{lem:contcompactsequenceapprox}.
Then, 
\be
\begin{split}
\chi\Big(f\big(g(c)\big)\Big)
&=\int\displaylimits_{\s(\mB)}\mathrm{ev}_{g(c)}(\psi)\;d\s^{f}_{\chi}(\psi)\quad\text{by definition of $\s^{f}_{\chi}$}\\
&=\int\displaylimits_{\s(\mB)}\psi\big(g(c)\big)\;d\s^{f}_{\chi}(\psi)\quad\text{by definition of $\mathrm{ev}_{g(c)}$}\\
&=\int\displaylimits_{\s(\mB)}\left(\;\int\displaylimits_{\s(\mC)}\mathrm{ev}_{c}\;d\s^{g}_{\psi}\right)d\s^{f}_{\chi}(\psi)\quad\text{by definition of $\s^{g}_{\psi}$}\\
&=\int\displaylimits_{\s(\mB)}\lim_{n\to\infty}\sum_{i_{n}}s_{n,i_n}\s^{g}_{\psi}(E_{n,i_n})\;d\s^{f}_{\chi}(\psi)\quad\text{by definition of Lebesgue integral}\\
&=\lim_{n\to\infty}\sum_{i_{n}}s_{n,i_n}\int\displaylimits_{\s(\mB)}\s^{g}_{\psi}(E_{n,i_n})\;d\s^{f}_{\chi}(\psi)\quad\text{by Monotone Convergence Theorem}\\
&=\lim_{n\to\infty}\sum_{i_{n}}s_{n,i_n}(\s^{g}\circ\s^{f})_{\chi}(E_{n,i_n})\quad\text{by definition of $(\s^{g}\circ\s^{f})_{\chi}$}\\
&=\int\displaylimits_{\s(\mC)}\lim_{n\to\infty}\sum_{i_{n}}s_{n,i_n}\chi_{E_{n,i_n}}\;d(\s^{g}\circ\s^{f})_{\chi}\quad\text{by definition of Lebesgue integral}\\
&=\int\displaylimits_{\s(\mC)}\mathrm{ev}_{c}\;d(\s^{g}\circ\s^{f})_{\chi}\quad\text{by definition of $s$.}
\end{split}
\ee
Again, the Monotone Convergence Theorem applies by similar
arguments as in step iii. in the proof of Proposition \ref{prop:compositionstochastic}. 
The case of more general $c\in\mC$ for which $\mathrm{ev}_{c}$ need
not be positive is handled by splitting up $\mathrm{ev}_{c}$  
using a Jordan decomposition and applying a similar procedure. 
\eprf

\begin{figure}
\centering
\begin{tikzpicture}[scale=2.65]
\def\n{1}
\node at (0.45,0.65) {$n=1$};
\draw[-] (-1,0)--(1,0);
\draw[-] (0,0)--(0,1.20);
\draw[-] (-0.05,1) -- (0.05,1);
	\draw[blue,dotted] ({-1+(2-1)/\n},0) -- ({-1+(2-1)/\n},{(erf(10*(-1+(2)/\n))+erf(10*(2-2+\n)/\n))/(2*erf(10))});
	\node[blue] at ({-1+(2-1)/\n},{(erf(10*(-1+(2)/\n))+erf(10*(2-2+\n)/\n))/(2*erf(10))}) {\small $\bullet$};
\draw[red,thick] plot[domain=-1:-0.125,samples=100] ({\x},{(exp(-100*pow(\x,2)))/(0.177245)});
\draw[red,thick] plot[domain=0.125:1,samples=100] ({\x},{(exp(-100*pow(\x,2)))/(0.177245)});
\node[blue] at (0,0) {\phantom{x}};
\end{tikzpicture}
\quad
\begin{tikzpicture}[scale=2.65]
\def\n{10}
\node at (0.45,0.65) {$n=10$};
\draw[-] (-1,0)--(1,0);
\draw[-] (0,0)--(0,1.20);
\draw[-] (-0.05,1) -- (0.05,1);
\foreach \k in {1,2,...,\n} {
	\draw[blue,dotted] ({-1+(2*\k-1)/\n},0) -- ({-1+(2*\k-1)/\n},{(erf(10*(-1+(2*\k)/\n))+erf(10*(2-2*\k+\n)/\n))/(2*erf(10))});
	\node[blue] at ({-1+(2*\k-1)/\n},{(erf(10*(-1+(2*\k)/\n))+erf(10*(2-2*\k+\n)/\n))/(2*erf(10))}) {\small $\bullet$};
	}
\draw[red,thick] plot[domain=-1:-0.125,samples=100] ({\x},{(exp(-100*pow(\x,2)))/(0.177245)});
\draw[red,thick] plot[domain=0.125:1,samples=100] ({\x},{(exp(-100*pow(\x,2)))/(0.177245)});
\end{tikzpicture}
\quad
\begin{tikzpicture}[scale=2.65]
\def\n{31}
\node at (0.45,0.65) {$n=31$};
\draw[-] (-1,0)--(1,0);
\draw[-] (0,0)--(0,1.20);
\draw[-] (-0.05,1) -- (0.05,1);
\foreach \k in {1,2,...,\n} {
	\draw[blue,dotted] ({-1+(2*\k-1)/\n},0) -- ({-1+(2*\k-1)/\n},{(erf(10*(-1+(2*\k)/\n))+erf(10*(2-2*\k+\n)/\n))/(2*erf(10))});
	\node[blue] at ({-1+(2*\k-1)/\n},{(erf(10*(-1+(2*\k)/\n))+erf(10*(2-2*\k+\n)/\n))/(2*erf(10))}) {\footnotesize $\bullet$};
	}
\draw[red,thick] plot[domain=-1:-0.125,samples=100] ({\x},{(exp(-100*pow(\x,2)))/(0.177245)});
\draw[red,thick] plot[domain=0.125:1,samples=100] ({\x},{(exp(-100*pow(\x,2)))/(0.177245)});
\end{tikzpicture}
\begin{tikzpicture}[scale=2.65]
\def\n{50}
\node at (0.45,0.65) {$n=50$};
\draw[-] (-1,0)--(1,0);
\draw[-] (0,0)--(0,1.20);
\draw[-] (-0.05,1) -- (0.05,1);
\foreach \k in {1,2,...,\n} {
	\draw[blue,dotted] ({-1+(2*\k-1)/\n},0) -- ({-1+(2*\k-1)/\n},{(erf(10*(-1+(2*\k)/\n))+erf(10*(2-2*\k+\n)/\n))/(2*erf(10))});
	\node[blue] at ({-1+(2*\k-1)/\n},{(erf(10*(-1+(2*\k)/\n))+erf(10*(2-2*\k+\n)/\n))/(2*erf(10))}) {\tiny $\bullet$};
	}
\draw[red] plot[domain=-1:-0.125,samples=100] ({\x},{(exp(-100*pow(\x,2)))/(0.177245)});
\draw[red] plot[domain=0.125:1,samples=100] ({\x},{(exp(-100*pow(\x,2)))/(0.177245)});
\end{tikzpicture}
\begin{tikzpicture}[scale=2.65]
\def\n{71}
\node at (0.45,0.65) {$n=71$};
\draw[-] (-1,0)--(1,0);
\draw[-] (0,0)--(0,1.20);
\draw[-] (-0.05,1) -- (0.05,1);
\foreach \k in {1,2,...,\n} {
	\draw[blue,dotted] ({-1+(2*\k-1)/\n},0) -- ({-1+(2*\k-1)/\n},{(erf(10*(-1+(2*\k)/\n))+erf(10*(2-2*\k+\n)/\n))/(2*erf(10))});
	\node[blue] at ({-1+(2*\k-1)/\n},{(erf(10*(-1+(2*\k)/\n))+erf(10*(2-2*\k+\n)/\n))/(2*erf(10))}) {\tiny $\bullet$};
	}
\draw[red,thick] plot[domain=-1:-0.125,samples=100] ({\x},{(exp(-100*pow(\x,2)))/(0.177245)});
\draw[red,thick] plot[domain=0.125:1,samples=100] ({\x},{(exp(-100*pow(\x,2)))/(0.177245)});
\end{tikzpicture}
\begin{tikzpicture}[scale=2.65]
\def\n{100}
\node at (0.45,0.65) {$n=100$};
\draw[-] (-1,0)--(1,0);
\draw[-] (0,0)--(0,1.20);
\draw[-] (-0.05,1) -- (0.05,1);
\foreach \k in {1,2,...,\n} {
	\draw[blue,dotted] ({-1+(2*\k-1)/\n},0) -- ({-1+(2*\k-1)/\n},{(erf(10*(-1+(2*\k)/\n))+erf(10*(2-2*\k+\n)/\n))/(2*erf(10))});
	\node[blue] at ({-1+(2*\k-1)/\n},{(erf(10*(-1+(2*\k)/\n))+erf(10*(2-2*\k+\n)/\n))/(2*erf(10))}) {\tiny $\bullet$};
	}
\draw[red,thick] plot[domain=-1:-0.125,samples=100] ({\x},{(exp(-100*pow(\x,2)))/(0.177245)});
\draw[red,thick] plot[domain=0.125:1,samples=100] ({\x},{(exp(-100*pow(\x,2)))/(0.177245)});
\end{tikzpicture}
\caption{The curve drawn in red is the graph of $\rho$ in (\ref{eq:rhoheat})
(cut off from above to fit on the page)
while the (larger) blue bullets represent the heights of the Dirac point distributions. 
The heights for the point measures decrease because the area over which
they are used to approximate functions decreases as $n$ increases.}
\label{fig:heatconvex}
\end{figure}
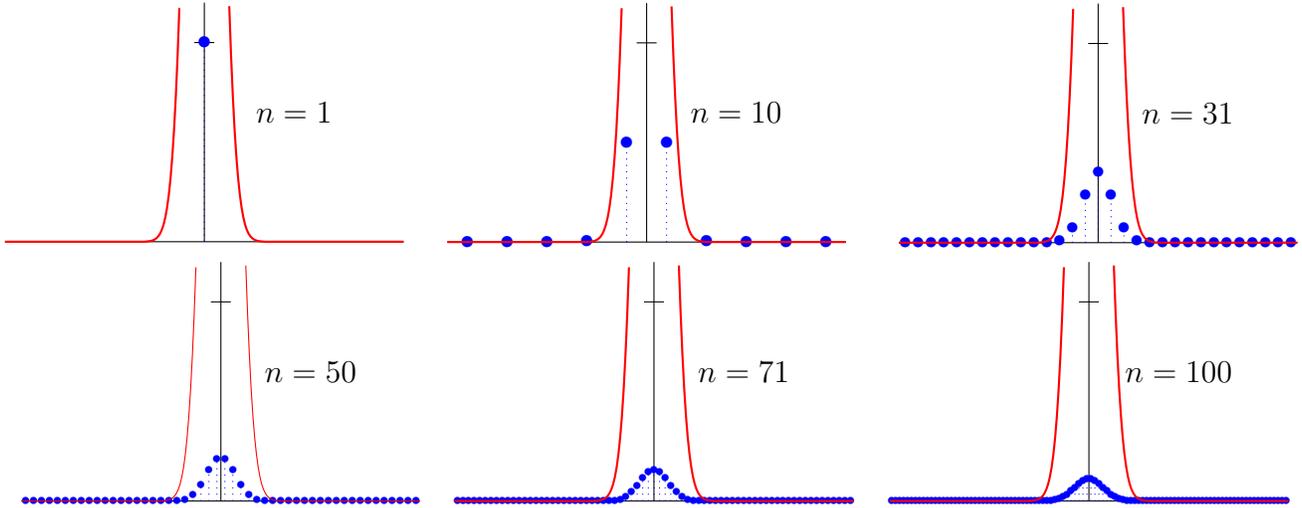

\bx
Let $\rho:[-1,1]\to\R$ be the restriction of the Gaussian used to describe
the distribution of heat after some time with initial condition a source
of heat at the origin, namely
\be
\label{eq:rhoheat}
[-1,1]\ni x\mapsto\rho(x):=\frac{\exp\left(-100x^2\right)}{\int_{-1}^{1}\exp\left(-100y^2\right)dy},
\ee
which has been properly normalized so that the integral of $\rho$ with respect
to the Lebesgue measure is $1.$ 
Consider the sequence of probability measures given by 
\be
\N\ni n\mapsto\mu_{n}:=\sum_{k=1}^{n}\left(\int\limits_{-1+\frac{2(k-1)}{n}}^{-1+\frac{2k}{n}}\rho(x)\;dx\right)\de_{-1+\frac{2k-1}{n}}.
\ee
The distribution and resulting convex sum approximation is depicted in
Figure \ref{fig:heatconvex} for several values of $n\in\N.$ 
Although the sequence of measures $d\mu_{n}$ does not converge
to $\rho\;\!d\mu_{\text{Lebesgue}}$ in the total variation norm topology, 
the convergence holds in the vague topology, i.e. 
\be
\lim_{n\to\infty}\int_{-1}^{1}\vf\;d\mu_{n}=\int_{-1}^{1}\vf(x)\rho(x)\;dx
\qquad\forall\;\vf\in C\big([-1,1]\big).
\ee
This gives an example of a probability measure on $[-1,1]$ that can be
approximated by a sequence of convex sums of Dirac point measures 
with respect to the vague topology. Again, it is not true that 
$\ds\lim_{n\to\infty}\mu_{n}(E)=\int_{E}\rho(x)\;dx$ for all Borel measurable
sets $E.$ For example, take $E:=[-1,1]\cap(\R\setminus\Q),$ all irrationals
in $[-1,1].$ The term on the left is $0$ while the term on the right is $1.$ 
\ex

\subsection{The stochastic commutative Gelfand-Naimark Theorem}
\label{sec:stochasticGN}

Notice that when $\mA=C(X)$ with $X$ a compact Hausdorff space and 
$\w\in\mS\big(C(X)\big),$ Theorem 
\ref{thm:regularprobmeasonspectrumfromstate}
says that there exists a unique regular probability measure $\mu$ 
on $\s\big(C(X)\big)$ such that 
\be
\int\displaylimits_{\s\big(C(X)\big)}\mathrm{ev}_{\vf}\;d\mu=\w(\vf)\qquad\forall\;\vf\in C(X). 
\ee
Meanwhile, by the RMK Representation Theorem, there exists a unique regular 
probability measure $\nu$ on $X$ such that 
\be
\w(\vf)=\int_{X}\vf\;d\nu\qquad\forall\;\vf\in C(X).
\ee
By the Gelfand-Naimark Theorem, we know that the spaces
$X$ and $\s\big(C(X)\big)$ are homeomorphic and one expects
these two measures to be related. 
In fact, the function $h_{X}:X\to\s\big(C(X)\big)$ 
pushes forward the measure $\nu$ to $\mu$ and the inverse
$h_{X}^{-1}:\s\big(C(X)\big)\to X$ pushes forward the measure $\mu$
to $\nu.$ 
This is explained in more detail in the following theorem.

\bt
\label{thm:stochasticGN}
The stochastic spectrum functor together with the
natural isomorphisms 
\be
\xy0;/r.25pc/:
(-17.5,7.5)*+{\cHStoch}="3";
(17.5,7.5)*+{\cHStoch}="1";
(0,-7.5)*+{\cCAlgPos^{\op}}="2";
{\ar"1";"2"^{C}};
{\ar"2";"3"^{\s}};
{\ar"1";"3"_{\id}};
{\ar@{=>}(0,6.5);"2"^(0.37){h}};
\endxy
\aand
\xy0;/r.25pc/:
(-17.5,7.5)*+{\cCAlgPos}="3";
(17.5,7.5)*+{\cCAlgPos}="1";
(0,-7.5)*+{\cHStoch^{\op}}="2";
{\ar"1";"2"^{\s}};
{\ar"2";"3"^{C}};
{\ar"1";"3"_{\id}};
{\ar@{=>}(0,6.5);"2"^(0.37){\G}};
\endxy
\ee
from the commutative Gelfand-Naimark Theorem
form an adjoint equivalence of categories. 
\et

\bprf
The only thing that has not been shown is that 
the transformations $\G$ and $h$ still satisfy naturality.
Note that each is an isomorphism in its respective category
on objects due to the usual Gelfand-Naimark Theorem.  
\begin{enumerate}[i.]
\item
\emph{Proof that $h$ is natural.}
Let $\mu:X\stoch Y$ be a stochastic map and let 
$\s^{C(\mu)}\!:\!\s\big(C(X)\big)\!\stoch\!\s\big(C(Y)\big)$
denote the stochastic map obtained from applying the functor $C$
and then $\s$ to $\mu.$  
The goal is to prove that the diagram 
\be
\xy0;/r.25pc/:
(-12.5,7.5)*+{X}="X";
(-12.5,-7.5)*+{Y}="Y";
(12.5,7.5)*+{\s(C(X))}="sCX";
(12.5,-7.5)*+{\s(C(Y))}="sCY";
{\ar"X";"sCX"^(0.4){h_{X}}};
{\ar"Y";"sCY"_(0.4){h_{Y}}};
{\ar@{~>}"X";"Y"_{\mu}};
{\ar@{~>}"sCX";"sCY"^{\s^{C(\mu)}}};
\endxy
\ee
commutes. Note that for any Borel set $E\in\mathcal{M}_{\s(C(Y))}$ 
and for any $x\in X,$ 
the lower left composition $\de_{h_{Y}}\circ\mu$ applied to
$x$ evaluated at $E$ gives
\be
(\de_{h_{Y}}\circ \mu)_{x}(E)=\int_{Y}\de_{h_{Y}(y)}(E)\;d\mu_{x}(y)
=\int_{Y}\chi_{h_{Y}^{-1}(E)}(y)\;d\mu_{x}(y)
=\mu_{x}\big(h_{Y}^{-1}(E)\big).
\ee
In other words, $(\de_{h_{Y}}\circ \mu)_{x}$ is the pushforward of the measure
$\mu_{x}$ along the function $h_{Y}.$ Now we check that this equals
$\big(\s^{C(\mu)}\circ h_{X}\big)_{x}$ for every $x\in X$ by integrating
$\mathrm{ev}_{E}$ with respect to this pushforward measure and 
$\vf\in C(Y).$ 
By the change of variables 
formula for integrating measurable functions with respect to 
pushforward measures (see Theorem 3.6.1 in \cite{Bo07i}),
\be
\label{eq:stochasticgelfandnatural}
\int\displaylimits_{\s(C(Y))}\mathrm{ev}_{\vf}\;d(\de_{h_{Y}}\circ \mu)_{x}
=\int\displaylimits_{\s(C(Y))}\mathrm{ev}_{\vf}\;d\big(\mu_{x}\circ h_{Y}^{-1}\big)
=\int_{Y}(\mathrm{ev}_{\vf}\circ h_{Y})\;d\mu_{x}
\ee
since $\mathrm{ev}_{\vf}$ is continuous (and hence measurable). 
Therefore, 
\be
\begin{split}
\int\displaylimits_{\s(C(Y))}\mathrm{ev}_{\vf}\;d(\de_{h_{Y}}\circ \mu)_{x}
&=\int_{Y}\mathrm{ev}_{\vf}\big(h_{Y}(y)\big)\;d\mu_{x}(y)\quad\text{ by (\ref{eq:stochasticgelfandnatural})}\\
&=\int_{Y}h_{Y}(y)\big(\vf\big)\;d\mu_{x}(y)\quad\text{ by definition of $\mathrm{ev}_{\vf}$}\\
&=\int_{Y}\vf(y)\;d\mu_{x}(y)\quad\text{ by definition of $h_{Y}$ (\ref{eq:evaluation})}\\
&=\big(C^{\mu}(\vf)\big)(x)\quad\text{ by definition of $C^{\mu}$ (\ref{eq:curlyEonstochasticmaps})}\\
&=h_{X}(x)\big(C^{\mu}(\vf)\big)\quad\text{ by definition of $h_{X}$ (\ref{eq:evaluation})}\\
&=\int\displaylimits_{\s(C(Y))}\mathrm{ev}_{\vf}\;d\s^{C(\mu)}_{h_{X}(x)}
\quad\text{ by definition of $\s^{C(\mu)}_{h_{X}(x)}$.}
\end{split}
\ee
Since this is true for all $\vf\in C(Y),$ the uniqueness of measures on $\s(C(Y))$ 
from Theorem \ref{thm:regularprobmeasonspectrumfromstate} implies
$\s^{C(\mu)}_{h_{X}(x)}=(\de_{h_{Y}}\circ \mu)_{x}$ for all $x\in X.$ 
Therefore, the diagram above commutes 
and $h:\id_{\cHStoch}\Rightarrow\s\circ C$ is a natural isomorphism. 
\item
\emph{Proof that $\G$ is natural.}
Let $f:\mB\stoch\mA$ be a positive map. The goal is to show that the 
diagram 
\be
\xy0;/r.25pc/:
(-12.5,7.5)*+{\mB}="B";
(-12.5,-7.5)*+{\mA}="A";
(12.5,7.5)*+{C(\s(\mB))}="ScB";
(12.5,-7.5)*+{C(\s(\mA))}="ScA";
{\ar"B";"ScB"^(0.4){\G_{\mB}}};
{\ar"A";"ScA"_(0.4){\G_{\mA}}};
{\ar@{~>}"B";"A"_{f}};
{\ar@{~>}"ScB";"ScA"^{C^{\s(f)}}};
\endxy
\ee
commutes, i.e. for arbitrary $b\in\mB,$ the two functions 
$C^{\s(f)}\big(\G_{B}(b)\big)$ and 
$\G_{\mA}\big(f(a)\big)$ are equal on $\s(\mA).$ 
To check equality of these functions, evaluating them on an arbitrary
$\chi\in\s(\mA)$ gives
\be
\Big(C^{\s(f)}\big(\G_{B}(b)\big)\Big)(\chi)
=\int\displaylimits_{\s(\mB)}\mathrm{ev}_{b}\;d\s^{f}_{\chi}
=\chi\big(f(b)\big)
=\Big(\G_{\mA}\big(f(a)\big)\Big)(\chi),
\ee
which is the desired result. 
\end{enumerate}
\eprf

We should also check that this equivalence agrees with the
usual Gelfand-Naimark equivalence for deterministic processes,
i.e. for {$*$-homomorphisms} on the $C^*$-algebra side and
measure-preserving functions on the analytic side. 

\bt
Let $\mB$ and $\mA$ be commutative $C^*$-algebras, let 
$f:\mB\stoch\mA$ be a {$*$-homomorphism}, and let $\chi\in\s(\mA)$
be a character on $\mA.$ 
Then the corresponding regular probability measure $\s^{f}_{\chi}$ 
is the Dirac measure at $\chi\circ f,$ i.e. $\s^{f}_{\chi}=\de_{\chi\circ f}.$ 
Therefore, $\s(f)$ uniquely determines
a continuous function $\s(\mA)\to\s(\mB)$ 
agreeing with the usual spectrum functor on $\cCAlg,$ i.e. the diagram 
\be
\xy0;/r.25pc/:
(-17.5,7.5)*+{\cH}="1";
(-17.5,-7.5)*+{\cHStoch}="2";
(17.5,7.5)*+{\cCAlg^{\op}}="3";
(17.5,-7.5)*+{\cCAlgPos^{\op}}="4";
{\ar"3";"1"_{\s}};
{\ar"4";"2"^{\s}};
{\ar"1";"2"_{\de}};
{\ar"3";"4"};
\endxy
\ee
commutes. 
\et

\bprf
This follows from the fact that $\chi\circ f$ is a character on $\s(\mB)$ 
and that it immediately satisfies the conditions of Theorem 
\ref{thm:representingpullbackstatesasmeasures}
so that the measure is uniquely specified. 
\eprf

\br
The full power of Theorem 
\ref{thm:representingpullbackstatesasmeasures}
is not needed to prove this last result. It follows from a simpler result
due to Bauer (see Proposition 1.4 in \cite{Ph01}).
\er

\section{Closing discussion}
\label{sec:closing}
We briefly discuss the broad physical scope of our results in the context
of quantum mechanics and then the relationship to similar 
conclusions by other authors. 

\subsection{Categories of classical and quantum dynamics}
\label{sec:theclassicalcube}

In this paper, we have collected functors 
\be
\xy0;/r.35pc/:
(-15,20)*+{\cH^{\op}}="cHaus";
(-15,-5)*+{\cHStoch^{\op}}="cHausStoch";
(15,20)*+{\cCAlg}="cCAlg";
(15,-5)*+{\cCAlgPos}="cCAlgPos";
(-30,5)*+{\FinSetFun^{\op}}="FinSetFun";
(-30,-20)*+{\FinSetStoch^{\op}}="FinSetStoch";
(3,5)*+{\fdcCAlg}="fdcCAlg";
(3,-20)*+{\fdcCAlgPos}="fdcCAlgPos";
{\ar@<1.0ex>"cHaus";"cCAlg"^{C}};
{\ar@<1.0ex>"cCAlg";"cHaus"^{\s}};
{\ar@<1.0ex>|(0.6){\hole}"cHausStoch";"cCAlgPos"^{C}};
{\ar@<1.0ex>|(0.4){\hole}"cCAlgPos";"cHausStoch"^{\s}};
{\ar@<1.0ex>"FinSetFun";"fdcCAlg"};
{\ar@<1.0ex>"fdcCAlg";"FinSetFun"};
{\ar@<1.0ex>"FinSetStoch";"fdcCAlgPos"};
{\ar@<1.0ex>"fdcCAlgPos";"FinSetStoch"};
{\ar@{^{(}->}"FinSetFun";"cHaus"};
{\ar@{^{(}->}"FinSetStoch";"cHausStoch"};
{\ar@{^{(}->}"fdcCAlg";"cCAlg"};
{\ar@{^{(}->}"fdcCAlgPos";"cCAlgPos"};
{\ar"FinSetFun";"FinSetStoch"_{\de}};
{\ar|(0.55){\hole}|(0.65){\hole}"cHaus";"cHausStoch"^(0.4){\de}};
{\ar"fdcCAlg";"fdcCAlgPos"};
{\ar"cCAlg";"cCAlgPos"};
\endxy
\ee
between a variety of categories. The left face of the cube
consists of dynamics on classical systems (phase space). 
The right face describes the algebraic counterparts in terms of observables
of those classical systems (functions on phase space) and how those evolve
under the dynamics. The top face describes 
deterministic dynamics in both space and algebraic settings. 
The bottom face describes stochastic dynamics 
in both space and algebraic settings. 
The front face describes finite classical systems (systems whose phase
space is finite). 
The back face describes classical systems whose phase space is
possibly infinite.
The functors between the left face and the right face are all equivalences.
The equivalences on the top face are due to Gelfand and Naimark
while the back edge of the bottom face has been described in this work. 

This cube has a natural extension to the right for
arbitrary $C^*$-algebras. $*$-homomorphisms are therefore
non-commutative analogues of deterministic processes
while completely positive maps are non-commutative analogues of 
non-deterministic processes. 
\be
\xy0;/r.35pc/:
(-15,20)*+{\cH^{\op}}="cHaus";
(-15,-5)*+{\cHStoch^{\op}}="cHausStoch";
(15,20)*+{\cCAlg}="cCAlg";
(15,-5)*+{\cCAlgPos}="cCAlgPos";
(-30,5)*+{\FinSetFun^{\op}}="FinSetFun";
(-30,-20)*+{\FinSetStoch^{\op}}="FinSetStoch";
(3,5)*+{\fdcCAlg}="fdcCAlg";
(3,-20)*+{\fdcCAlgPos}="fdcCAlgPos";
(45,20)*+{\CAlg}="CAlg";
(45,-5)*+{\CPCAlg}="CPCAlg";
(33,5)*+{\fdCAlg}="fdCAlg";
(33,-20)*+{\CPfdCAlg}="fdCAlgPos";
{\ar@<1.0ex>"cHaus";"cCAlg"^{C}};
{\ar@<1.0ex>"cCAlg";"cHaus"^{\s}};
{\ar@<1.0ex>|(0.6){\hole}"cHausStoch";"cCAlgPos"^{C}};
{\ar@<1.0ex>|(0.4){\hole}"cCAlgPos";"cHausStoch"^{\s}};
{\ar@<1.0ex>"FinSetFun";"fdcCAlg"};
{\ar@<1.0ex>"fdcCAlg";"FinSetFun"};
{\ar@<1.0ex>"FinSetStoch";"fdcCAlgPos"};
{\ar@<1.0ex>"fdcCAlgPos";"FinSetStoch"};
{\ar@{^{(}->}"FinSetFun";"cHaus"};
{\ar@{^{(}->}"FinSetStoch";"cHausStoch"};
{\ar@{^{(}->}"fdcCAlg";"cCAlg"};
{\ar@{^{(}->}"fdcCAlgPos";"cCAlgPos"};
{\ar@{^{(}->}"cCAlg";"CAlg"};
{\ar@{^{(}->}"fdcCAlg";"fdCAlg"};
{\ar@{^{(}->}|(0.6){\hole}"cCAlgPos";"CPCAlg"};
{\ar@{^{(}->}"fdcCAlgPos";"fdCAlgPos"};
{\ar@{^{(}->}"fdCAlg";"CAlg"};
{\ar@{^{(}->}"fdCAlgPos";"CPCAlg"};
{\ar"FinSetFun";"FinSetStoch"};
{\ar|(0.55){\hole}|(0.65){\hole}"cHaus";"cHausStoch"};
{\ar"fdcCAlg";"fdcCAlgPos"};
{\ar|(0.6){\hole}"cCAlg";"cCAlgPos"};
{\ar"fdCAlg";"fdCAlgPos"};
{\ar"CAlg";"CPCAlg"};
\endxy
\ee

The theorems presented in this work offer more motivation
for thinking of completely positive maps and the dynamics
associated with them in quantum theory as a natural extension of 
stochastic dynamics in classical theory.
This can be seen more concretely if one adds more data to these categories. 
The categories above only describe dynamics of the observables of a
system, but not of individual states. In other words, we have so far
described the Heisenberg picture of quantum mechanics. 
However, this is a red herring. Evolution of states is also described
in this setting because a state \emph{is} a completely positive map. 
If a completely positive map $f:\mB\stoch\mA$ describes the evolution of
observables described by $\mA$ to observables described by $\mB,$
and $\w:\mA\stoch\C$ is a state on $\mA,$ then $\w$ evolves to 
$\w\circ f:\mB\stoch\C$ under the evolution described by $f.$ 
Therefore, one can include \emph{both} evolution of observables \emph{and}
states by taking
the slice category of $\CPCAlg$ over $\C,$ denoted by 
$\CPCAlg/_{\C}.$ The objects of this category are 
$C^*$-algebras equipped with states $\w:\mA\stoch\C.$
A morphism from $\xi:\mB\stoch\C$ to $\w:\mA\stoch\C$ consists of a
completely positive map $f:\mB\stoch\mA$ such that $\xi=\w\circ f.$ 
This is a non-commutative generalization of spaces equipped with 
probability measures. 
This situation was described in Example \ref{ex:unitaryevolution}
though more general dynamics is encompassed in Lindblad's formalism \cite{Li76}.

\subsection{Categorical probability theory}
\label{sec:relationshiptoCTPT}

As mentioned earlier, Sections \ref{sec:stochasticmatrices} and 
\ref{sec:GelfandNaimark} were entirely a review though phrased in 
a categorical setting that was used to state and prove various results in 
Section \ref{sec:stochasticGelfandNaimark}. 
Many of the results in Section \ref{sec:stochasticGelfandNaimark} 
relied on well known facts in analysis though our categorical framework
brought some questions to surface which seem to have not been 
addressed earlier. 
Similar results were obtained by
Furber and Jacobs \cite{FuJa13}, who
proved an equivalence
of categories between $\cHStoch^{\op}$ and $\cCAlgPos$ by
realizing the category $\cHStoch$ 
as the Kleisli category of a particular monad, a slight variant of the Giry monad, 
on the category $\cH.$ 
There are three main differences between our findings.
\begin{enumerate}[i.]
\item
Rather than constructing the Giry monad, we worked directly with the
objects and morphisms of $\cHStoch$ bypassing the monad entirely. 
This saves us quite a bit of space from explaining background material 
on monads, Kleisli categories, and how these are related to adjunctions,
which is all beautifully explained in \cite{Se73}.
Furthermore, it makes our presentation more accessible to those with a
minimal categorical background.
However, as a result, some proofs become a bit more
involved compared to \cite{FuJa13}. For instance, given two stochastic maps 
$\mu:X\stoch Y$ and $\nu:Y\stoch Z,$
proving that the composition $\nu\circ\mu$ is well-defined, 
continuous, regular, and a Markov kernel was quite involved
(see Proposition \ref{prop:compositionstochastic} and the preceding lemmas). 
From the Giry monad perspective, $\mathrm{ProbMeas}$ can
be viewed as a functor from the category $\cH$ to itself (one must prove this). 
In addition, there are two natural transformations that provide
$\mathrm{ProbMeas}$ with the structure of a monad. 
The inclusion of a space $X$ into $\mathrm{ProbMeas}(X)$ via
the Dirac measure defines a natural transformation
$\de:\id_{\cH}\Rightarrow\mathrm{ProbMeas}.$
In addition, there is a natural transformation 
$b:\mathrm{ProbMeas}\circ \mathrm{ProbMeas}\Rightarrow\mathrm{ProbMeas}$ 
given on a space $X$ by the assignment
\be
\begin{split}
\mathrm{ProbMeas}\big(\mathrm{ProbMeas}(X)\big)
&\xrightarrow{b_{X}}\mathrm{ProbMeas}(X)\\
\W&\mapsto\left(\mathcal{M}_{X}\ni E\xmapsto{b_{X}(\W)}\hspace{-3mm}\int\displaylimits_{\mathrm{ProbMeas}(X)}\hspace{-5mm}\mathrm{ev}_{E}\;d\W\right).
\end{split}
\ee
Notice that one must still show that $\mathrm{ev}_{E}$ is measurable
to define this map. 
$b_{X}$ is surprisingly easy to interpret and has significant physical meaning.
$\mathrm{ProbMeas}(X)$ is a convex subset of the vector space of all 
finite measures. 
A probability measure on it therefore has a barycenter, 
its center of mass. 
The map $b_{X}$ is precisely the assignment of that probability measure to 
its barycenter. 
This follows from the fact that for any Borel set $E\in\mathcal{M}_{X},$
\be
\int\displaylimits_{\mathrm{ProbMeas}(X)}\hspace{-4mm}\mathrm{ev}_{E}\;d\W
\quad=\int\displaylimits_{\mathrm{ProbMeas}(X)}\hspace{-6mm}\mathrm{ev}_{E}(\mu)\;d\W
(\mu)
\quad=\int\displaylimits_{\mathrm{ProbMeas}(X)}\hspace{-6mm}\mu(E)\;d\W(\mu).
\ee
A physically relevant instance where the barycenter naturally occurs 
will be discussed after point iii. 
Going back to the Giry monad and $\cHStoch$, using the barycenter map, 
the definition of the composition 
$\nu\circ\mu$ can then be taken to be the 
composition of the functions
\be
X\xrightarrow{\mu}\mathrm{ProbMeas}(Y)\xrightarrow{\mathrm{ProbMeas}(\nu)}\mathrm{ProbMeas}\big(\mathrm{ProbMeas}(Z)\big)\xrightarrow{b_{Z}}\mathrm{ProbMeas}(Z)
\ee
from which it immediately follows that the composition is continuous, since
each function in this composition is continuous. 
The relationship between our formula for the composition and this composition
of functions follows from the change of variables formula
since applying the above composition to a point gives
\be
X\ni x\mapsto\mu_{x}\mapsto\nu_{*}\mu_{x}\mapsto b_{Z}(\nu_{*}\mu_{x}),
\ee
where $\nu_{*}(\mu_{x})$ denotes the pushforward measure of $\mu_{x}$
along the (measurable) function $\nu:Y\to\mathrm{ProbMeas}(Z).$ 
This evaluates on measurable sets $E\in\mathcal{M}_{Z}$ to
\be
\big(b_{Z}(\nu_{*}\mu_{x})\big)(E)
=\int\displaylimits_{\mathrm{ProbMeas}(Z)}\hspace{-6mm}\mathrm{ev}_{E}\;d(\nu_{*}\mu_{x})
=\int_{Y}(\mathrm{ev}_{E}\circ\nu)\;d\mu_{x},
\ee
which agrees with our definition of composition
(\ref{eq:compositionofstochasticmaps}). 
By using an intermediate amount of 
category theory, our work offers a bridge
towards more categorical approaches to probability theory \cite{Gi82}. 

\item
In \cite{FuJa13}, Furber and Jacobs use the RMK Representation Theorem immediately 
so that a stochastic map $X\stoch Y$ is defined to be a continuous function 
$X\to \mS\big(C(Y)\big)$ into states on the $C^*$-algebra of continuous
functions on $Y$ instead of a continuous function 
$X\to\mathrm{ProbMeas}(Y)$ into regular probability measures on $Y$
(this is how the Kleisli category associated to their version of the 
Giry monad was defined).
There is no significant difference between these perspectives. However, 
we prefer to think of the left side of our cube in Section 
\ref{sec:theclassicalcube} as purely analytical, topological, and 
measure-theoretic. In particular, we provided an explicit 
definition of the composition using Markov kernels and without
any reference to the RMK theorem. 
The right face of this cube is more algebraic and utilizes algebras of
function spaces, states, positive maps, etc.
Of course, the point of the equivalence between these two faces
is that this distinction is just psychological. Nevertheless, it is 
helpful to distinguish the different tools used in both subjects. 
As a result, many of the proofs in our paper are completely different in nature
from those of Furber and Jacobs.

\item
Finally, Furber and Jacobs showed that 
$C:\cHStoch^{\op}\to\cCAlgPos$ is 
fully faithful and essentially surjective \cite{FuJa13}.
However, 
any abstract nonsense construction of an explicit inverse will necessarily
relate any commutative $C^*$-algebra to one
of continuous functions on some space and will therefore use the 
RMK theorem in its construction. 
We have therefore provided an explicit construction of an inverse
without using the Gelfand transform nor the RMK theorem. 
Furthermore, we showed that the usual Gelfand transform is
part of the adjoint equivalence for the stochastic Gelfand-Naimark Theorem. 
Our construction of the stochastic spectrum functor applies to 
positive maps and take characters to states. Choquet theory 
was used to ``spread'' this state out onto the spectrum in terms
of a probability measure which was obtained as a limit of 
successive approximations via decompositions of the state into convex
combinations of Dirac measures supported on the space of characters. 
This perspective may be useful for 
explorations in the category of non-commutative $C^*$-algebras. 
\end{enumerate}

One instance where the barycenter map is important in physics occurs in
preparing systems on which to perform experiments. 
In an experiment, one has a protocol for preparing a state, which will
then be probed in some way, say, by measuring some observable. 
Because the probe may alter the state, 
one needs the protocol to be as accurate as possible so that the
same state can be constructed so that it can be probed in subsequent experiments. 
Indeed, this is how a state is defined mathematically---via its expectation 
values.
However, in practice, protocols, as well as instruments, are rarely
perfect. As a result, the prepared states may differ slightly from one another
so that they are more accurately represented by a probability measure
on the set of states. Hence, an experimentalist is not always obtaining
expectation values for the same state. Nevertheless,
after viewing the set of states as probability measures via the RMK
Representation Theorem,
the barycenter map provides one with a canonical state that 
represents the average state prepared by the protocol. 
\emph{This} is the state very often referred to in most textbooks, 
but one should keep in mind that the barycenter represents only one part
of all the data. 

\subsection{The Baire approach}
\label{sec:otherapproaches}
One may have circumvented many of the difficulties we encountered by using
Baire sets instead of Borel sets. This perspective is emphasized in works
such as \cite{EFHN15}.
In this case, one has the following analogue of 
Lemma \ref{lem:evaluationofmeasurecontinuous}.

\blem
\label{lem:IddosLemma}
Let $Y$ be a compact Hausdorff space,
let $E\subseteq Y$ be a Baire set, 
and let $\mathrm{Pr}(Y)$ denote all probability measures 
on $Y.$ Then the evaluation function 
\be
\begin{split}
\mathrm{Pr}(Y)&\xrightarrow{\mathrm{ev}_{E}}\R\\
\mu&\mapsto\mu(E)
\end{split}
\ee
is Borel measurable. 
\elem

\bprf
Let $P$ be the collection of all compact $G_{\de}$ subsets of $Y$
and let
\be
L:=
\left\{E\subseteq Y\;:\;\mathrm{Pr}(Y)\xrightarrow{\mathrm{ev}_{E}}[0,1]\text{ is Borel measurable}\right\}
\ee
\begin{enumerate}[i.]
\setlength{\itemsep}{0pt}
\item
The intersection of compact sets is compact and
the finite intersection of a countable collection of open sets is 
still a countable intersection
of open sets. Hence $P$ is a $\pi$-system. 

\item
This is the same as in the proof of 
Lemma \ref{lem:evaluationofmeasurecontinuous} so $L$ is a $\l$-system. 

\item
Let $K=U_1\cap U_2\cap U_3\cap\cdots$ be a compact $G_{\de}$ set in 
$Y$ with $\{U_{n}\}$ an at most countable collection of open sets in $Y.$ 
It suffices to assume that 
$U_{1}\supseteq U_{2}\supseteq U_{3}\supseteq\cdots$
since if this were not the case, the sequence $U_{1}, U_{1}\cap U_{2}, 
U_{1}\cap U_{2}\cap U_{3}, \dots$ would also have intersection equal to $K.$ 
For each $n,$ by Urysohn's Lemma, there exists
a continuous function $f_{n}:Y\to[0,1]$ such that $f_{n}(y)=1$ for $y\in K$ and 
$f_{n}(y)=0$ for $y\in U_{n}^{c}.$ Similarly, it suffices to assume that 
$f_{1}\ge f_{2}\ge f_{3}\ge\cdots$ since if this were not the case, 
the sequence 
$\N\times Y\ni(n,y)\mapsto\min\big\{f_{1}(y),f_{2}(y),\dots,f_{n}(y)\big\}$
would also have the same pointwise limit. From these assumptions, it
follows that 
\be
\lim_{n\to\infty}f_{n}(y)=\chi_{K}(y)\qquad\forall\;y\in Y. 
\ee
Hence, for any $a\in\R,$ the preimage of $(a,\infty)$ under $\mathrm{ev}_{K}$
is 
\be
\begin{split}
\mathrm{ev}_{K}^{-1}\big((a,\infty)\big)&=\big\{\mu\in\mathrm{Pr}(Y)\;:\;\mu(K)>a\big\}\\
&=\left\{\mu\in\mathrm{Pr}(Y)\;:\;\lim_{n\to\infty}\int_{Y}f_{n}\;d\mu>a\right\}\\
&\qquad\text{by the Monotone Convergence Theorem}\\
&=\bigcap_{n=1}^{\infty}\left\{\mu\in\mathrm{Pr}(Y)\;:\;\int_{Y}f_{n}\;d\mu>a\right\}\quad\text{since $\cdots\ge f_{n}\ge f_{n+1}\ge\cdots$}.
\end{split}
\ee
By definition of the vague topology, the set in curly brackets in the last equality
is open. Hence, this is a countable intersection of open sets and is therefore
a Borel subset of $\mathrm{Pr}(Y).$ 
This shows that $\mathrm{ev}_{K}$ is Borel measurable for all 
compact $G_{\de}$ sets $K.$ 
\end{enumerate}
Since all the conditions of Dynkin's $\pi$-$\l$ theorem are 
satisfied, the sigma algebra generated by 
$P,$
which is the Baire $\s$-algebra on $Y,$ is contained $L.$ 
\eprf

This allows one to define the composition
of stochastic maps on Baire sets. One can then use the theorem that
every Baire measure extends uniquely to a regular measure 
under these assumptions (see Corollary 7.3.4. in \cite{Bo07ii}).
However, this would not guarantee that formula 
(\ref{eq:compositionofstochasticmaps}) for the composition of stochastic 
maps is still a valid mathematical expression. The construction
of the regular Borel extension of this measure is somewhat formal
and does not provide one with the explicit formula we have produced.

\pagebreak
\section*{Index of notation}
\label{indexofnotation}
\addcontentsline{toc}{section}{\nameref{indexofnotation}}
\begin{longtable}{c|c|c|c}
\hline
Notation & Name and/or description & Location & Page \\
\hline
$\FinSetFun$&category of finite sets and functions& Notation \ref{not:finprob}&\pageref{not:finprob}\\
\hline
$\FinProb$& \begin{tabular}{c}category of probability measures on finite sets \\and measure-preserving functions\end{tabular}& Notation \ref{not:finprob}&\pageref{not:finprob}\\
\hline
$\mathrm{Pr}(Y)$&probability measures on a finite set $Y$&Definition \ref{defn:prandstochasticmaps}&\pageref{defn:prandstochasticmaps}\\
\hline
$X\stoch Y$&stochastic map from $X$ to $Y$&Definition \ref{defn:prandstochasticmaps}&\pageref{defn:prandstochasticmaps}\\
\hline
$f_{yx}$&\begin{tabular}{c}$yx$-component of a stochastic map\\$f:X\stoch Y$\end{tabular}&Definition \ref{defn:prandstochasticmaps}&\pageref{defn:prandstochasticmaps}\\
\hline
$\de_{yy'}$&\begin{tabular}{c}the value of the Kronecker delta function\\on $(y,y')\in Y\times Y$\end{tabular}&Equation \ref{eq:Kroneckerdelta}&\pageref{eq:Kroneckerdelta}\\
\hline
$\FinSetStoch$&category of finite sets and stochastic maps& Notation \ref{not:finstochandfinsetstoch}&\pageref{not:finstochandfinsetstoch}\\
\hline
$\mA,\mB,\mC$&an algebra (often a unital $C^*$-algebra)& Definition \ref{defn:algebra}&\pageref{defn:algebra}\\
\hline
$1_{\mA}$&multiplicative unit in an algebra $\mA$& Definition \ref{defn:algebra}&\pageref{defn:algebra}\\
\hline
$\lVert\;\cdot\;\rVert$&a norm or seminorm& Definition \ref{defn:normedalgebra}&\pageref{defn:normedalgebra}\\
\hline
iff&``if and only if'' (used only in definitions)& Definition \ref{defn:normedalgebra}&\pageref{defn:normedalgebra}\\
\hline
$\C^{X}$&set of all functions from $X$ to $\C$& Example \ref{ex:CXalgebra}&\pageref{ex:CXalgebra}\\
\hline
$M_{n}(\C)$&set of all complex $n\times n$ matrices& Example \ref{ex:nbybmatrices}&\pageref{ex:nbybmatrices}\\
\hline
$\CAlg$&category of $C^*$-algebras and $*$-homomorphisms& Definition \ref{defn:CAlg}&\pageref{defn:CAlg}\\
\hline
$\cCAlg$&\begin{tabular}{c}the subcategory of $\CAlg$ of commutative\\$C^*$-algebras and $*$-homomorphisms\end{tabular}& Definition \ref{defn:CAlg}&\pageref{defn:CAlg}\\
\hline
$\fdcCAlg$&\begin{tabular}{c}the subcategory of $\cCAlg$ of finite-dimensional\\$C^*$-algebras and $*$-homomorphisms\end{tabular}& Definition \ref{defn:CAlg}&\pageref{defn:CAlg}\\
\hline
CP map&completely positive (unital) map& Definition \ref{defn:CPU}&\pageref{defn:CPU}\\
\hline
$\mB\to\mA$&$*$-homomorphism of $C^*$-algebras& Notation \ref{not:starhomoandCP}&\pageref{not:starhomoandCP}\\
\hline
$\mB\stoch\mA$&completely positive map of $C^*$-algebras& Notation \ref{not:starhomoandCP}&\pageref{not:starhomoandCP}\\
\hline
$\CPCAlg$&\begin{tabular}{c}category of $C^*$-algebras\\ and completely positive (CP) maps\end{tabular}& Notation \ref{defn:CPCAlg}&\pageref{defn:CPCAlg}\\
\hline
$\cCAlgPos$&\begin{tabular}{c}the subcategory of $\CPCAlg$ of\\commutative $C^*$-algebras and CP maps\end{tabular}& Notation \ref{defn:CPCAlg}&\pageref{defn:CPCAlg}\\
\hline
$\fdcCAlgPos$&\begin{tabular}{c}the subcategory of $\cCAlgPos$ of\\finite-dimensional $C^*$-algebras and CP maps\end{tabular}& Notation \ref{defn:CPCAlg}&\pageref{defn:CPCAlg}\\
\hline
$C$&\begin{tabular}{c}functor $C:\FinSetStoch^{\op}\to\fdcCAlgPos$\\and its restriction to $\FinSetFun$\end{tabular}&\begin{tabular}{c}Theorem \ref{thm:FinProbtoCP}\\Theorem \ref{thm:restrictEtofinprob}\end{tabular}&\begin{tabular}{c}\pageref{thm:FinProbtoCP}\\\pageref{thm:restrictEtofinprob}\end{tabular}\\
\hline
$C(X)$&continuous functions on $X$&Example \ref{ex:ConcHaus}&\pageref{ex:ConcHaus}\\
\hline
$C^{f}\equiv C(f)$&pullback $C(Y)\to C(X)$ associated to $f:X\to Y$&Proposition \ref{prop:Conmorphisms}&\pageref{prop:Conmorphisms}\\
\hline
$\cH$&\begin{tabular}{c}category of compact Hausdorff spaces\\ and continuous maps\end{tabular}& Notation \ref{not:cHaus}&\pageref{not:cHaus}\\
\hline
$C$&
functor $\cH^{\op}\to\cCAlg$& Proposition \ref{prop:cHaustocCAlg}&\pageref{prop:cHaustocCAlg}\\
\hline
$\mathfrak{B}$&base on a set (for a topology)& Definition \ref{defn:basefortopology}&\pageref{defn:basefortopology}\\
\hline
$\ov n$&the set $\{1,2,\dots,n\}$ with $n\in\N$& Notation \ref{not:nset}&\pageref{not:nset}\\
\hline
$\lVert\;\cdot\;\rVert_{\a}$&seminorm indexed by $\a$ in some set&Proposition \ref{prop:seminormsondualfromoriginal}&\pageref{prop:seminormsondualfromoriginal}\\
\hline
$V^{\vee}$&algebraic (linear) dual of vector space $V$& Proposition \ref{prop:topologicaldual}&\pageref{prop:topologicaldual}\\
\hline
$V^{*}$&topological dual of normed space $V$&Proposition \ref{prop:topologicaldual}&\pageref{prop:topologicaldual}\\
\hline
$\Theta$&ordered index set for nets&Proposition \ref{prop:weak*topologyconvergence}&\pageref{prop:weak*topologyconvergence}\\
\hline
$\s(\mA)$&spectrum of commutative Banach algebra $\mA$&Definition \ref{def:spectrum}&\pageref{def:spectrum}\\
\hline
$\s(f)\equiv\s^{f}$&spectrum of a homomorphism $f$&Proposition \ref{prop:spectrumofmap}&\pageref{prop:spectrumofmap}\\
\hline
$\s$&spectrum functor $\cCAlg^{\op}\to\cH$&Proposition \ref{prop:cCAlgtocHaus}&\pageref{prop:cCAlgtocHaus}\\
\hline
$\G$&Gelfand transform $\G:\id_{\cCAlg}\Rightarrow C\circ\s$&Proposition \ref{prop:Gelfandtransform}&\pageref{prop:Gelfandtransform}\\
\hline
$h$&natural transformation $h:\id_{\cH}\Rightarrow\s\circ C$&Proposition \ref{prop:hnaturaltransform}&\pageref{prop:hnaturaltransform}\\
\hline
$(Y,\mathcal{M}_{Y})$&\begin{tabular}{c}compact Hausdorff space $Y$\\with its Borel $\s$-algebra $\mathcal{M}_{Y}$\end{tabular}&Section \ref{sec:chausmeasures}&\pageref{sec:chausmeasures}\\
\hline
$\varnothing$&the empty set&Definition \ref{defn:complexmeasure}&\pageref{defn:complexmeasure}\\
\hline
$\mu,\nu$&measure (positive, real, or complex)&Definition \ref{defn:complexmeasure}&\pageref{defn:complexmeasure}\\
\hline
$\C\mathrm{Meas}(Y)$&set of complex measures on $Y$&Definition \ref{defn:complexmeasure}&\pageref{defn:complexmeasure}\\
\hline
$\chi_{E}$&characteristic/indicator function on $E$&Definition \ref{defn:integral}&\pageref{defn:integral}\\
\hline
$s$&simple function&Definition \ref{defn:integral}&\pageref{defn:integral}\\
\hline
$\int\;\cdot\;d\mu$&integral with respect to a measure $\mu$&Definition \ref{defn:integral}&\pageref{defn:integral}\\
\hline
$\mathrm{ProbMeas}(Y)$&\begin{tabular}{c}set of regular probability measures on $Y$\\(equipped with the vague topology)\end{tabular}&Notation \ref{not:ProbMeas}&\pageref{not:ProbMeas}\\
\hline
$\mu:X\stoch Y$&\begin{tabular}{c}stochastic map described by a\\continuous function $\mu:X\to\mathrm{ProbMeas}(Y)$\end{tabular}&Definition \ref{defn:continuousstochasticprocess}&\pageref{defn:continuousstochasticprocess}\\
\hline
$\de_{f}$&\begin{tabular}{c}Dirac delta stochastic map\\associated to a continuous function $f$\end{tabular}&Example \ref{ex:diracstochastic}&\pageref{ex:diracstochastic}\\
\hline
$\nu\circ\mu$&\begin{tabular}{c}composition of $\mu:X\stoch Y$\\followed by $\nu:Y\stoch Z$\end{tabular}&Proposition \ref{prop:compositionstochastic}&\pageref{prop:compositionstochastic}\\
\hline
$\cHStoch$&\begin{tabular}{c}category of compact Hausdorff spaces\\and stochastic maps\end{tabular}&Theorem \ref{thm:cHStoch}&\pageref{thm:cHStoch}\\
\hline
$\de$&Dirac measure functor $\de:\cH\to\cHStoch$&Theorem \ref{thm:cHStoch}&\pageref{thm:cHStoch}\\
\hline
$C$&functor $\cHStoch^{\op}\to\cCAlgPos$&Theorem \ref{thm:stochasticcontinuousfunctor}&\pageref{thm:stochasticcontinuousfunctor}\\
\hline
RMK&short for ``Riesz-Markov-Kakutani''&Theorem \ref{thm:Riesz}&\pageref{thm:Riesz}\\
\hline
$\mathrm{ex}(C)$&extreme points of a convex set $C$&Definition \ref{defn:convex}&\pageref{defn:convex}\\
\hline
$\mS(\mA)$&states on commutative $C^*$-algebra $\mA$&Example \ref{ex:statesconvexspectrumextreme}&\pageref{ex:statesconvexspectrumextreme}\\
\hline
$\mathrm{ch}(X)$&convex hull of $X$&Definition \ref{defn:convexhull}&\pageref{defn:convexhull}\\
\hline
$\s$&\begin{tabular}{c}stochastic spectrum functor\\$\cCAlgPos^{\op}\to\cHStoch$\end{tabular}&Theorem \ref{thm:stochasticspectrum}&\pageref{thm:stochasticspectrum}\\
\hline
\end{longtable}

\pagebreak
\bibliographystyle{plain}
\bibliography{states}

{\small
\textsc{Mathematics Department, University of Connecticut, 
341 Mansfield Road U1009, 
Storrs, CT 06269, USA}

{{\em Email:} {\tt\href{mailto:arthur.parzygnat@uconn.edu}
{arthur.parzygnat@uconn.edu}}}}
\end{document}